\documentclass[11pt]{amsart}

\usepackage{amssymb}
\usepackage{mathtools}

\usepackage{color}

\newtheorem{theo}{Theorem} 

\newtheorem{lemma}{Lemma}[section]
\newtheorem{prop}[lemma]{Proposition}
\newtheorem{corol}[lemma]{Corollary}
\newtheorem{theoin}[lemma]{Theorem}
\newtheorem{claim}[lemma]{Claim}

\theoremstyle{remark}
\newtheorem{remark}[lemma]{Remark}

\newtheorem{notation}[lemma]{Notation}

\theoremstyle{definition}
\newtheorem{defi}[lemma]{Definition}

\newcommand{\coltwo}[2]{\begin{pmatrix}#1\\#2\end{pmatrix}}

\newcommand{\lin}{\textsc{l}}

\newcommand{\dd}{\mathsf{d}}

\newcommand{\NN}{\mathbb{N}}
\newcommand{\RR}{\mathbb{R}}

\newcommand{\eps}{\varepsilon}

\newcommand{\AAA}{\mathcal{A}}

\newcommand{\CCC}{\mathcal{C}}

\newcommand{\JJJ}{\mathcal{J}}

\newcommand{\OOO}{\mathcal{O}}
\newcommand{\PPP}{\mathcal{P}}

\newcommand{\lf}{\left}
\newcommand{\rg}{\right}

\newcommand{\tlambda}{\widetilde{\lambda}}

\newcommand{\tw}{\widetilde{w}}

\newcommand{\ttau}{\widetilde{\tau}}

\newcommand{\trho}{\tilde{\rho}}

\newcommand{\tbeta}{\widetilde{\beta}}

\newcommand{\tU}{\widetilde{U}}

\newcommand{\tu}{\tilde{u}}

\newcommand{\tv}{\tilde{v}}

\newcommand{\tx}{\widetilde{x}}
\newcommand{\Stt}{\texttt{S}}

\newcommand{\rhob}{\overline{\rho}}

\newcommand{\tK}{\widetilde{K}}

\newcommand{\loc}{\rm loc}
\newcommand{\hdot}{\dot{H}^1}
\newcommand{\Hdot}{\dot{H}^1(\RR^N)}

\newcommand{\EMPH}[1]{\medskip\noindent\textit{#1}.}

\newcommand{\ssstyle}{\scriptscriptstyle}

\DeclareMathOperator{\supp}{supp}

\DeclareMathOperator{\ext}{ext}
\DeclareMathOperator{\Div}{div}

\newcommand{\ds}{\displaystyle}

\numberwithin{equation}{section} 

\title[Blow-up for energy critical wave]{Universality of blow-up profile for small radial type II blow-up solutions of energy-critical wave equation}
\author[T.~Duyckaerts]{Thomas Duyckaerts$^1$}
\author[C.~Kenig]{Carlos Kenig$^2$}
\author[F.~Merle]{Frank Merle$^3$}
\thanks{$^1$Cergy-Pontoise (UMR 8088). Partially supported by ANR Grants ONDNONLIN and ControlFlux}
\thanks{$^2$University of Chicago. Partially supported by NSF Grant DMS-0456583}
\thanks{$^3$Cergy-Pontoise (UMR 8088), IHES, CNRS. Partially supported by ANR Grant ONDNONLIN}

\date{\today}

\begin{document}
\begin{abstract}
 Consider the energy critical focusing wave equation on the Euclidian space. A blow-up type II solution of this equation is a solution which has finite time of existence but stays bounded in the energy space. The aim of this work is to exhibit universal properties of such solutions. 

Let $W$ be the unique radial positive stationary solution of the equation. Our main result is that in dimension $3$, under an appropriate smallness assumption, any type II blow-up radial solution is essentially the sum of a rescaled $W$ concentrating at the origin and a small remainder which is continuous with respect to the time variable in the energy space. This is coherent with the solutions constructed by Krieger, Schlag and Tataru. One ingredient of our proof is that the unique radial solution which is compact up to scaling is equal to $W$ up to symmetries.
\end{abstract}

\maketitle

\tableofcontents

\section{Introduction}
Consider the focusing energy-critical wave equation on an interval $I$ ($0\in I$)
\begin{equation}
\label{CP}
\left\{ 
\begin{gathered}
\partial_t^2 u -\Delta u-|u|^{\frac{4}{N-2}}u=0,\quad (t,x)\in I\times \RR^N\\
u_{\restriction t=0}=u_0\in \hdot,\quad \partial_t u_{\restriction t=0}=u_1\in L^2,
\end{gathered}\right.
\end{equation}
where $u$ is real-valued, $N\in \{3,4,5\}$, and $\hdot:=\Hdot$. We will denote by $S(I):=L^{\frac{2(N+1)}{N-2}}(I\times \RR^N)$. We will often restrict ourselves to the case of radial solutions in space dimension $N=3$.

The Cauchy problem \eqref{CP} is locally well-posed in $\hdot\times L^2$. This space is invariant by the scaling of the equation: if $u$ is a solution to \eqref{CP}, $\lambda>0$ and
$$ u_{\lambda}=\frac{1}{\lambda^{\frac{N-2}{2}}}u\left(\frac{t}{\lambda},\frac{x}{\lambda}\right),$$
then $u_{\lambda}$ is also a solution and $\|u_{\lambda}(0)\|_{\hdot}=\|u_0\|_{\hdot}$, $\|\partial_t u_{\lambda}(0)\|_{L^2}=\|u_1\|_{L^2}$.

Let $T_+\in(0,+\infty]$ be the maximal positive time of definition for the solution $u$. It satisfies the following finite time blow-up criterion
\begin{equation}
\label{FBUC}
T_+<\infty\Longrightarrow \|u\|_{S(0,T_+)}=+\infty.
\end{equation} 
Note that this criterion does not rule out type II blow-up, i.e. solutions such that $T_+<\infty$ and
\begin{equation}
\label{BlowUpII}
\sup_{t\in [0,T_+)} \|\partial_t u(t)\|^2_{L^2}+\|\nabla u(t)\|_{L^2}^2<\infty.
\end{equation}
This is different to the case of lower order non-linearity (of the form $|u|^{p-1}u$ with $p<\frac{N+2}{N-2}$), where the finite time blow-up implies the blow-up of the energy norm. 

Energy arguments of Levine type \cite{Levine74} are not expected to give
directly type II blow-up. Examples of radial type II blow-up solutions of \eqref{CP} were constructed in space dimension $N=3$ by Krieger, Schlag and Tataru \cite{KrScTa09}. 
The aim of this article is to exhibit universal properties of this type of solutions.

Let
\begin{equation}
\label{defW}
W:=\frac{1}{\left(1+\frac{|x|^2}{N(N-2)}\right)^{\frac{N-2}{2}}},
\end{equation}
which is a stationary solution of \eqref{CP}.  
The construction of \cite{KrScTa09} relies on an elaborate fixed point argument which yields the following description of the solution
\begin{equation}
\label{sol_KST}
u(t)=\frac{1}{\lambda^{1/2}(t)}W\left(\frac{x}{\lambda(t)}\right)+\varepsilon(t),
\end{equation}
where $\lambda(t)=(T_+-t)^{1+\nu}$, $\nu>0$  and
\begin{equation}
\label{cond_KST}
\lim_{t \to  T_+}\int_{|x|\leq T_+-t} |\nabla \eps(t)|^2\,dx+\int_{|x|\leq T_+-t} |\partial_t  \eps(t)|^2\,dx+\int_{|x|\leq T_+-t} |\eps(t)|^6\,dx =0.
\end{equation} 
In this work, we investigate the converse problem: if we consider an arbitrary type II radial blow-up solution, does such a decomposition hold?

We will obtain this result in an appropriate smallness case.
From \cite{KeMe08}, if,
$$ \sup_{t\in [0,T_+)}\|\nabla u(t)\|_{L^2}^2+\|\partial_t u(t)\|_{L^2}^2< \|\nabla W\|_{L^2}^2,$$
then $T_+=+\infty$ and the solution scatters forward in time, and in particular does not blow up. 
In this work, the authors introduce a general road map to tackle such critical problems in focusing and defocusing situations. From a concentration-compactness result (in this case \cite{BaGe99}), one reduces the proof to some rigidity property of solutions of \eqref{CP} that are compact in the energy space up to the invariances of the equation. This rigidity property has to be shown by independent arguments.

The threshold $\|\nabla W\|^2_{L^2}$ is sharp. Indeed from \cite{KrScTa09}, for all $\eta_0>0$ there exists a type II blow-up solution such that
\begin{equation}
\label{bound}
 \sup_{t\in [0,T_+)} \|\nabla u(t)\|_{L^2}^2+\|\partial_t u(t)\|_{L^2}^2\leq \|\nabla W\|_{L^2}^2+\eta_0.
\end{equation} 
In the present article, we consider type II blow-up solutions such that \eqref{bound} holds. Our main result is the following.
\begin{theo}
\label{T:classification}
Assume that $N=3$. There exists $\eta_0>0$ such that for any \emph{radial} solution $u$ of \eqref{CP} such that $T_+(u)=T_+<\infty$ that satifies \eqref{bound},
there exist a solution $v(t)$ of \eqref{CP} defined in a neighborhood of $t=T_+$, a sign $\iota_0\in \{\pm 1\}$, and a $C^0$ positive function $\lambda(t)$ on $(0,T_+)$ such that, as $t\overset{\scriptscriptstyle{<}}{\to} T_+$,
\begin{gather}
\label{dev_u}
u(t)=v(t)+\frac{\iota_0}{\lambda^{1/2}(t)}W\left(\frac{x}{\lambda(t)}\right)+o(1)\text{ in }\hdot,\\
\label{dev_dtu}
\partial_t u(t)=\partial_t v(t)+o(1)\text{ in }L^2,\\
\label{lambda_small}
\lambda(t)=o(T_+-t).
\end{gather}
\end{theo}
Note that \eqref{dev_u}, \eqref{dev_dtu} imply that $u$ is of the form \eqref{sol_KST} with $\eps$ satisfying \eqref{cond_KST}: any radial blow-up solution satisfying \eqref{bound} is of the type of the solutions constructed by Krieger, Schlag and Tataru.

As it is now well-known from previous works on similar problems (see remarks below), the result is based on classification of solutions of \eqref{CP} that are compact up to the symmetry of the equation. We state this result for its own interest. 
\begin{theo}
\label{T:compact=W}
Let $u$ be a nonzero radial solution of \eqref{CP} in space dimension $N=\{3,4,5\}$. Assume that there exists a function $\lambda(t)$ of $t\in (T_-(u),T_+(u))$ such that
$$ K=\left\{\left(\lambda^{N/2-1}(t)u\left(t,\lambda(t)\cdot\right),\lambda^{N/2}(t)\partial_t u\left(t,\lambda(t)\cdot\right)\right),\; t\in \RR\right\}$$
has compact closure in $\hdot\times L^2$. Then there exist $\lambda_0>0$ and a sign $\iota_0\in \{\pm 1\}$ such that
$$u(t,x)=\frac{\iota_0}{\lambda_0^{N/2-1}}W\left(\frac{x}{\lambda_0}\right).$$
\end{theo}

\begin{remark}
 The proof of Theorem \ref{T:compact=W} (see Section \ref{S:radial_compact}) uses the material of \cite{DuMe08}, where a first classification result of this type was obtained. Namely, at the energy threshold $E(u_0,u_1)=E(W,0)$, all solutions such that $\int |\nabla u_0|^2+\int |u_1|^2\leq \int |\nabla W|^2$ are globally defined, and the only ones that do not scatter are (up to the transformations of the equation) $W$ and a solution $W^-$, which scatters backward in time and tends to $W$ exponentially as $t$ goes to $+\infty$.
\end{remark}

\begin{remark}
These results are essential to understand type II blow-up. After one has exhibited an universal profile for blow-up, one can hope using local dynamics near $W$, or linearization around $W$ (see e.g. \cite{KrSc07}) to understand the possible blow-up speeds, which will complete the program to understand type II blow-up. Moreover, this is the first step to prove that the boundary of the set of initial data that lead to blow-up is given by type II blow-up solutions.

The proof of Theorem \ref{T:classification} highlights, through the mechanism of profile decomposition and the finite speed of propagation, why the only candidates to be type II blow-up profiles are compact solutions. The only case where such a striking fact was established was for GKdV by Martel and Merle \cite{MaMe00,MaMe01}.
\end{remark}

\begin{remark}
In the case of nonlinear wave maps, all blow-up solutions are of type II: the equation is defocusing, in the sense that the energy provides a bound on the energy norm. An analogue of Theorem \ref{T:classification} is known locally in space for a sequence of times (without size condition due to the defocusing nature of the equation). Namely, if $u$ is a blow-up solution, there exist a sequence of times $t_n\to T_+$ and a sequence $\lambda_n \to 0^+$ such that $u(t_n,x/\lambda_n)$ tends to a nonlinear object. This follows from a remarkable paper of Christodoulou  and Tahvildar-Zadeh \cite{ChTZ93}. See also the article of Shatah and Tahvildar-Zadeh \cite{ShTZ97} which established a result similar to Theorem \ref{T:compact=W} in this context, the articles of Struwe \cite{Struwe02,Struwe03}, and the recent preprint of Sterbenz and Tataru \cite{StTa09P} for the general case of solutions without any special invariant properties. We also refer to the works of Rodnianski and Sterbenz \cite{RoSt06P}, Krieger, Schlag and Tataru \cite{KrScTa08} and Rapha\"el and Rodnianski \cite{RaRo09P} for the construction of blow-up solutions.
\end{remark}
\begin{remark}
 Universality of blow-up profiles for a critical equation, as $t$ goes to the blow-up time $T_+$ (without restriction to a sequence of times) was established, also under smallness condition, in two cases:
\begin{itemize}
\item[-] for the critical KdV equation
$$u_t=(u_{xx}+u^5)_x, \quad x\in \mathbb{R}$$
by classification of compact solutions of the GKdV equation: see Martel Merle \cite{MaMe00,MaMe01,MaMe03};
\item[-] for the mass-critical NLS equation 
$$ iu_t=\Delta u+|u|^{\frac{4}{N}}u,\quad 1\leq N\leq 5$$
by Merle and Rapha\"el, by classification of solutions  that are nondispersive (in a weak sense) \cite{MeRa04}.
\end{itemize}
See also the subcritical wave equation in dimension one where all blow-up profiles were found by Merle and Zaag \cite{MeZa07,MeZa08P} for general data (see the work of Caffarelli and Friedman for specific data \cite{CaFr86}).
\end{remark}

We next sketch the proof of Theorem \ref{T:classification}.

Consider a radial, blow-up solution $u$ of \eqref{CP} in space dimension $N=3$ that satisfies \eqref{bound} and assume for simplicity that the blow-up time $T_+(u)$ is $1$. As is shown in Section \ref{S:reg_sing}, $u$ may be decomposed as the sum of a solution to \eqref{CP} which is well-defined around the blow-up time, and a singular part $a(t,x)$ which is supported in the light cone $|x|\leq 1-t$. Consider a sequence $t_n\to 1^-$ and a Bahouri-G\'erard \cite{BaGe99} profile decomposition associated to the sequence $(a(t_n),\partial_t a(t_n))_n$. According to the result of \cite{KeMe08}, the bound \eqref{bound} implies that there is one large profile and that the other profiles are small (see Remark \ref{R:bound_dN}). This contrasts with \cite{KeMe08} where the minimality of the solution imposes automatically that there is only one profile (the large one). In our case, we must show by another mechanism that the small profiles do not exist, which would imply by Theorem \ref{T:compact=W} that the large profile is $W$, yielding Theorem \ref{T:classification}.

The main idea to exclude the small profiles is that any small block of energy norm decoupled from the main profile would yield, for each time $t_n$, a non negligible amount of energy norm localized on a light cone.\footnote{This follows from a property of the radial three-dimensional linear wave equation that does not hold in the non-radial setting or in higher dimensions (see Lemma \ref{L:linear_behavior}), which is the main reason why we restricted ourselves to radial solutions in space dimension $N=3$.} By finite speed of propagation one can show that these small energy blocks, localized in disjoint light cones, sum up, implying the blow-up of the energy norm, which contradicts the bound \eqref{bound}.
A similar phenomenon is highlighted in the context of the nonlinear Schr\"odinger equation in \cite{MeRa05,MeRa08}. Unfortunately, this strategy can be implemented only for a class of profiles that are very small and exterior in a certain sense (see Proposition \ref{P:ext_profile}), and we must exclude the other small profiles by indirect means (see Sections \ref{S:non_self_similar}, \ref{S:global} and \ref{S:radial}).

By Proposition \ref{P:sum_of_W}, the existence of a sequence $\tau_n\to 1^-$ such that $a(\tau_n)$ concentrates at a speed faster than self-similar implies, at least for another sequence of times $t_n\to 1^-$, that all profiles are equal to the stationary solution $W$. This property follows from rigidity arguments involving virial type identities. In particular, there cannot be small profiles, and the bound \eqref{bound} implies that $W$ is the only profile for this particular sequence $t_n$. This yields the strong condition that the energy of the singular part $E(a,\partial_ta)$ tends to $E(W,0)$ as $t\to 1^-$ (see Corollary \ref{C:energy}) which can be combined with the results of \cite{KeMe08} to complete the proof (see \S \ref{SS:proof_compact} and \S \ref{SS:proof_W}). 

It remains to exclude the case of self-similar concentration, which is the object of Proposition \ref{P:noss}. For this, we argue by contradiction, showing (as a consequence of the non-existence of small exterior profiles) that this self-similar concentration, if it exists, must concern the large profile. The solution of \eqref{CP} corresponding to this profile is globally defined and non-scattering backward in time, satisfies a global bound similar to \eqref{bound} for negative times, and is partially located around the light cone $|t|=|x|$, as $t\to -\infty$. This type of solution is excluded by Proposition \ref{P:norm_infinity}, using the non-existence, shown in \cite{KeMe08}, of self-similar blow-up solutions of \eqref{CP} which are compact up to scaling.

\smallskip

The outline of the paper is as follows.

After some preliminaries (Section \ref{S:preliminaries}), we give in Section \ref{S:reg_sing} general results on type II blow-up solutions of \eqref{CP} in space dimensions $N=3,4,5$. 
In the two next sections we restrict ourselves to radial solutions in space dimension $3$.
In Section \ref{S:general_radial} we show the nonexistence of small exterior profiles for a radial type II blow-up solution.  In Section \ref{S:non_self_similar}, we assume that the solution does not concentrate at a self-similar rate, and show that in this case, there exists a sequence $t_n\to T_+$ such that $u(t_n)$ decomposes as a sum of rescaled stationary solutions. 

In the two following sections, we consider solutions of \eqref{CP} that do not blow up in finite time.  Section \ref{S:radial_compact} is devoted to the proof of Theorem \ref{T:compact=W}, which is a consequence of the classification result of \cite{DuMe08}. Section \ref{S:global} is concerned with the localization of the energy for globally defined, bounded, non-scattering solutions of \eqref{CP}. 

Section \ref{S:radial} gathers the results of all previous Sections to prove Theorem \ref{T:classification}. In Appendix \ref{A:ortho} we prove some technical properties of profile decompositions. Appendix \ref{A:faraway} shows a simple result on a family of sequences of positive numbers which is needed in some parts of the proof.

In all the article, for sequences of positive numbers $\{\alpha_n\}_n$ and $\{\beta_n\}_n$, we will write $\alpha_n\ll \beta_n$ when $\alpha_n/\beta_n\to 0$ as $n\to \infty$, and $\alpha_n\approx \beta_n$ when $C^{-1}\alpha_n\leq \beta_n\leq C\alpha_n$ for some large constant $n$. We will denote by $o_n(1)$ a sequence that goes to $0$ as $n$ goes to $\infty$.

\section{Preliminaries}
\label{S:preliminaries}

\subsection{Cauchy problem}
The Cauchy problem for equation \eqref{CP} was developped in \cite{Pecher84,GiSoVe92,LiSo95,ShSt94,ShSt98,Sogge95,Kapitanski94}.
If $I$ is an interval, we denote by
$$  S(I)=L^{\frac{2(N+1)}{N-2}}\left(I\times\RR^N\right), \;W(I)=L^{\frac{2(N+1)}{N-1}}\left(I\times\RR^N\right),\;N(I)=L^{\frac{2(N+1)}{N+3}}\left(I\times\RR^N\right).$$
Let $\Stt_{\lin}(t)$ be the one-parameter group associated to the linear wave equation. By definition, if $(v_0,v_1)\in \hdot\times L^2$ and $t\in \RR$, $v(t)=\Stt_{\lin}(t)(v_0,v_1)$ is the solution of
\begin{gather}
 \label{lin_wave}
\partial_t^2 v-\Delta v=0,\\
\label{lin_wave_IC}
v_{\restriction t=0}=v_0,\quad
\partial_t v_{\restriction t=0}=v_1.
\end{gather}
We have
$$ \Stt_{\lin}(t)(v_0,v_1)=\cos(t\sqrt{-\Delta})v_0+\frac{1}{\sqrt{-\Delta}}\sin(t\sqrt{-\Delta})v_1.$$
By Strichartz and Sobolev estimates,
\begin{equation}
\label{strichartz}
\lf\|v\rg\|_{S(\RR)}+\lf\|D_x^{1/2}v\rg\|_{W(\RR)}\leq C_S\left(\|v_0\|_{\hdot}+\|v_1\|_{L^2}\right).
\end{equation} 
A solution of \eqref{CP} on an interval $I$, where $0\in I$, is a function $u\in C^0(I,\hdot)$ such that  $\partial_t u\in C^0(I,L^2)$, 
\begin{equation}
\label{finite_norms}
J\Subset I\Longrightarrow \|D_x^{1/2}u\|_{W(J)}+\|u\|_{S(J)}<\infty
\end{equation}
satisfying the Duhamel formulation
\begin{equation}
\label{solution}
u(t)=\Stt(t)u_0+\int_0^t \frac{\sin\big((t-s)\sqrt{-\Delta}\big)}{\sqrt{-\Delta}}\left(|u(s)|^{\frac{4}{N-2}}u(s)\right)ds.
\end{equation} 
We recall there exists a small $\delta_0>0$ such that for any interval $I$ containing $0$ and any $(u_0,u_1)\in \hdot\times L^2$ such that 
\begin{equation}
\label{small_data}
\left\|\Stt_{\lin}(t)(u_0,u_1)\right\|_{S(I)}<\delta_0, 
\end{equation} 
there exists an unique solution $u$ of \eqref{CP} on $I$. Furthermore if $\delta_0$ is chosen small enough, this solution satisfies:
\begin{equation}
\label{small_data_bis}
\|u\|_{S(I))}\leq 2\left\|\Stt_{\lin}(t)(u_0,u_1)\right\|_{S(I)}.
\end{equation}
Sticking together these local solutions, we get that for any initial condition $(u_0,u_1)$ in the energy space, there exists an unique solution $u$ of \eqref{CP}, which is defined on a maximal interval of definition $$I_{\max}=I_{\max}(u_0,u_1)=\big(T_-(u_0,u_1),T_+(u_0,u_1)\big).$$
We will often write $I_{\max}(u)$, $T_{\pm}(u)$, instead of $I_{\max}(u_0,u_1)$, $T_{\pm}(u_0,u_1)$. 

If $\left\|\Stt_{\lin}(t)(u_0,u_1)\right\|_{S(I)}=\delta<\delta_0$, then $u$ is close to the linear solution with initial condition $(u_0,u_1)$ in the following sense: if $A=\left\|D_x^{1/2}\Stt_{\lin}(t)(u_0,u_1)\right\|_{W(I)}$, we have  
\begin{multline}
\label{almost_linear}
\left\|u(\cdot)-\Stt_{\lin}(\cdot)(u_0,u_1)\right\|_{S(I)}\\+\sup_{t\in I}\left(\left\|u(t)-\Stt_{\lin}(t)(u_0,u_1)\right\|_{\hdot}+\left\|\partial_tu(t)-\partial_t(\Stt_{\lin}(t)(u_0,u_1))\right\|_{L^2}\right)
\leq CA\delta^{\frac{4}{N-2}},
\end{multline} 
(see for example \cite{KeMe06}, proof of Theorem 2.7).

Any solution $u$ of \eqref{CP} satisfies the blow-up criterion \eqref{FBUC}, and the analogue for negative time. As a consequence, if 
$\|u\|_{S(0,T_+)}<\infty$, then $T_+=+\infty$. Furthermore in this case, the solution scatters forward in time in $\hdot\times L^2$: there exists a solution $v$ of the linear equation \eqref{lin_wave} such that
$$\lim_{t\to +\infty}\|u(t)-v(t)\|_{L^2}+\|\partial_t u(t)-\partial_t  v(t)\|_{L^2}=0.$$
 Of course an analoguous statement holds backward in time also.

We next recall a long-time perturbation theory result for \eqref{CP} (see Theorem  2.20 of \cite{KeMe08}).
\begin{theoin}
\label{T:LTPT}
Let $M>0$. There exists $\eps_0=\eps_0(M)$ with the following property. Let $I\subset \RR$ be a time interval such that $0\in I$, and $\tilde{u}$ be defined on $I\times \RR^N$ such that
$$ \|\tu\|_{S(I)}+\sup_{t\in I} \lf(\left\|\tu(t)\rg\|_{\hdot}+\lf\|\partial_t\tu(t)\right\|_{L^2}\rg)\leq M,\quad J\Subset I\Rightarrow \left\|D_x^{1/2}\tu\right\|_{W(J)}<\infty,$$
Denote by $(\tu_0,\tu_1)=(\tu(0),\partial_t \tu(0))$. Consider $(u_0,u_1)\in \hdot\times L^2$ and  $\eps\in(0,\eps_0)$. Assume 
\begin{gather*}
\partial_t^2 \tu-\Delta \tu-|\tu|^{\frac{4}{N-2}}\tu=e,\quad (t,x)\in I\times \RR^N\\
\text{and }
\left\|u_0-\tu_0\rg\|_{\hdot}+\lf\|u_1-\tu_1\right\|_{L^2}+\left\|D_x^{1/2} e\right\|_{N(I)}\leq \eps.
\end{gather*} 
Then the solution $u$ of \eqref{CP} with initial condition $(u_0,u_1)$ satisfies $I_{\max}(u)\subset I$ and for a $\beta_0>0$,
$$ \|u\|_{S(I)}\leq C(M),\quad \lf(\sup_{t\in I}\left\|u(t)-\tu(t)\rg\|_{\hdot}+\lf\|\partial_t u(t)-\partial_t\tu(t)\right\|_{L^2}\rg)\leq C(M)\eps^{\beta_0}.$$
\end{theoin}
\subsection{Remarks on stationary solutions of \eqref{CP}}
Recall from  \eqref{defW} the definition of the stationary solution $W$. It is known from the works of T.~Aubin \cite{Au76} and G.~Talenti \cite{Ta76} that $W$ is the unique minimizer, up to translation, scaling and multiplication by a scalar constant, for the Sobolev inequality on $\RR^N$
\begin{equation*}
\|f\|_{L^{\frac{2N}{N-2}}}\leq C_N\|\nabla f\|_{L^2}.
\end{equation*} 
By a classical ODE argument, we also have the following uniqueness result:
\begin{claim}
 \label{C:uniqueness_W}
Let $U$ be a (real) $\hdot(\RR^N)$ radial solution of 
$$\Delta U+|U|^{\frac{4}{N-2}}U=0.$$
Then 
$$U=0 \quad \text{or}\quad \exists \lambda_0>0,\;U=\pm\frac{ 1}{\lambda_0^{\frac{N-2}{2}}}W\left(\frac{x}{\lambda_0}\right).$$
\end{claim}

Note that the equation $\Delta W+W^{\frac{N+2}{N-2}}=0$ implies $E(W,0)=\frac{1}{N}\int |\nabla W|^2>0$. This fact is used to prove the following variational properties of $W$ which will be needed throughout the paper.
\begin{claim}
\label{C:variationnal}
Let $v\in \hdot$. Then
\begin{equation}
 \label{variationnal}
\|\nabla v\|_{L^2}^2\leq \|\nabla W\|_{L^2}^2\text{ and }E(v,0)\leq E(W,0)\Longrightarrow \|\nabla v\|_{L^2}^2\leq \frac{\|\nabla W\|_{L^2}^2}{E(W,0)}E(v,0)=N E(v,0).
\end{equation}
Furthermore, if $\|\nabla v\|_{L^2}^2\leq \left(\frac{N}{N-2}\right)^{\frac{N-2}{2}}\|\nabla W\|_{L^2}^2$, then $E(v,0)\geq 0$.
\end{claim}
\begin{proof}
The first part of the claim is shown in \cite{DuMe08} (see the proof of Claim 2.4). For the second part, write
$$ E(v,0)=\frac{1}{2}\int |\nabla v|^2-\frac{N-2}{2N}\int |v|^{\frac{2N}{N-2}}\geq \frac{1}{2}\int |\nabla v|^2-\frac{N-2}{2N}C_N^{\frac{2N}{N-2}} \left(\int |\nabla v|^2\right)^{\frac{N}{N-2}},$$
where $C_N$ is the best constant in the Sobolev inequality $\|v\|_{\frac{2N}{N-2}}\leq C_N\|\nabla v\|_{L^2}$. Let $y=\int|\nabla v|^2$ and assume that $E(v,0)$ is negative. Then
$$0> \frac{1}{2} y-\frac{N-2}{2N} C_N^{\frac{2N}{N-2}}y^{\frac{N}{N-2}}.$$
This shows that $y\geq y_*$, where $y_*$ is the unique positive solution of $\frac{1}{2} y-\frac{N-2}{2N} C_N^{\frac{2N}{N-2}}y^{\frac{N}{N-2}}=0$. Using that $C_N^{-N}=\int |\nabla W|^2$, we obtain $y*=\left(\frac{N}{N-2}\right)^{\frac{N-2}{2}}\int|\nabla W|^2$, which concludes the proof.
\end{proof}

\subsection{Profile decomposition}
We recall here the profile decomposition of H.~Bahouri and P.~G\'erard \cite{BaGe99}. This paper is written in space dimension $N=3$ but the results stated below hold in all dimension $N\geq 3$. See also \cite{BrCo85} and \cite{Li85Reb} for the elliptic case and \cite{MeVe98} for the Schr\"odinger equation.

Consider a sequence $(v_{0,n},v_{1,n})_n$ which is bounded in $\hdot\times L^2$. 
Let $(U^j_{\lin})_{j\geq 0}$ be a sequence of solutions of the linear equation \eqref{lin_wave},
with initial data $(U^j_0,U^j_1)\in \hdot\times L^2$, and $(\lambda_{j,n};x_{j,n};t_{j,n})\in (0,+\infty)\times \RR^N\times \RR$, $j,n\in \NN$, be a family of parameters satisfying the pseudo-orthogonality relation
\begin{equation}
\label{ortho_param}
j\neq k\Longrightarrow \lim_{n\to \infty} \frac{\lambda_{j,n}}{\lambda_{k,n}}+\frac{\lambda_{k,n}}{\lambda_{j,n}}+\frac{|t_{j,n}-t_{k,n}|}{\lambda_{j,n}}+\frac{\left|x_{j,n}-x_{k,n}\right|}{\lambda_{j,n}}=+\infty.
\end{equation} 
We say that $(v_{0,n},v_{1,n})_n$ admits a profile decomposition  $\lf\{U_{\lin}^j\rg\}_j$, $\lf\{\lambda_{j,n},x_{j,n};t_{j,n}\rg\}_{j,n}$ when 
\begin{equation}
\label{decompo_profil}
\left\{\begin{aligned}
 v_{0,n}&=\sum_{j=1}^J \frac{1}{\lambda_{j,n}^{\frac{N-2}{2}}}U_{\lin}^j\left(\frac{-t_{j,n}}{\lambda_{j,n}},\frac{x-x_{j,n}}{\lambda_{j,n}}\right)+w_{0,n}^J(x),\\
v_{1,n}&=\sum_{j=1}^J \frac{1}{\lambda_{j,n}^{\frac{N}{2}}}\partial_t U_{\lin}^j\left(\frac{-t_{j,n}}{\lambda_{j,n}},\frac{x-x_{j,n}}{\lambda_{j,n}}\right)+w_{1,n}^J(x),
\end{aligned}\right.
\end{equation}
with
\begin{equation}
\label{small_w}
\lim_{n\rightarrow+\infty}\limsup_{J\rightarrow+\infty} \left\|w_{n}^J\right\|_{S(\RR)}=0,
\end{equation}
where $w_{n}^J$ is the solution of \eqref{lin_wave} with initial conditions $(w_{0,n}^J,w_{1,n}^J)$. 
Then:
\begin{theoin}[\cite{BaGe99}]
\label{T:profile}
If the sequence $(v_{0,n},v_{1,n})_n$ is bounded in the energy space $\hdot\times L^2$, there always exists a subsequence of $(v_{0,n},v_{1,n})_n$ which admits a profile decomposition. Furthermore, 
\begin{equation}
\label{weak_CV_wJ}
j\leq J\Longrightarrow 
\left(\lambda_{j,n}^{\frac{N-2}{2}} w_n^J\left(t_{j,n},x_{j,n}+\lambda_{j,n}y\right),\lambda_{j,n}^{\frac{N}{2}} \partial_tw_n^J\left(t_{j,n},x_{j,n}+\lambda_{j,n}y\right)\right)\xrightharpoonup[n\to \infty]{}0,
\end{equation}
weakly in $\hdot_y\times L^2_y$, and the following Pythagorean expansions hold for all $J\geq 1$
\begin{gather}
\label{pythagore1a} 
\left\|v_{0,n}\rg\|_{\hdot}^2=\sum_{j=1}^J \left\|U^j_{\lin}\left(\frac{-t_{j,n}} {\lambda_{j,n}}\right)\rg\|_{\hdot}^2+\left\|w_{0,n}^J\rg\|_{\hdot}^2+o_n(1)\\
\label{pythagore1b} 
\lf\|v_{1,n}\right\|^2_{L^2}=\sum_{j=1}^J \lf\|\partial_t U^j_{1}\left(\frac{-t_{j,n}} {\lambda_{j,n}}\right)\right\|^2_{L^2}+\lf\|w_{1,n}^J\right\|^2_{L^2}+o_n(1)\\
\label{pythagore2}
E(v_{0,n},v_{1,n})=\sum_{j=1}^J E\left(U^j_{\lin}\left(-\frac{t_{j,n}}{\lambda_{j,n}}\right),\partial_t U^j_{\lin}\left(-\frac{t_{j,n}} {\lambda_{j,n}}\right)\right)+E\left(w_{0,n}^J, w_{1,n}^J\right)+o_n(1).
\end{gather}
 \end{theoin}
Replacing $U^j_{\lin}(t,x)$ by  $V^j_{\lin}(t,x)=\frac{1}{\lambda_{j}^{\frac{N-2}{2}}}U_{\lin}^j\left(\frac{t-t_{j}}{\lambda_{j}},\frac{x-x_{j}}{\lambda_{j}}\right)$ for some good choice of the parameters $\lambda_j,\,t_j,\,x_j$, and extracting subsequences, we can always assume that one of the following two cases occurs
\begin{equation}
\label{choice_param}
\forall n,\;t_{j,n}=0\quad\text{or}\quad\lim_{n\to\infty} \frac{t_{j,n}}{\lambda_{j,n}}\in \{-\infty,+\infty\}. 
\end{equation} 
We will need the following bound on the parameters:
\begin{lemma}
\label{L:bound_param}
Let $v_n$ be as above and $\{\mu_n\}_n$ be a sequence of positive numbers. Assume
\begin{equation}
\label{cpct_infty}
\lim_{R\rightarrow +\infty} \limsup_{n\rightarrow+\infty} \int_{|x|\geq R\mu_n} \left(|\nabla v_{0,n}|^2+v_{1,n}^2\right)dx=0.
\end{equation} 
Then for all $j$, the sequences $\left\{\frac{\lambda_{j,n}}{\mu_n}\right\}_n$, $\left\{\frac{t_{j,n}}{\mu_n}\right\}_n$ and $\left\{\frac{x_{j,n}}{\mu_n}\right\}_n$ are bounded. Furthermore, there is at most one $j$ such that $\left\{\frac{\lambda_{j,n}}{\mu_n}\right\}_n$ does not converge to $0$.
\end{lemma}
\begin{proof}
 The case $\mu_n=1$ follows from \cite{BaGe99} (see p.154-155 for the proof). For the general case, apply the result of the case $\mu_n=1$ to the rescaled sequence $(\tilde{v}_{0,n},\tilde{v}_{1,n})_n$ defined by  $$\tilde{v}_{0,n}(t,x)=\mu_n^{N/2-1}v_{0,n}(\mu_n x),\quad \tilde{v}_{1,n}(t,x)=\mu_n^{N/2}v_{1,n}(\mu_n x),$$
\end{proof}

\begin{notation}
\label{N:nonlinearprofiles}
For any profile decomposition with profiles $\left\{U^j_{\lin}\right\}$ and parameters $\left\{\lambda_{j,n},t_{j,n},x_{j,n}\right\}$, we will denote  by $\left\{U^j\right\}$ the non-linear profiles associated with $\left\{U^j_{\lin}\left(\frac{-t_{j,n}}{\lambda_{j,n}}\right),\partial_t U^j_{\lin}\left(\frac{-t_{j,n}}{\lambda_{j,n}}\right)\right\}$, which are the unique solutions of \eqref{CP} such that for all $n$, $\frac{-t_{j,n}}{\lambda_{j,n}}\in I_{\max}\left(U^j\right)$ and
$$ \lim_{n\rightarrow +\infty} \left\|U^j\left(\frac{-t_{j,n}}{\lambda_{j,n}}\right)-U^j_{\lin}\left(\frac{-t_{j,n}}{\lambda_{j,n}}\right)\right\|_{\hdot}
+\left\|\partial_t U^j\left(\frac{-t_{j,n}}{\lambda_{j,n}}\right)-\partial_t U^j_{\lin}\left(\frac{-t_{j,n}}{\lambda_{j,n}}\right)\right\|_{L^2}=0.$$
\end{notation}
Assuming \eqref{choice_param}, the proof of the existence of $U^j$ follows from the local existence for \eqref{CP} if $t_{j,n}=0$ and from the existence of wave operators for equation \eqref{CP} if $t_{j,n}/\lambda_{j,n}$ tends to $\pm\infty$. 
By the Strichartz inequalities on the linear problem and the small data Cauchy theory (see \eqref{small_data_bis}), if 
$\lim_{n\rightarrow +\infty} \frac{-t_{j,n}}{\lambda_{j,n}}=+\infty$, then 
$T_+\left(U^j\right)=+\infty$ and
\begin{equation}
\label{finite_norm2}
 s_0>T_-\left(U^j\right) \Longrightarrow \|U^j\|_{S(s_0,+\infty)}<\infty,
\end{equation} 
an analoguous statement holds in the case $\lim_{n\rightarrow +\infty} \frac{t_{j,n}}{\lambda_{j,n}}=+\infty$.
\begin{notation}
\label{N:rescaled}
We will often write, for the sake of simplicity
\begin{equation*}
U^j_n(t,x)=\frac{1}{\lambda_{j,n}^{N/2-1}}U^j\left(\frac{t-t_{j,n}}{\lambda_{j,n}},\frac{x-x_{j,n}}{\lambda_{j,n}}\right),\quad
U^j_{\lin,n}(t,x)=\frac{1}{\lambda_{j,n}^{N/2-1}}U^j_{\lin}\left(\frac{t-t_{j,n}}{\lambda_{j,n}},\frac{x-x_{j,n}}{\lambda_{j,n}}\right).
\end{equation*}
\end{notation}

We will need the following approximation result, which follows from Theorem \ref{T:LTPT} and is an adaptation to the focusing case of the result of Bahouri-G\'erard (see the Main Theorem p. 135 in \cite{BaGe99}).

\begin{prop}
\label{P:lin_NL}
 Let $\{(v_{0,n},v_{1,n})\}_n$ be a bounded sequence in $\hdot\times L^2$, which admits the profile decomposition \eqref{decompo_profil}. Let $\theta_n\in (0,+\infty)$. Assume that for all $j,\,n$, 
\begin{equation}
\label{bounded_strichartz}
\forall j\geq 1, \quad\forall n,\;\frac{\theta_n-t_{j,n}}{\lambda_{j,n}}<T_+(U^j)\text{ and } \limsup_{n\rightarrow +\infty} \left\|U^j\right\|_{S\big(-\frac{t_{j,n}}{\lambda_{j,n}},\frac{\theta_n-t_{j,n}}{\lambda_{j,n}}\big)}<\infty.
\end{equation}
Let $u_n$ be the solution of \eqref{CP} with initial data $(v_{0,n},v_{1,n})$.
Then for large $n$, $u_n$ is defined on $[0,\theta_n)$,
\begin{equation}
\label{NL_bound}
\limsup_{n\rightarrow +\infty}\|u_n\|_{S(0,\theta_n)}<\infty,
\end{equation} 
and
\begin{equation}
\label{NL_profile} 
\forall t\in [0,\theta_n),\quad
u_n(t,x)=\sum_{j=1}^J U^j_n\left(t,x\right)+w^J_{n}(t,x)+r^J_n(t,x),
 \end{equation}
where 
\begin{equation}
\label{cond_rJn}
\lim_{n\rightarrow +\infty} \limsup_{J\rightarrow +\infty} \|r^J_n\|_{S(0,\theta_n)}+\sup_{t\in (0,\theta_n)} \left(\|\nabla r^J_n(t)\|_{L^2}+\|\partial_t r^J_n(t)\|_{L^2}\right)=0.
\end{equation}
An analoguous statement holds if $\theta_n<0$.
\end{prop}
\begin{remark}
\label{R:lin_NL}
Assume that for all $j$, at least one of the following occurs:
\begin{enumerate}
 \item \label{small_profile}
$\ds \|U^j_0\|_{\hdot}+\|U^j_1\|_{L^2}<\frac{\delta_0}{C_S}$, where the constant $C_S$ is given by the Strichartz estimate \eqref{strichartz} and $\delta_0$ by the small data theory;
\item \label{dispersive_profile}$\ds \lim_{n\rightarrow+\infty} \frac{-t_{j,n}}{\lambda_{j,n}}=+\infty$,
\item \label{finite_limit_profile}
$\ds \limsup_{n\rightarrow +\infty} \frac{\theta_n-t_{j,n}}{\lambda_{j,n}}<T_+\left(U^j\right).$
\end{enumerate}
Then \eqref{bounded_strichartz} holds. Indeed in case \eqref{small_profile}, it follows from \eqref{strichartz} and the small data theory. In case \eqref{dispersive_profile}, it follows directly from the small data theory: see \eqref{finite_norm2}. It remains to treat case \eqref{finite_limit_profile}, when $t_{j,n}=0$ or $-t_{j,n}/\lambda_{j,n}\to -\infty$. If $t_{j,n}=0$, then by definition 
$T_-\left(U^{j}\right)<0$  and \eqref{bounded_strichartz} is a consequence of \eqref{finite_norms} and \eqref{finite_limit_profile}. If $-t_{j,n}/\lambda_{j,n}\to -\infty$, then the analogue of \eqref{finite_norm2} for negative times and \eqref{finite_limit_profile} imply \eqref{bounded_strichartz}.
\end{remark}
\begin{remark}
\label{R:pythag} 
When $N$ is odd, under the assumptions of Proposition \ref{P:lin_NL}, we have localized pseudo-orthogonality properties for all time of the interval $(0,\theta_n)$ as follows: let $\tau_n\in (0,\theta_n)$ for all $n$, $\{\mu_n\}$ be any sequence of positive numbers and $\chi\in C^{\infty}_0(\RR^N)$ be radial and such that $\chi=1$ in a neighborhood of $0$. Then, if $\varphi=1$, $\varphi=\chi$ or $\varphi=1-\chi$, one can show after extraction the following Pythagorean expansion:
\begin{equation*}
\int \varphi\left(\frac{x}{\mu_n}\right)|\nabla_{t,x} u(\tau_n)|^2dx=\sum_{j=1}^J \int \varphi\left(\frac{x}{\mu_n}\right)\left|\nabla_{t,x} U_n^j(\tau_n)\right|^2dx +\int \varphi\left(\frac{x}{\mu_n}\right) \left|\nabla_{t,x}w_n^J(\tau_n)\right|^2dx+o_n(1),
\end{equation*}
where $|\nabla_{t,x}u|^2=(\partial_t u)^2+|\nabla_x u|^2$.
This follows easily from Claim \ref{C:ortho} in the appendix and we omit the proof.
\end{remark}

\begin{proof}[Sketch of proof of Proposition \ref{P:lin_NL}]

Denote by $F(u)=|u|^{\frac{4}{N-2}}u$ the nonlinearity. Let
$$ \tilde{u}_{n}^J(t,x)=\sum_{j=1}^J \frac{1}{\lambda_{j,n}^{\frac{N-2}{2}}}U^j\left(\frac{t-t_{j,n}}{\lambda_{j,n}},\frac{x-x_{j,n}}{\lambda_{j,n}}\right)+w^J_{n}(t,x).$$

We will apply Theorem \ref{T:LTPT} to $\tu_n$ and $u_n$ for large $n$.

We notice that there exists $J_0>0$ such that for all $j\geq J_0+1$, 
\begin{equation}
\label{small_strichartz}
\forall j\geq J_0+1,\quad \left\|U_{\lin}^j\right\|_{S(\RR)}<\delta_0,
 \end{equation} 
where $\delta_0$ is given by the small data theory for \eqref{CP}. Indeed, it is an immediate consequence of the Pythagorean expansions \eqref{pythagore1a}, \eqref{pythagore1b} and Strichartz estimates. Thus we can use the small data theory which implies by \eqref{small_data_bis}
\begin{equation}
\label{small_strichartzNL}
\forall j\geq J_0+1,\quad \left\|U^j\right\|_{S(\RR)}\leq C\lf(\|U^j_0\|_{\hdot}+\|U^j_1\|_{L^2}\rg).
 \end{equation} 
Fixing a large $J$, one can show, as a consequence of the orthogonality \eqref{ortho_param} of the parameters,
\begin{equation*}
\lf\| \sum_{j=1}^J \frac{1}{\lambda_{j,n}^{\frac{N-2}{2}}}U^j\left(\frac{t-t_{j,n}}{\lambda_{j,n}},\frac{x-x_{j,n}}{\lambda_{j,n}}\right) \rg\|_{S(0,\theta_n)}^{\frac{2(N+1)}{N-2}}=\sum_{j=1}^J \lf\|\frac{1}{\lambda_{j,n}^{\frac{N-2}{2}}}U^j\left(\frac{t-t_{j,n}}{\lambda_{j,n}},\frac{x-x_{j,n}}{\lambda_{j,n}}\right) \rg\|_{S(0,\theta_n)}^{\frac{2(N+1)}{N-2}} +o_n(1).
 \end{equation*} 
Combining with \eqref{small_strichartzNL} and the Pythagorean expansions \eqref{pythagore1a}, \eqref{pythagore1b}, we get
\begin{equation*}
\label{bound_S}
\limsup_{J\to \infty}\limsup_{n\to\infty}
\|\tu_n^J\|_{S(0,\theta_n)}<\infty.
 \end{equation*} 
Let $J\geq J_0$ and $e_{n}^J=(\partial_t^2-\Delta)\tu_{n}^J-F(\tu_{n}^J)$. Then
\begin{equation*}
e_{n}^J(t,x)=\sum_{j=1}^J F\left(\frac{1}{\lambda_{j,n}^{\frac{N-2}{2}}}U^j\left(\frac{t-t_{j,n}}{\lambda_{j,n}}\right)\right) -F\left(\sum_{j=1}^J \frac{1}{\lambda_{j,n}^{\frac{N-2}{2}}}U^j\left(\frac{t-t_{j,n}}{\lambda_{j,n}}\right)+w^J_{n}\right).
\end{equation*}
From \eqref{bounded_strichartz} and again the orthogonality \eqref{ortho_param} of the parameters, we can deduce
\begin{gather*}
\lim_{J\rightarrow +\infty}\limsup_{n\rightarrow +\infty}\left\|D_x^{1/2}e_{n}^J\right\|_{N(0,\tau_n)}=0.
\end{gather*}
Furthermore 
$$\tilde{u}_{n}^J(0)=u_{n}^J(0),\quad \partial_t \tilde{u}_{n}^J(0)=\partial_t u_{n}^J(0),$$
which yields by Theorem \ref{T:LTPT} the conclusion of the proposition.
\end{proof}

We will also need the following technical claim. The proof is postponed to Appendix \ref{A:ortho}.
\begin{claim}
\label{C:dispersive_tw}
Assume that $N$ is odd.
Let $w_n$ be a sequence of the radial solutions to the linear wave equation \eqref{lin_wave} with bounded energy and such that 
\begin{equation}
\label{dispersive_w}
\lim_{n\to\infty} \|w_n\|_{S(\RR)}=0.
\end{equation} 
Let $(w_{0,n},w_{1,n})$ be the initial data of $w_n$, $\chi\in C_0^{\infty}(\RR^N)$, radial and such that $\chi=1$ around the origin, and $\tlambda_n$ be a sequence of positive numbers and consider the solution $\tw_n$ to \eqref{lin_wave} with initial data $(\tw_{0,n},\tw_{1,n})=\left(\varphi\big(|x|/\tlambda_n\big)w_{0,n},\varphi\big(|x|/\tlambda_n\big)w_{1,n}\right)$, where $\varphi=\chi$ or $\varphi=1-\chi$. Then
\begin{equation}
\label{dispersive_tw}
\lim_{n\to\infty} \|\tw_n\|_{S(\RR)}=0.
\end{equation} 
\end{claim}

\section{Description of general type II blow-up solutions}
\label{S:reg_sing}
In this section we consider a general type II blow-up solution of \eqref{CP} in space dimension $N\in\{3,4,5\}$, that is a solution $u$ bounded in the energy space and such that $T_+(u)<\infty$. We do not assume that $u$ is spherically symetric. 
\begin{defi}
\label{D:regular}
Let $x_0\in \RR^N$. We will say that the point $x_0$ is \emph{regular} if
\begin{equation}
\label{regular}
\forall \eps>0,\; \exists R,\; \forall t\in [0,T_+(u)),\quad \int_{|x-x_0|\leq R} |\nabla u|^2+\frac{u^2}{|x-x_0|^2}+(\partial_t u)^2 \leq \eps.
\end{equation}
If $x_0$ is not regular, we will say that it is \emph{singular}. We will denote by $S$ the set of singular points.
\end{defi}

\begin{theoin}
\label{T:general_blowup}
Let $u$ be a solution of \eqref{CP} with type II blow-up forward in time, and $T_+=T_+(u)$ the blow-up time. Then there exists $K\in \NN^*$ and $K$ distinct points $m_1,\ldots,m_K$ of $\RR^N$ such that $S=\{m_1,\ldots,m_K\}$. Furthermore there exists $(v_0,v_1)\in \hdot\times L^2$ such that 
\begin{equation}
\label{weak_limit}
(u(t),\partial_t u(t))\xrightharpoonup[t\to T_+]{} (v_0,v_1)\text{ weakly in }\hdot\times L^2.
\end{equation}
If 
$\varphi\in C^{\infty}_0(\RR^N)$ is equal to $1$ around each singular point, we have
\begin{equation}
\label{CV_to_v}
\lim_{t\to T^+}\left\|\left(1-\varphi\right)(u(t)-v_0)\right\|_{\hdot}+\left\|\left(1-\varphi\right)(\partial_t u(t)-v_1)\right\|_{L^2}=0.
\end{equation} 
Furthermore, if $k\in 1\ldots K$,
\begin{gather}
\label{limsup}
\limsup_{t\rightarrow T_+} \int_{|x-m_k|\leq |t-T_+|} |\nabla u(t,x)|^2+|\partial_t u(t,x)|^2\geq \int |\nabla W|^2\\
\label{liminf}
\liminf_{t\rightarrow T_+} \int_{|x-m_k|\leq |t-T_+|} |\nabla u(t,x)|^2+|\partial_t u(t,x)|^2\geq \frac{2}{N}\int |\nabla W|^2.
\end{gather}
\end{theoin}
\begin{defi}
\label{D:regular_part}
Under the assumptions of Theorem \ref{T:general_blowup}, le $v$ be the solutions of \eqref{CP} such that $(v(T_+),\partial_t v(T_+))=(v_0,v_1)$. We will call $v$ the \emph{regular part} of $u$ at the blow-up time $T_+$, and $a=u-v$ the \emph{singular part} of $u$. Note that \eqref{CV_to_v} implies, together with the finite speed of propagation, that 
$$\supp a\subset \bigcup_{k=1}^K \Big\{(t,x), \;|x-m_k|\leq \lf|t-T_+\rg|\Big\}.$$
\end{defi}

This section is divided into two parts. In \S \ref{SS:generality}, we perform a first analysis of the behaviour of $u$ around each singular point, showing \eqref{weak_limit} and \eqref{CV_to_v}. In \S \ref{SS:bound_below}, we write a profile decomposition of the solution around each singular point to show \eqref{limsup} and \eqref{liminf}.
 
We will assume in all the sequel without loss of generality that the blow-up time is $T_+(u)=1$.

\subsection{Generality on regular and singular points}
\label{SS:generality}
\begin{lemma}
\label{L:compact}
There exists a constant $\delta_1>0$ with the following properties:
\begin{enumerate}
\item \label{compact_local} for all $x_0\in \RR^N$, $t_0\in (0,1)$ and $R>0$, if 
$$ \int_{|x-x_0|\leq |t_0-1|+R}|\nabla u(t_0)|^2+|\partial_t u(t_0)|^2+\frac{1}{|x-x_0|^2}|u(t_0)|^2\leq \delta_1,$$
and $\varphi\in C^{\infty}_0(\RR^N)$ has compact support in $\{|x-x_0|\leq R\}$, then
$(\varphi u(t),\varphi\partial_t u(t))$ has a limit in $\hdot\times L^2$ as $t\xrightarrow{\scriptscriptstyle <}1$;
\item \label{compact_infty} for all $t_0\in(0,1)$, and $R>0$, if
$$ \int_{|x|\geq R}|\nabla u(t_0)|^2+|\partial_t u(t_0)|^2+\frac{1}{|x|^2}|u(t_0)|^2\leq \delta_1,$$
and $\varphi \in C^{\infty}(\RR^N)$ is equal to $1$ at infinity and is supported in the set $|x|\geq R+|1-t_0|$, then $(\varphi u(t),\varphi\partial_t u(t))$ has a a limit in $\hdot\times L^2$ as $t\xrightarrow{\scriptscriptstyle <}1$.
\end{enumerate}
\end{lemma}
\begin{proof}
Let us prove \eqref{compact_local}. Assume that for some parameter $\eta_0>0$ to be determined later,
\begin{equation*}
\int_{|x-x_0|\leq |t_0-1|+R}|\nabla u(t_0)|^2+|\partial_t u(t_0)|^2+\frac{1}{|x-x_0|^2}|u(t_0)|^2\leq \eta_0E(W,0).
\end{equation*}
If $\eta_0$ is chosen small enough, then, by a standard extension theorem, there exist $\tilde{u}_0\in \hdot$, $\tilde{u}_1\in L^2$ compactly supported on $\RR^N$ and such that
\begin{gather}
\label{data_equal}
\tilde{u}_0(x)=u(t_0,x) \text{ and } \tilde{u}_1=\partial_tu(t_0,x)\text{ if } |x-x_0|\leq |t_0-1|+R,\\
\label{small_norm}
\int_{\RR^N}|\nabla \tilde{u}_0|^2+|\tilde{u}_1|^2+\frac{1}{|x-x_0|^2}|\tilde{u}_0|^2\leq C\eta_0 E(W,0)<E(W,0).
\end{gather} 

Consider the solution $\tilde{u}$ of \eqref{CP} with initial condition $(\tilde{u}_0,\tilde{u}_1)$ at $t=t_0$. By \eqref{small_norm}, we have
\begin{align*}
E(\tilde{u}_0,\tilde{u}_1)&=\frac{1}{2}\int_{\RR^N} |\nabla \tilde{u}_0|^2+ \frac{1}{2}\int_{\RR^N} |\tilde{u}_1|^2-\frac{N-2}{2N}\int_{\RR^N} |\tilde{u}_0|^{\frac{2N}{N-2}}<E(W,0)\\
\|\nabla\tilde{u}_0\|_{L^2}^2&\leq E(W,0)<\int |\nabla W|^2.
\end{align*}
By the result of Kenig-Merle \cite{KeMe08}, $\tilde{u}$ is globally defined. The mapping $t\mapsto (\tilde{u}(t),\partial_t\tilde{u}(t))$ is continuous from $\RR$ to $\hdot\times L^2$. 
By the finite speed of propagation and \eqref{data_equal},
$$ \forall t\in [t_0,1], \; \forall x\in \RR^N,\quad |x-x_0|\leq |t-1|+R\Longrightarrow u(t,x)=\tilde{u}(t,x),\; \partial_tu(t,x)=\partial_t\tilde{u}(t,x),$$
In particular, $(\varphi u(t),\varphi \partial_t u(t))=(\varphi \tilde{u}(t),\varphi \partial_t \tilde{u}(t))$ has a limit as $t\rightarrow 1$, which concludes the proof of case \eqref{compact_local}.

Case \eqref{compact_infty} is similar. Indeed in this case, if $\delta_1$ is small enough, there exist $\tilde{u}_0$ and $\tilde{u}_1$ such that 
\begin{gather*}
\tilde{u}_0(x)=u(t_0,x) \text{ and } \tilde{u}_1=\partial_tu(t_0,x)\text{ if } |x|\geq R,\\
\int_{\RR^N}|\nabla \tilde{u}_0|^2+|\tilde{u}_1|^2+\frac{1}{|x|^2}|\tilde{u}_0|^2<E(W,0).
\end{gather*} 
Consider the solution $\tilde{u}$ with initial data $(\tilde{u}_0,\tilde{u}_1)$ at $t=t_0$. By the finite speed of propagation, $u$ and $\tilde{u}$ coincide if $|x|>|t_0-t|+R$, and the result follows again by the global existence of $\tilde{u}$.
\end{proof}
\begin{corol}
\label{C:norm_singular}
For any singular point $m$, for all $t\in I_{\max}=I_{\max}(u)$,
\begin{equation}
\label{bound_singular1}
\delta_1 \leq \int_{|x-m|\leq |t-1|}|\nabla u(t)|^2+|\partial_t u(t)|^2+\frac{1}{|x-m|^2}|u(t)|^2,
\end{equation} 
where $\delta_1$ is given by Lemma \ref{L:compact}. Furthermore, the set $S$ of singular points is finite. 
\end{corol}
\begin{proof}
The finiteness of $S$ follows immediately from \eqref{bound_singular1} and the fact that the blow-up is of type II.
 
Let us show \eqref{bound_singular1}. We argue by contradiction. Consider a singular point $m$, and assume that for some $t_0\in I_{\max}$ and $\eps>0$,
\begin{equation*}
\int_{|x-m|\leq |t_0-1|+\eps}|\nabla u(t_0)|^2+|\partial_t u(t_0)|^2+\frac{1}{|x-m|^2}|u(t_0)|^2<\delta_1.
\end{equation*}
Let $\varphi\in C^{\infty}_0(\RR^N)$ such that $\varphi(x)=0$ if $|x-m|\geq \eps$ and $\varphi(x)=1$ if $|x-m|\leq \frac{\eps}{2}$. By Lemma \ref{L:compact}, $(\varphi u,\varphi \partial_t u)$ converges in $\hdot\times L^2$ as $t$ tends to $1$, contradicting, in view of the continuous embedding of $\hdot$ into $L^2\left(\frac{1}{|x-m|^2}dx\right)$ the assumption that $m$ is a singular blow-up point.

We have proven that for all $t\in I_{\max}$, for all $\eps>0$,
$$
\delta_1\leq \int_{|x-m|\leq |t-1|+\eps}|\nabla u(t)|^2+|\partial_t u(t)|^2+\frac{1}{|x-m|^2}|u(t)|^2,
$$
concluding the proof of \eqref{bound_singular1}.
\end{proof}
We are now ready to prove \eqref{weak_limit} and \eqref{CV_to_v} of Theorem \ref{T:general_blowup}. Let us first show that $(u(t),\partial_t u(t))$ has a weak limit in $\hdot\times L^2$ as $t\xrightarrow{\ssstyle <}1$. It is equivalent to show that all weak limits of sequences $\big\{(u(t_n),\partial_t u(t_n))\big\}_n$ where $t_n\xrightarrow{\ssstyle <}1$, coincide. For this, notice that the definition of a regular point and Lemma \ref{L:compact} show that if $(v_0,v_1)$ and $(\tv_0,\tv_1)$ are such limits, then they must coincide around any regular point. As the set of singular point is finite, this shows as desired that $(v_0,v_1)=(\tv_0,\tv_1)$. Denote by 
$$ (v_0,v_1)=\underset{t\to 1}{w\mbox{-}\lim}\, (u(t),\partial_t u(t)).$$ 
By point \eqref{compact_local} of Lemma \ref{L:compact}, $(u,\partial_t u)$ has a limit in $\hdot_{\loc}\lf(\RR^N\setminus S\rg)\times L^{2}_{\loc}\lf(\RR^N\setminus S\rg)$ as $t$ goes to $1$. The uniqueness of limits shows that this limit must be $(v_0,v_1)$. Using point \eqref{compact_infty} of Lemma \ref{L:compact}, we get that the convergence to $v$ is also global, hence \eqref{CV_to_v}.

We finish this part by noting that there is at least one singular point. If not, \eqref{CV_to_v} shows that $(u(t),\partial_t u(t))$ has a limit as $t\rightarrow 1$, which shows that $1$ is not the maximal positive time of existence, a contradiction.

\subsection{Bounds from below on the norm of the main profile}
In this subsection we will complete the proof of Theorem \ref{T:general_blowup} by
studying the behavior of $u$ in the neighborhood of singular points by using a profile decomposition. We assume that 
$$ 0\in S.$$
Consider an increasing sequence $\{\tau_n\}\in (t_0,1)^{\NN}$ that tends to $1$ and a function $\psi\in C^{\infty}_0\lf(\RR^N\rg)$ such that $\psi=1$ close to $0$ and $\supp\psi\cap S=\{0\}$. After extracting a subsequence, we can assume that there exists a profile decomposition
 \begin{equation}
\label{decompo_solution1}
\left\{
\begin{aligned}
 \psi u(\tau_n)-\psi v(\tau_n)&=\sum_{j=1}^J \frac{1}{\lambda_{j,n}^{N/2-1}}U^j_{\lin}\left(\frac{- t_{j,n}}{\lambda_{j,n}},\frac{x-x_{j,n}}{\lambda_{j,n}}\right) +w_{0,n}(x)\\
\psi\frac{\partial u}{\partial t}(\tau_n)-\psi\frac{\partial u}{\partial t}(\tau_n)&=\sum_{j=1}^J \frac{1}{\lambda_{j,n}^{N/2}}\frac{\partial U^j_{\lin}}{\partial t} \left(\frac{- t_{j,n}}{\lambda_{j,n}},\frac{x-x_{j,n}}{\lambda_{j,n}}\right)+w_{1,n}(x),
 \end{aligned}\right.
 \end{equation}
where $U^j_{\lin}$ is a solution of the linear wave equation \eqref{lin_wave} with initial conditions $\left(U_0^j,U_1^j\right)$.

As $\psi(u-v)$ is supported in $\{|x|\leq 1-t\}$, when $t$ is close to $1$, Lemma \ref{L:bound_param} implies
\begin{equation}
\label{concentration_param1}
\forall j\geq 1,\;\exists C_j,\; \forall n,\quad \left|\lambda_{jn}\right|+\left|t_{j,n}\right|+\left|x_{j,n}\right|\leq C_j(1-\tau_n).
\end{equation} 
Let us first show:
\label{SS:bound_below}
\begin{lemma}
\label{L:profiles}
Reorder the decomposition \eqref{decompo_solution1} so that 
\begin{equation}
\label{reorder}
\left\|\nabla U_0^{1}\right\|_{L^2}^2+\left\|U_1^{1}\right\|_{L^2}^2=\sup_{j\geq 1}\left(\left\| \nabla U_0^{j}\right\|^2_{L^2}+\left\| U_1^{j}\right\|_{L^2}^2\right).
\end{equation} 
Then
\begin{equation}
\label{bb_big_profile}
\left\|\nabla U_0^{1}\right\|_{L^2}^2+\left\|U_1^{1}\right\|_{L^2}^2\geq \frac{2}{N}\|\nabla W\|_{L^2}^2.
\end{equation}
\end{lemma}
Lemma \ref{L:profiles}, together with the Pythagorean expansions \eqref{pythagore1a} and \eqref{pythagore1b} implies immediately \eqref{liminf}.
\begin{remark}
\label{R:bound_d3}
In space dimension $N=3$, we have the following immediate corollary of Lemma \ref{L:profiles}. Assume
$$ \liminf_{t\to 1^-}\left(\|\nabla u(t)\|_{L^2}^2+\|\partial_t u(t)\|_{L^2}^2\right)< \frac{4}{3}\|\nabla W\|^2_{L^2},$$
then there is only one singular point. See Remark \ref{R:bound_dN} below for an improvement.
\end{remark}

\begin{proof}[Proof of Lemma \ref{L:profiles}]
Assume that 
$$\|\nabla U_0^{1}\|_{L^2}^2+\|U_1^{1}\|_{L^2}^2< \frac{2}{N}\|\nabla W\|_{L^2}^2,$$
and thus for all $j\geq 1$, 
$$\|\nabla U^j_0\|_{L^2}^2+\|U^j_1\|_{L^2}^2< \frac{2}{N}\|\nabla W\|_{L^2}^2.$$
Using that $2E(f,g)\leq \|\nabla f\|^2_{L^2}+\|g\|^2_{L^2}$ and that $E(W,0)=\frac{1}{N}\|\nabla W\|_{L^2}^2$, we get that there exists an $\eps_0>0$ such that for all $j,n$
$$E\left(U_{\lin}^j\left(-t_{j,n}/\lambda_{j,n}\right),\partial_t U_{\lin}^j\left( -t_{j,n}/\lambda_{j,n}\right)\right) \leq E(W,0)- \eps_0,\quad \|\nabla U_0^j\|_{L^2}^2\leq \|\nabla W\|_{L^2}^2-\eps_0.$$
Then according to \cite{KeMe08}, for all $j$, $U^j$ is globally define and scatters. By Proposition \ref{P:lin_NL} the solution with initial condition $(\psi u(\tau_n),\psi \partial_t u(\tau_n))$ is globally defined and scatters for large $n$. Using the finite speed of propagation, we get a contradiction with the fact that $0$ is singular. Hence \eqref{bb_big_profile}.
\end{proof}

It remains to show \eqref{limsup}.
We first recall the following scattering result (see \cite[Corollary 7.4]{KeMe08}):
\begin{prop}
\label{P:scatteringKeMe}
Let $u$ be a solution of \eqref{CP} such that
$$ \limsup_{t\to T^+(u)}\|\nabla u(t)\|_{L^2}^2+\|\partial_t u(t)\|_{L^2}^2<\|\nabla W\|_{L^2}^2.$$
Then $u$ is globally defined and scatters.
\end{prop}
The following proposition implies \eqref{limsup} by the Pythagorean expansions \eqref{pythagore1a}, \eqref{pythagore1b}:
\begin{prop}
\label{P:concentration}
Let $\eps_0>0$. There exists a sequence $\{\ttau_n\}\in (t_0,1)^{\NN}$ that tends to $1$ such that $\lf(\psi a(\ttau_n),\psi\partial_t a(\ttau_n)\rg)$ admits a profile decomposition $\lf\{\tU_{\lin}^j\rg\}_j$, $\lf\{\tlambda_{j,n},\tx_{j,n},\tilde{t}_{j,n}\rg\}_{j,n}$ 
such that $\tilde{t}_{1,n}=0$ and 
\begin{equation}
\label{large_profile}
\left\|\nabla \tU^{1}_0\right\|_{L^2}^2+\left\|\tU^1_1\right\|_{L^2}^2\geq \|\nabla W\|^2_{L^2}-\eps_0.
\end{equation} 
\end{prop}
\begin{proof}
We follow the lines of the proof of \cite[Corollary 7.5]{KeMe08}. In all the proof, we will always work up to the extraction of a subsequence for sequences indexed by $n$. In particular, any real sequence indexed by $n$ will be assumed to have a limit in $\RR\cup\{\pm\infty\}$.

Consider an increasing sequence $\{\tau_n\}\in (t_0,1)^{\NN}$ that tends to $1$. Let $\tilde{u}_n$ and $\tilde{v}_n$ be the solutions of \eqref{CP} such that
$$ (\tilde{u}_n,\partial_t \tilde{u}_n)_{\restriction t=\tau_n}=(\psi u(\tau_n),\psi\partial_t u(\tau_n)),\quad (\tilde{v}_n,\partial_t \tilde{v}_n)_{\restriction t=\tau_n}=(\psi v(\tau_n),\psi\partial_t v(\tau_n)).$$
By finite speed of propagation, and the fact that $x=0$ is a singular point for $u$, $T_+(\tilde{u}_n)\leq 1$. Furthermore, $(\psi v(\tau_n),\psi\partial_t v(\tau_n))$ has a limit in $\hdot\times L^2$ as $n\to \infty$, which implies that there exists a small $t_0>0$ such that $\tilde{v}_n(\tau_n+t)$ is well defined for large $n$ and $|t|\leq t_0$.

After extracting a subsequence, there exists a profile decomposition with profiles $\lf\{U^j_{\lin}\rg\}$ and parameters $\{\lambda_{j,n},t_{j,n},x_{j,n}\}$ associated to the sequence $\big(\tilde{u}_n(\tau_n)-\tilde{v}_n(\tau_n),\partial_t\tilde{u}_n(\tau_n)-\partial_t\tilde{v}_n(\tau_n)\big)_n$.
The fact that $\tilde{\psi}(u-v)$ is supported in $\{|x|\leq 1-t\}$ and Lemma \ref{L:bound_param} imply
\begin{equation}
\label{concentration_param}
\forall j\geq 1,\;\exists C_j,\; \forall n,\quad \left|\lambda_{jn}\right|+\left|t_{j,n}\right|+\left|x_{j,n}\right|\leq C_j(1-\tau_n).
\end{equation}

Let us consider the associated nonlinear profiles $U^j$ (see Notation \ref{N:nonlinearprofiles}).
Reordering the profiles, we get a $J_0$ such that
\begin{equation*}
\forall j\leq J_0,\quad \|U^j\|_{S\left(0,T_+(U^j)\right)}=\infty,\quad \forall j\geq J_0+1,\quad \|U^j\|_{S\left(0,T_+(U^j)\right)}<\infty.
\end{equation*}

By the finite blow-up criterion, $T_+(U^j)=+\infty$ if $j\geq J_0+1$.
By Proposition \ref{P:lin_NL} there is at least one solution $U^j$ that does not scatter forward in time (otherwise we would have $T_+(u)>1$), and thus $J_0\geq 1$.

For $1\leq j\leq J_0$, 
$\lim_{n} \frac{-t_{j,n}}{\lambda_{j,n}}=\ell_j\in \left\{{-\infty}\right\}\cup \RR$
(the case $\ell_j=+\infty$ is excluded as the nonlinear profile does not scatter forward in time). If $\ell_j$ is finite, the corresponding profile is compact up to scaling and translation, and we may assume $t_{j,n}=0$. Thus
\begin{equation}
\label{cond_tjn}
 \forall j\in \{1,\ldots,J_0\},\quad t_{j,n}=0 \text{ or } \lim_{n\rightarrow +\infty} \frac{-t_{j,n}}{\lambda_{j,n}}=-\infty.
\end{equation} 
By Proposition \ref{P:scatteringKeMe}, for all $j\in \{1,\ldots, J_0\}$, there exists a time $T_j$ such that
\begin{equation}
\label{defTj}
T_-\left(U^j\right)<T_j <T_+\left(U^j\right) \text{ and } \|\nabla u(T_j)\|_{L^2}^2+\|\partial_t u(T_j)\|^2_{L^2}\geq \|\nabla W\|^2_{L^2}-\eps_0,
\end{equation} 
furthermore, using that $T_+\left(U^j\right)>0$ if $t_{j,n}=0$, we may choose $T_j$ such that
$$ \left(\forall n,\; t_{j,n}=0\right)\Longrightarrow T_j>0.$$
Extracting subsequences and reordering the profiles, we may assume
\begin{equation}
\label{t1n_min}
\forall n,\quad  t_{1,n}+\lambda_{1,n}T_1=\min_{1\leq j\leq J_0}\left(t_{j,n}+\lambda_{j,n}T_j\right).
\end{equation}
Denote by $\theta_n=t_{1,n}+\lambda_{1,n}T_1$. Note that $\theta_n\geq 0$ for large $n$ and, by \eqref{concentration_param} 
\begin{equation}
\label{limit0}
\lim_{n\rightarrow +\infty}  \theta_n =0.
\end{equation}


For all $j$, we have, by definition of $\theta_n$, $\frac{\theta_n-t_{j,n}}{\lambda_{j,n}}\leq T_j<T_+(U^j)$. According to Remark \ref{R:lin_NL} we can use 
Proposition \ref{P:lin_NL} which shows that $\tau_n+\theta_n<T_+(\tilde{u}_n)\leq 1$, that $\big\{\|\tilde{u}_n\|_{S(\tau_n,\tau_n+\theta_n)}\big\}_n$ is bounded and
\begin{multline}
\label{decompo_solution''}
 \tilde{u}_n(\tau_n+t)=\tilde{v}_n(\tau_n+t)+ \sum_{j=1}^J \frac{1}{\lambda_{j,n}^{\frac{N-2}{2}}}U^j\left(\frac{t- t_{j,n}}{\lambda_{j,n}},\frac{x-x_{j,n}}{\lambda_{j,n}}\right)\\
+w^J_{n}(t,x)+r^J_n(t,x), \quad t\in(0,\theta_n),
\end{multline}
where $r^J_n$ satisfies \eqref{cond_rJn}. If $j\geq 1$, there exists (extracting if necessary) a linear wave $\tU_{\lin}^j$ such that
$$ \lim_{n\to \infty} \lf\|U^{j}\lf(\frac{\theta_n-t_{j,n}}{\lambda_{j,n}}\rg)-\tU_{\lin}^{j}\lf(\frac{\theta_n-t_{j,n}}{\lambda_{j,n}}\rg) \rg\|_{\hdot}+ \lf\|\partial_tU^{j}\lf(\frac{\theta_n-t_{j,n}}{\lambda_{j,n}}\rg)-\partial_t\tU_{\lin}^{j}\lf(\frac{\theta_n-t_{j,n}}{\lambda_{j,n}}\rg) \rg\|_{L^2}=0.$$
Indeed, if $\lf\{\frac{\theta_n-t_{j,n}}{\lambda_{j,n}}\rg\}_n$ converges this is obvious, if it goes to $-\infty$, it implies that  $\lf\{\frac{\theta_n-t_{j,n}}{\lambda_{j,n}}\rg\}_n$ also goes to $-\infty$, and we can take $\tU^j_{\lin}=U^j_{\lin}$.
Writing $\ttau_n=\tau_n+\theta_n$ and $\tilde{t}_{j,n}=t_{j,n}-\theta_n$ we get by \eqref{decompo_solution''}, 
 \begin{equation}
\label{decompo_tilde_un}
\left\{
\begin{aligned}
 \left(\tilde{u}_n-\tilde{v}_n\right)(\ttau_n)&=\sum_{j=1}^J \frac{1}{\lambda_{j,n}^{\frac{N}{2}-1}}\tU^j_{\lin}\left(\frac{- \tilde{t}_{j,n}}{\lambda_{j,n}},\frac{x-x_{j,n}}{\lambda_{j,n}}\right) +w_{n}^J(\theta_n)+r_n^J(\theta_n)+o_n(1)\\
\partial_t\left(\tilde{u}_n-\tilde{v}_n\right)(\ttau_n)&=\sum_{j=1}^J \frac{1}{\lambda_{j,n}^{\frac N2}}\partial_t \tU^j_{\lin} \left(\frac{- \tilde{t}_{j,n}}{\lambda_{j,n}},\frac{x-x_{j,n}}{\lambda_{j,n}}\right)+\partial_t w_{n}^J(\theta_n)+\partial_t r_n^J(\theta_n)+o_n(1),
 \end{aligned}\right.
 \end{equation}
This a profile decomposition for the sequence $\big(\tilde{u}_n(\ttau_n)-\tilde{v}_n(\ttau_n),\partial_t \tilde{u}_n(\ttau_n)-\partial_t\tilde{v}_n(\ttau_n)\big)$, with profiles $\tU_{\lin}^j$ and parameters $\lambda_{j,n},x_{j,n},\tilde{t}_{j,n}$.
Note that the orthogonality of the parameters follows directly from the equality $\tilde{t}_{j,n}-\tilde{t}_{k,n}=t_{j,n}-t_{k,n}$.

Next notice that by finite speed of propagation and the definitions of $\tilde{u}_n$ and $\tilde{v}_n$, there exists a $r_0>0$ such that, if $n$ is large and $|x|<r_0$ then $\tilde{u}_n(\ttau_n)=u(\ttau_n)$, $\partial_t \tilde{u}_n(\ttau_n)=\partial_t u(\ttau_n)$, $\tilde{v}_n(\ttau_n)=v(\ttau_n)$ and $\partial_t \tilde{v}_n(\ttau_n)=\partial_t v(\ttau_n)$. Using that $u(\ttau_n)-v(\ttau_n)$ and $\partial_t u(\ttau_n)-\partial_t v(\ttau_n)$ are supported in the set $\{|x|\leq 1-\ttau_n\}$, one can replace, in the decomposition \eqref{decompo_tilde_un}, $\tilde{u}_n$ and $\tilde{v}_n$ by $\psi u$ and $\psi v$.

Finally,
$\frac{\theta_n- t_{1,n}}{\lambda_{1,n}}=T_1$.
Thus the first profile $\tU^1$ in the decomposition \eqref{decompo_tilde_un} is compact up to modulation, and we may assume $\tilde{t}_{1,n}=0$ as announced. The inequality \eqref{large_profile} follows from the choice of $T_1$.
\end{proof}

\begin{remark}
\label{R:bound_dN}
We can improve Remark \ref{R:bound_d3} as follows.
If for some $t_0\in (0,1)$,
$$ \sup_{t\in (t_0,1)}\|\nabla u(t)\|_{L^2}^2+\|\partial_t u(t)\|^2_{L^2}< \left(1+\frac{2}{N}\right)\|\nabla W\|_{L^2}^2,$$
then there is only one singular point. This is a direct consequence of \eqref{limsup} and \eqref{liminf}.
\end{remark}

%


\section{Finite speed of propagation and exclusion of small exterior profiles}
\label{S:general_radial}
In the two next sections, we assume that $N=3$ and that $u$ is spherically symmetric, blows up at time $T=1$ and satisfies 
\begin{equation}
\label{bound_nabla'}
\sup_{\tau_0\leq t<1} \sqrt{\|\nabla u(t)\|_{L^2}^2+\|\partial_t u(t)\|_{L^2}^2}\leq C_0.
\end{equation} 
In these two sections we will not make any further assumption on $C_0>0$.
By spherical symmetry $0$ is the  only singular point. We denote by 
$$ a(t,x)=u(t,x)-v(t,x)$$
the singular part of $u$ at the blow-up time $t=1$ (see Definition \ref{D:regular_part}).

The main result of this section (Proposition \ref{P:ext_profile}), shown in \S \ref{SS:No_exterior}, is that the norm of the most exterior profile of any profile decomposition of a sequence $(a(t_n),\partial_t a(t_n))$ is bounded from below by an universal constant independent of the solution. 

\subsection{Linear behavior}
\label{SS:linear}
We start by two preliminaries results on the linear problem, valid in odd dimension only, that will be needed in the sequel. The first one follows from Huygens principle:
\begin{lemma}
\label{L:lin_odd}
Assume that $N$ is odd. Let $v$ be a solution of the linear wave equation \eqref{lin_wave}, with initial conditions $(v_0,v_1)$, $\lf\{\lambda_n\rg\}_n$, $\lf\{t_n\rg\}_n$ be two real sequences, with $\lambda_n$ positive. 
$$ v_n(t,x)=\frac{1}{\lambda_n^{N/2-1}}v\lf(\frac{t}{\lambda_n},\frac{x}{\lambda_n}\rg).$$
and assume $\lim_{n\to \infty} \frac{t_n}{\lambda_n}=\ell\in [-\infty,+\infty]$.
Then, if $\ell=\pm \infty$.
$$  \lim_{R\to \infty} \limsup_{n\to \infty} \int_{\big||x|-|t_n|\big|\geq R\lambda_n}  |\nabla v_n(t_n)|^2+\frac{1}{|x|^2}|v_n(t_n)|^2+\lf(\partial_t v_n(t_n)\rg)^2dx=0$$
and if $\ell\in \RR$,
$$  \lim_{R\to \infty} \limsup_{n\to \infty} \int_{\substack{\{|x|\geq R \lambda_n\}\\ \cup\{|x|\leq \frac{1}{R}\lambda_n\}}} |\nabla v_n(t_n)|^2+\frac{1}{|x|^2}|v_n(t_n)|^2+\lf(\partial_t v_n(t_n)\rg)^2dx=0.$$
\end{lemma}
\begin{proof}
This is a classical property.  In the case $\ell\in \RR$, just notice that 
$$v_n(t_n,x)=\frac{1}{\lambda_n^{N/2-1}}v\lf(\ell,\frac{x}{\lambda_n}\rg)+o_n(1)\text{ in }\hdot,\quad \partial_t v_n(t_n,x)=\frac{1}{\lambda_n^{N/2}}\partial_t v\lf(\ell,\frac{x}{\lambda_n}\rg)+o_n(1)\text{ in }L^2,$$
which implies the announced estimate (in this case we do not need any assumption on the parity of $N$). 

Let us treat the case $\ell=\pm\infty$. 
Let $\eps>0$, $\chi\in C_0^{\infty}(\RR^N)$, such that $\chi(x)=1$ for $|x|\leq 1/2$ and $\chi(x)=0$ for $|x|\geq 1$. Then
$$\lim_{R\to \infty}\lf\|\nabla(v_0^R-v_0)\rg\|_{L^2}+ \lf\|v_1^R-v_1\rg\|_{L^2}=0,\text{ where }v_0^R(x)=\chi\lf(\frac{x}{R}\rg)v_0(x),\; v_1^R(x)=\chi\lf(\frac{x}{R}\rg)v_1(x).$$
Choose $R_{\eps}$ such that for $R\geq R_{\eps}$, $\sqrt{\lf\|\nabla(v_0^R-v_0)\rg\|_{L^2}^2+ \lf\|v_1^R-v_1\rg\|_{L^2}^2}\leq \eps$. Let $R\geq R_{\eps}$ and denote by $v^R_n$ the solution with initial condition $v_{0,n}^R=\lambda_n^{1-N/2}v^R_0(x/\lambda_n)$, $v_{1,n}^R=\lambda_n^{-N/2}v^R_1(x/\lambda_n)$. By conservation of the energy and the scaling of the equation,
$$ \forall n,\quad \sqrt{\lf\|\nabla v^R_n(t_n)-\nabla v_n(t_n)\rg\|_{L^2}^2+ \lf\|\partial_t v^R_n(t_n)-\partial_t v_n(t_n)\rg\|_{L^2}^2}\leq \eps.$$
By the strong Huygens principle, $(v_n^R(t_n),\partial_t v_n^R(t_n))$ is supported in the ring $\{|t_n|-R\lambda_n \leq |x|\leq |t_n|+R\lambda_n\}$. Hence for large $n$ (using Hardy's inequality),
\begin{align*}
&\lf(\int_{\big||t_n|-|x|\big|\geq R\lambda_n} |\nabla v_n(t_n)|^2+\frac{1}{|x|^2}|v_n(t_n)|^2+\lf(\partial_t v_n(t_n)\rg)^2dx\rg)^{1/2}
\\
&\qquad\leq \lf(\int_{\big||t_n|-|x|\big|\geq R\lambda_n} \lf|\nabla v_n^R(t_n)\rg|^2+\frac{1}{|x|^2}|v_n^R(t_n)|^2+\lf(\partial_t v_n^R(t_n)\rg)^2dx\rg)^{1/2}+C\eps=C\eps,
\end{align*}
which concludes the proof of the lemma.
\end{proof}

We next give, in Lemma \ref{L:linear_behavior}, a property of the energy of radial solutions to the linear equation in space dimension $N=3$. In Corollary \ref{C:small_data_behavior'} we deduce a similar property for solutions of the non-linear equation which are sum of small profiles.
\begin{lemma}
\label{L:linear_behavior}
Assume that $N=3$. Let $v$ be a \textbf{radial} solution of \eqref{lin_wave}, $t_0\in \RR$, $0<r_0<r_1$.
Then the following property holds for all $t\geq t_0$ or for all $t\leq t_0$
\begin{multline}
\label{pro_linear}
\int_{r_0+|t-t_0|<r<r_1+|t-t_0|} \left(\partial_r \big(r v(t,x)\big)\right)^2 +r^2(\partial_t v(t,x))^2dr\\
\geq \frac 12\int_{r_0<r<r_1} \left(\partial_r \big(r v(t_0,x)\big)\right)^2+r^2(\partial_t v(t_0,x))^2dr.
\end{multline} 
\end{lemma}
\begin{proof}
We can assume that $t_0=0$. Let $f= rv$, $f_0=f_{\restriction t=0}$, $f_1=\partial_t f_{\restriction t=0}$. Then
\begin{equation}
\label{eq_wave_d1}
\partial^2_t f=\partial^2_r f, \quad t\in \RR,\; r>0.
\end{equation} 
Furthermore, as $v(t)$ is in $\hdot$ for all $t$, by Hardy's inequality in dimension $3$,
$$\int \frac{1}{r^2}\big(f(t,r)\big)^2dr+\int \big(\partial_r f(t,r)\big)^2dr <\infty.$$
By Sobolev embeddings in dimension $1$, for all $t$, $f(t,\cdot)$ is continuous and satisfies the condition $f(t,0)=0$.  By explicit computation we get
\begin{equation*}
f(t,r)=F(t+r)-F(t-r),\quad t\in \RR,\; r>0
\end{equation*} 
where $F$ is defined by 
\begin{equation*}
F(s)=\begin{cases}
\ds\frac 12 f_0(s)+\frac{1}{2}\int_0^{s} f_1(\sigma)\,d\sigma,& s>0\\
\ds-\frac 12 f_0(-s)+\frac{1}{2}\int_0^{-s} f_1(\sigma)\,d\sigma,& s<0.
\end{cases}
\end{equation*} 
Thus, if $t\in \RR$,
\begin{equation}
\label{estim_norm_f}
\int_{r_0+|t|}^{r_1+|t|}(\partial_t f(t,r))^2+(\partial_r f(t,r))^2dr=2\int_{r_0+|t|}^{r_1+|t|} \big(F'(t+r)\big)^2+\big(F'(t-r)\big)^2dr.
\end{equation}
Consequently, if $t>0$,
$$ \int_{r_0+|t|}^{r_1+|t|}(\partial_t f)^2+(\partial_r f)^2\,dr\geq 2\int_{r_0}^{r_1}\big(F'(-r)\big)^2\,dr,$$
and if $t<0$,
$$ \int_{r_0+|t|}^{r_1+|t|}(\partial_t f)^2+(\partial_r f)^2\,dr\geq 2\int_{r_0}^{r_1}\big(F'(r)\big)^2\,dr.$$
By \eqref{estim_norm_f} at $t=0$ we get that the inequality
$$\int_{r_0+|t|}^{r_1+|t|}\big(\partial_t f(t,r)\big)^2+\big(\partial_r f(t,r)\big)^2\,dr \geq \frac{1}{2}\int_{r_0}^{r_1}\big(f_1(r)\big)^2+\big(\partial_r f_0(r)\big)^2dr$$
holds for all $t>0$ or for all $t<0$, hence \eqref{pro_linear}.
%
%
\end{proof}
\begin{corol}
\label{C:small_data_behavior'}
Assume that $N=3$. Let $C_0>0$. Then there exists a constant $\delta_1=\delta_1(C_0)>0$ with the following property. Consider $J>0$,  and let $\{\lambda_{1,n}\}_n$,\ldots,$\{\lambda_{J,n}\}_n$ be sequences of positive numbers such that 
$$ \lambda_{1,n}\ll \ldots\ll \lambda_{J,n} \text{ as }n\to \infty.$$
Consider $J$ \textbf{radial} solutions  $U^1$, \ldots, $U^J$ of \eqref{CP} with initial conditions $(U^j_0,U^j_1)$, $j=1\ldots J$ such that 
$$\forall j\in \{1,\ldots J\},\quad \sqrt{\|\nabla U^j_0\|_{L^2}^2+\|U^j_1\|_{L^2}^2}=\eta_j\leq\delta_1.$$
Consider a sequence $w_n$ of solutions of the linear wave equation \eqref{lin_wave} such that
$$ \forall j\in \{1,\ldots, J\}, \quad 
 \lf(\lambda_j^{N/2-1}w_{0,n}(\lambda_jx),\lambda_j^{N/2}w_{1,n}(\lambda_jx)\rg)\xrightharpoonup[n\to \infty]{} 0 \text{ weakly in }\hdot\times L^2,
$$
Let $\eta= \sqrt{\sum_{j=1}^J\eta_j^2}$ and assume that $\eta\leq C_0.$
Let
$$ U_n(t,x)=\sum_{j=1}^J\frac{1}{\lambda_{j,n}^{N/2-1}}U^j\lf(\frac{t}{\lambda_{j,n}},\frac{x}{\lambda_{j,n}}\rg)+w_n(t,x).$$
Then there exists $r_1>0$ such that for large $n$, the inequality
\begin{equation}
\label{cor_NL'}
 \sqrt{\int_{r_{1}\lambda_{1,n}+|t|<|x|} \big|\nabla U_n(t,x)\big|^2+\big(\partial_t U_n(t,x)\big)^2dx}\geq \frac{\eta}{4}
\end{equation} 
holds for all $t>0$ or for all $t<0$.
\end{corol}
\begin{proof}
Denote by 
$$ U_{0,n}(x)=U_{n}(0,x),\quad U_{1,n}(x)=\partial_t U_{n}(0,x).$$
Let $U^j_{\lin}$ be the solution of \eqref{lin_wave} with initial conditions $(U^j_0,U^j_1)$, $j=1\ldots J$ and 
$$U_{\lin,n}(t,x)=\sum_{j=1}^J\frac{1}{\lambda_{j,n}^{N/2-1}}U^j_{\lin}\lf(\frac{t}{\lambda_{j,n}},\frac{x}{\lambda_{j,n}}\rg)+w_n(t,x)$$

\EMPH{Step 1}
We first show that if $\delta_1=\delta_1(C_0)$ is chosen small enough, then
\begin{equation}
\label{small_NL'}
\sup_{t\in \RR}\sqrt{\left\|U_n(t)-U_{\lin,n}(t)\right\|_{\hdot}^2+\left\|\partial_tU_n(t)-\partial_t U_{\lin,n}(t)\right\|_{L^2}^2}\leq \frac{\eta}{4},
\end{equation} 
Indeed by \eqref{almost_linear}, if $\delta_1^3\leq \frac{1}{C C_0}$, for some large constant $C$, then $\sqrt{\|\nabla U_0^j\|_{L^2}^2+\|U_1^j\|_{L^2}^2}=\eta_j\leq \delta_1$ implies
\begin{equation*}
\sup_{t\in \RR}\sqrt{\left\|U^j(t)-U^j_{\lin}(t)\right\|_{\hdot}^2+\left\|\partial_tU^j(t)-\partial_t U^j_{\lin}(t)\right\|_{L^2}^2}\leq C\eta_j^5\leq \frac{\eta_j^2}{4C_0}.
\end{equation*} 
By the triangle inequality and the fact that $\eta\leq C_0$, 
\begin{equation*}
\sup_{t\in \RR}\left(\sqrt{\left\|U_n(t)-U_{\lin,n}(t)\right\|_{\hdot}^2+\left\|\partial_tU_n(t)-\partial_t U_{\lin,n}(t)\right\|_{L^2}^2}\right)\leq \frac{\eta^2}{4C_0}\leq \frac{\eta}{4}.
\end{equation*} 
Hence \eqref{small_NL'}.

\EMPH{Step 2}
We next show that there exists $r_1>0$  such that
\begin{equation}
\label{minor_U0n}
\liminf_{n\to \infty}\int_{r_1 \lambda_{1,n}}^{+\infty}\left(\partial_r \big(r U_{0,n}(r)\big)\right)^2+\big(r U_{1,n}(r)\big)^2dr\geq \frac{\eta^2}{2}.
\end{equation} 
Indeed, if $f\in \hdot$ is a radial function and $0<R_0<R_1$,
\begin{multline*}
\int_{R_0}^{R_1} \left(\partial_r \big(r f(r)\big)\right)^2dr
= \int_{R_0}^{R_1} f^2+r^2(\partial_r f)^2+2rf\partial_r f dr
=\int_{R_0}^{R_1}f^2+r^2(\partial_r f)^2+r\partial_r \left(f^2\right)dr\\
=\int_{R_0}^{R_1}r^2(\partial_r f)^2dr+R_1f^2(R_1)-R_0f^2(R_0).
\end{multline*}
By Hardy's inequality, $\int f^2(t,r)dr<\infty$, which implies that there exist sequences $R_n\to +\infty$ and $\tilde{R}_n\to 0$ such that $R_nf^2(R_n)\to 0$ and $\tilde{R}_nf^2(\tilde{R}_n)\to 0$. Letting $R_1=R_n$ and $n\to +\infty$, we get
\begin{equation}
\label{int_f_1}
\int_{R_0}^{+\infty} \left(\partial_r \big(r f(r)\big)\right)^2dr=\int_{|x|\geq R_0} |\nabla f|^2dx-R_0f^2(R_0)\leq \int_{|x|\geq R_0} |\nabla f|^2dx.
\end{equation}
Letting $R_0=\tilde{R}_n$ and $n\to +\infty$ we get
\begin{equation}
\label{int_f_2}
\int_{0}^{+\infty} \left(\partial_r \big(r f(r)\big)\right)^2dr=\int_{\RR^3} |\nabla f|^2dx.
\end{equation}
By \eqref{int_f_2}, there exists $r_1>0$ such that
\begin{equation}
\label{minor_U1}
\int_{r_1}^{+\infty}\big(\partial_r(rU_0^1(r))\big)^2+\big(rU_1^1(r)\big)^2 dr\geq \frac{\eta_1^2}{2}.
\end{equation}

Let $g^j=\partial_r(rU^j_0(r))\in L^2(dr)$. Then
\begin{multline}
\label{scalar_p_Uj}
A_n^{j,k}:=\lf|\int_{r_1\lambda_{1,n}}^{+\infty}\partial_r\left(\frac{r}{\lambda_{j,n}^{1/2}}U^j_0\lf(\frac{r}{\lambda_{j,n}}\rg)\rg)\partial_r\left(\frac{r}{\lambda_{k,n}^{1/2}}U^k_0\lf(\frac{r}{\lambda_{k,n}}\rg)\rg)dr\rg|\\=\lf|\int _{r_1\lambda_{1,n}}^{+\infty} \frac{1}{\lambda_{j,n}^{1/2}}g_j\left(\frac{r}{\lambda_{j,n}}\rg)\frac{1}{\lambda_{k,n}^{1/2}}g_k\left(\frac{r}{\lambda_{k,n}}\rg)dr\rg|
=\lf|\int _{\frac{r_1\lambda_{1,n}}{\lambda_{j,n}}}^{+\infty} g_j\left(\rho\rg)\frac{\lambda_{j,n}^{1/2}}{\lambda_{k,n}^{1/2}}g_k\left(\frac{\lambda_{j,n}}{\lambda_{k,n}}\rho\rg)d\rho\rg|.
\end{multline}
Letting $j=k$ in \eqref{scalar_p_Uj} we get that if $j>1$, $A_n^{j,j}\to \int_0^{\infty} g_j^2(\rho)d\rho$ as $n\to \infty$. Furthermore if $1\leq j<k$ we obtain, using that $\lambda_{k,n}/\lambda_{j,n}\to +\infty$, that for all $\eps>0$, 
\begin{equation*}
\left|A_n^{j,k}\right|\leq C_j\sqrt{\int_0^{R_{\eps}} \frac{\lambda_{j,n}}{\lambda_{k,n}}g_k^2\left(\frac{\lambda_{j,n}}{\lambda_{k,n}}\rho\rg)d\rho}+C_k\sqrt{\int_{R_{\eps}}^{+\infty} g_j^2\left(\rho\rg)d\rho}\leq o_n(1)+\eps.
\end{equation*}
and hence $A_n^{j,k}\to 0$ as $n\to\infty$.
Similarly, noting that $h_{n}=\partial_r(rw_{0,n})=w_{0,n}+r\partial_r w_{0,n}$ is such that $\lambda_{j,n}^{1/2}h_{n}(\lambda_{j,n}\cdot)$ converges weakly to $0$ in $L^2(dr)$ we get
\begin{multline*}
\int_{r_1\lambda_{1,n}}^{+\infty}\partial_r\left(\frac{r}{\lambda_{j,n}^{1/2}}U^j_0\lf(\frac{r}{\lambda_{j,n}}\rg)\rg)\partial_r(rw_{0,n})dr
=\int_{r_1\frac{\lambda_{1,n}}{\lambda_{j,n}}}^{+\infty}g_j(\rho)\lambda_{j,n}^{1/2}h_n\lf(\lambda_{j,n}\rho\rg) d\rho\\
=\begin{cases}
\int_{r_1}^{+\infty}g_1(\rho)\lambda_{1,n}^{1/2}h_n\lf(\lambda_{1,n}\rho\rg) d\rho+o_n(1)\text{ if }j=1\\
\int_{0}^{+\infty}g_j(\rho)\lambda_{j,n}^{1/2}h_n\lf(\lambda_{j,n}\rho\rg) d\rho+o_n(1)\text{ if }j>1,
 \end{cases}
\end{multline*}
which tends to $0$ as $n\to \infty$. Using similar estimates on $U_{1,n}$ and $w_{1,n}$ and combining with \eqref{minor_U1} we get \eqref{minor_U0n}. 

\EMPH{Step 3: end of the proof}
In view of Step 2 and Lemma \ref{L:linear_behavior}, if $n$ is large, then the following holds for all $t>0$ or for  all $t<0$:
\begin{equation*}
 \int_{r_1\lambda_{1,n}+|t|}^{+\infty} (\partial_r (rU_{n,\lin}))^2+(\partial_t (rU_{n,\lin}))^2dr\geq \frac{\eta^2}{4}.
\end{equation*} 
By \eqref{int_f_1}, we get that for all $t>0$ or for all $t<0$,
\begin{equation*}
 \int_{|x|\geq r_1\lambda_{1,n}+|t|} |\nabla U_{n,\lin}|^2+|\partial_t U_{n,\lin}|^2dx\geq \frac{\eta^2}{4}.
\end{equation*}
By Step 1 and the triangle inequality, 
\begin{equation*}
\sqrt{\int_{|x|\geq r_1\lambda_{1,n}+|t|} |\nabla U_{n}|^2+|\partial_t U_{n}|^2dx}\geq \frac{\eta}{2}-\frac{\eta}{4}=\frac{\eta}{4}, 
\end{equation*}
which concludes the proof.
\end{proof}

\subsection{No small exterior profile}
\label{SS:No_exterior}
Before stating the main result of this section, we introduce some notations. In all the sequel we assume $N=3$. Let $\tau_n\to 1^-$ and consider a profile decomposition of $(a(\tau_n),\partial_t a(\tau_n))$ with profiles $\{U^j_{\lin}\}$ and parameters $\{\lambda_{j,n},t_{j,n}\}$. We will consider as usual the nonlinear profiles $\{U^j\}$ associated to $\{U^j_{\lin}\}$, $-t_{j,n}/\lambda_{j,n}$, and will write, for the sake of simplicity
\begin{equation*}
 U^j_{\lin,n}(t,x)=\frac{1}{\lambda_{j,n}^{1/2}}U^j_{\lin}\left(\frac{t-t_{j,n}}{\lambda_{j,n}},\frac{x}{\lambda_{j,n}}\right),\quad U^j_n(t,x)=\frac{1}{\lambda_{j,n}^{3/2}}U^j\left(\frac{t-t_{j,n}}{\lambda_{j,n}},\frac{x}{\lambda_{j,n}}\right).
\end{equation*}
The second expression is defined as long as $(t-t_{j,n})/\lambda_{j,n}$ is in $(T_-(U^j),T_+(U^j))$. We will also write
$$ U^j_{0,n}= U^j_{\lin,n}(0,x),\quad U^j_{1,n}=\left(\partial_tU^{j}_{\lin,n}\right)(0,x).    $$
Let $j\in \NN^*$. Extracting subsequences and time-translating the profiles if necessary we can assume that 
\begin{equation}
\label{Hprofiles1}
\forall n,\quad t_{j,n}=0\text{ or }\lim_{n\to +\infty}\frac{t_{j,n}}{\lambda_{j,n}}\in \{-\infty,+\infty\}.
\end{equation}
We will denote by 
$$ \rho_{j,n}=|t_{j,n}|\quad \text{if} \quad\frac{|t_{j,n}|}{\lambda_{j,n}}\to \infty$$
and
$$ \rho_{j,n}=\lambda_{j,n}\quad \text{ if }\quad t_{j,n}=0.$$
According to Lemma \ref{L:lin_odd} the sequence $(U^j_{0,n},U^j_{1,n})_n$ is localized, for large $n$, around $|x|\approx \rho_{j,n}$. 

Reordering the profiles and extracting subsequences, we can find a $J_0 \in \NN$ such that (here $\delta_1(C_0)$ is given by Corollary \ref{C:small_data_behavior'}, and $C_0$ is the constant in assumption \eqref{bound_nabla'}):
\begin{align}
\label{Hprofiles2}
j>J_0 &\iff \left(E(U^j_0,U^j_1)\leq \frac{1}{N}\left(\delta_1(C_0)\right)^2\text{ and }\big\|\nabla U_0^j\big\|_{L^2}<\big\|\nabla W\big\|_{L^2}\right)\\
\notag &\text{ or }\Big(\lim_{n\to \infty} \frac{-t_{j,n}}{\lambda_{j,n}}\in \{\pm \infty\} \text{ and }E(U^j_0,U^j_1)<E(W,0)
\Big)\\
\label{Hprofiles3}
\rho_{J_0,n}&\lesssim \rho_{J_0-1,n}\lesssim \ldots\lesssim \rho_{1,n}.
\end{align}
In particular if $j>J_0$ and $t_{j,n}=0$ for all $n$, we have by Claim \ref{C:variationnal} that $\big\|\nabla U^j_0\big\|_{L^2}^2+\big\|U^j_1\big\|^2_{L^2}\leq \lf(\delta_1(C_0)\rg)^2$.

In this section we show:
\begin{prop}
\label{P:ext_profile}
Under the above assumptions, 
\begin{equation}
\label{no_E_outside'}
\lim_{R\rightarrow +\infty} \limsup_{n\rightarrow+\infty} \int_{|x|\geq R\rho_{1,n}} \left(|\nabla a(\tau_n)|^2+(\partial_t a(\tau_n))^2\right)dx=0.
\end{equation} 
\end{prop}
\begin{proof}
We argue by contradiction. If \eqref{no_E_outside'} does not hold, there exists $\eps_0>0$ and a sequence $\overline{\rho}_n$ such that
\begin{equation}
\label{exterior}
\int_{|x|\geq \overline{\rho}_n} \left(|\nabla a(\tau_n)|^2+(\partial_t a(\tau_n))^2\right)dx\geq \eps_0,\quad \lim_{n\to +\infty}\frac{\overline{\rho}_n}{\rho_{1,n}}=+\infty.
\end{equation} 
Since $\supp a(\tau_n)\subset \big\{|x|\leq 1-\tau_n\big\}$, we have that $\overline{\rho}_n\leq 1-\tau_n$.
Moreover, by Claim \ref{C:faraway}, we get, extracting subsequences in $n$, a sequence $\{\trho_n\}_n$ such that
\begin{equation}
\label{trho_1}
\rho_{1,n}\ll \trho_n \ll \overline{\rho}_n.
\end{equation} 
and 
\begin{equation}
\label{trho_j}
\forall j,\quad \trho_n\ll \rho_{j,n}\text{ or }\rho_{j,n}\ll \trho_n.
\end{equation}
Let $\chi\in C^{\infty}(\RR^N)$, such that $\chi(x)=1$ if $|x|\geq 2$ and $\chi(x)=0$ if $|x|\leq 1$. Then
\begin{align}
\label{dev_chi0}
\chi\lf(\frac{x}{\tilde{\rho}_n}\rg)u(\tau_n,x)&=\chi\lf(\frac{x}{\tilde{\rho}_{n}}\rg)v(\tau_n,x)+\sum_{j=1}^J \chi\lf(\frac{x}{\tilde{\rho}_{n}}\rg)U_{0,n}^j+\chi\lf(\frac{x}{\tilde{\rho}_{n}}\rg)w_{0,n}^{J}\\
\label{dev_chi1}
\chi\lf(\frac{x}{\trho_{n}}\rg)\partial_t u(\tau_n,x)&=\chi\lf(\frac{x}{\trho_{n}}\rg)\partial_t v(\tau_n,x)+\sum_{j=1}^J \chi\lf(\frac{x}{\trho_{n}}\rg)U_{1,n}^j+\chi\lf(\frac{x}{\trho_{n}}\rg)w_{1,n}^{J}.
\end{align}
\begin{claim}
\label{C:exterior_only}
If $\rho_{j,n}\ll\tilde{\rho}_n$ then
$$\lim_{n\to +\infty}\lf\|\chi\lf(\frac{x}{\tilde{\rho}_{n}}\rg)U_{0,n}^j\rg\|_{\hdot}+\lf\|\chi\lf(\frac{x}{\tilde{\rho}_{n}}\rg)U_{1,n}^j\rg\|_{L^2}=0.
$$
If $\tilde{\rho}_n\ll \rho_{j,n}$ then
$$\lim_{n\to +\infty}\lf\|\chi\lf(\frac{x}{\tilde{\rho}_{n}}\rg)U_{0,n}^j-U_{0,n}^j\rg\|_{\hdot}+\lf\|\chi\lf(\frac{x}{\tilde{\rho}_{n}}\rg)U_{1,n}^j-U_{1,n}^j\rg\|_{L^2}=0.$$
\end{claim}
\begin{proof}
Indeed by Lemma \ref{L:lin_odd}, 
$$\lim_{R\to +\infty}\liminf_{n\to \infty} \int_{\frac{1}{R}\rho_{j,n}\leq |x|\leq R\rho_{j,n}}\lf|\nabla U_{0,n}^j\rg|^2+\lf|U_{1,n}^j\rg|^2=\int_{\RR^N}\lf|\nabla U_{0,n}^j\rg|^2+\lf|U_{1,n}^j\rg|^2.$$
In the case $\rho_{j,n}\ll\tilde{\rho}_n$, choose $\eps>0$ and $R=R({\eps})$ such that
$$ \limsup_{n\to +\infty}\int_{R\rho_{j,n}\leq |x|}\lf|\nabla U_{0,n}^j\rg|^2+\lf|U_{1,n}^j\rg|^2\leq \eps.$$
As $R\rho_{j,n}\ll \tilde{\rho}_n$, we get that for large $n$,
$$ \lf\|\chi\lf(\frac{x}{\tilde{\rho}_{n}}\rg)U_{0,n}^j\rg\|_{\hdot}^2+\lf\|\chi\lf(\frac{x}{\tilde{\rho}_{n}}\rg)U_{1,n}^j\rg\|_{L^2}^2 \leq \int_{R\rho_{j,n}\leq |x|}\lf|\nabla U_{0,n}^j\rg|^2+\lf|U_{1,n}^j\rg|^2\leq \eps,$$
which shows the first estimate of the claim. The proof of the second one is similar and we skip it.
\end{proof}
Let us denote by $\JJJ_{\ext}$ the set of indexes $j$ such that $\tilde{\rho}_n\ll \rho_{j,n}$. Note that for $j\in \JJJ_{\ext}$, $j>J_0$ and thus 
\begin{multline*}
j\in \JJJ_{\ext}\Longrightarrow \bigg(E(U^j_0,U^j_1)\leq \frac{\left(\delta_1(C_0)\right)^2}{N}\text{ and }  \big\|\nabla U^j_0\big\|_{L^2}<\big\|\nabla W\|_{L^2}\bigg)\\
\text{ or } \bigg(\lim_{n\to\infty}\frac{-t_{j,n}}{\lambda_{j,n}}=\pm\infty\text{ and }E(U^j_0,U^j_1)<E(W,0)\bigg)
\end{multline*}
so the corresponding nonlinear profile $U^j$ is globally defined and scatters in both time directions. In view of Claim \ref{C:exterior_only}, we rewrite \eqref{dev_chi0}, \eqref{dev_chi1} as
\begin{align}
\label{dev_chi0'}
\chi\lf(\frac{x}{\tilde{\rho}_n}\rg)u(\tau_n,x)&=v(\tau_n,x)+\sum_{\substack{j \in \JJJ_{\ext}\\j\leq J}} U_{0,n}^j(x)+\tilde{w}_{0,n}^J(x)\\
\label{dev_chi1'}
\chi\lf(\frac{x}{\trho_{n}}\rg)\partial_t u(\tau_n,x)&=\partial_t v(\tau_n,x)+\sum_{\substack{j\in \JJJ_{\ext}\\ j\leq J}} U_{1,n}^j(x)+\tw_{1,n}^J(x).
\end{align}
where
\begin{equation*}
\tw_{0,n}^J=\chi\lf(\frac{x}{\trho_{n}}\rg)w_{0,n}^{J}+o_n(1)\text{ in }\hdot,\quad \tw_{1,n}^J=\chi\lf(\frac{x}{\trho_{n}}\rg)w_{1,n}^{J}+o_n(1)\text{ in }L^2.
\end{equation*}
By Claim \ref{C:dispersive_tw}, 
\begin{equation}
\label{small_tw}
\lim_{n\rightarrow+\infty}\limsup_{J\rightarrow+\infty} \left\|\tw_{n}^J\right\|_{S(\RR)}=0.
\end{equation}
Indeed if \eqref{small_tw} does not hold, one can find, in view of \eqref{small_w}, sequences $\{n_k\}_k$, $\{J_k\}_k$ and $\eps>0$ such that
\begin{equation*}
\forall k,\;\left\|\tw_{n_k}^{J_k}\right\|_{S(\RR)}\geq \eps\text{ and } \lim_{k\to\infty} \left\|w_{n_k}^{J_k}\right\|_{S(\RR)}=0,
\end{equation*} 
a contradiction with Claim \ref{C:dispersive_tw}.

By \eqref{small_tw}, the decomposition \eqref{dev_chi0'}, \eqref{dev_chi1'} is a profile decomposition of the sequence 
$$\chi\lf(\frac{x}{\tilde{\rho}_n}\rg)\lf(u(\tau_n,x),\partial_t u(\tau_n,x)\rg).$$ 
Denote by $\tu_n$ the solution of \eqref{CP} such that
$$ \tu_{n\restriction t=\tau_n}=\chi\lf(\frac{x}{\tilde{\rho}_n}\rg)u(\tau_n,x) ,\quad \partial_t\tu_{n\restriction t=\tau_n}=\chi\lf(\frac{x}{\tilde{\rho}_n}\rg)\partial_t u(\tau_n,x).$$
Using that all the non-linear solutions $U^j$, $j\in \JJJ_{\ext}$ are globally defined and scatter, we get by Proposition \ref{P:lin_NL} that $\tu_n$ is globally defined for large $n$ and
\begin{equation}
\label{approx_tu}
\tu_n(\tau_n+t,x)=v(\tau_n+t,x)+\sum_{\substack{j\in \JJJ_{\ext}\\j\leq J}} U^j_n(t,x)+\tw_{n}^J(t,x)+r_n^J(t,x),
\end{equation} 
where $r_n^J$ satisfies \eqref{cond_rJn}. By the definition of $\tu_n$,
\begin{equation}
\label{eq_tau_n}
\tu_n(\tau_n,x)=u(\tau_n,x),\quad\partial_t \tu_n(\tau_n,x)=\partial_t u(\tau_n,x)\text{ if }|x|\geq 2\trho_n,
\end{equation} 
By finite speed of propagation, as long as $0\leq \tau_n+t<1$, we have
\begin{equation}
\label{eq_tau_n+t}
\tu_n(\tau_n+t,x)=u(\tau_n+t,x),\quad\partial_t \tu_n(\tau_n+t,x)=\partial_t u(\tau_n+t,x)\text{ if }|x|\geq 2\trho_n+|t|.
\end{equation}
The key point of the proof is the following claim:
\begin{claim}
\label{C:exterior}
The set $\JJJ_{\ext}$ is empty.
\end{claim}
\begin{proof}
The proof takes several steps.

\EMPH{Step 1. No profile dispersing backward in time}
Let $k\in \JJJ_{\ext}$. We first show by contradiction that we cannot have $\frac{-t_{k,n}}{\lambda_{k,n}}\to -\infty$. Let us assume that $\frac{-t_{k,n}}{\lambda_{k,n}}\to -\infty$. Then $\rho_{k,n}=|t_{k,n}|$. Furtermore $U^k$ scatters backward in time. 
As a consequence, by Lemma \ref{L:lin_odd}, if $M$ is large enough, there exists $\eps_k>0$ such that  for all large $n$,
\begin{equation*}
\int_{|x|\geq t_{k,n}+\tau_n-M\lambda_{k,n}} \lf|\nabla U^{k}_n\lf(-\tau_n,x\rg)\rg|^2+\lf|\partial_t U^{k}_n\lf(-\tau_n,x\rg)\rg|^2 dx\geq \eps_k.
\end{equation*} 
As $k\in \JJJ_{\ext}$, we know that $t_{k,n}=\rho_{k,n}\gg \trho_n $. Furthermore $\lambda_{k,n}=o(|t_{k,n}|)$. Thus for large $n$, $t_{k,n}+\tau_n-M\lambda_{k,n}\gg 2\trho_n+\tau_n$, and the preceding inequality implies
\begin{equation}
\label{channel}
\int_{|x|\geq 2\trho_n+\tau_n} \lf|\nabla U^{k}_n\lf(-\tau_n,x\rg)\rg|^2+\lf|\partial_t U^{k}_n\lf(-\tau_n,x\rg)\rg|^2 dx\geq \eps_k.
\end{equation} 
Using again that $U^k$ scatters backward in time and that $\trho_n\ll t_{k,n}$, we get by Lemma \ref{L:lin_odd} 
\begin{equation}
\label{channel_bis}
\int_{|x|\leq 2\trho_n+\tau_n} \lf|\nabla U^{k}_n\lf(-\tau_n,x\rg)\rg|^2+\lf|\partial_{t} U^{k}_n\lf(-\tau_n,x\rg)\rg|^2 dx=o_n(1).
 \end{equation} 
Let $j\in \JJJ_{\ext}\setminus\{k\}$. Then $U^j$ scatters in both time directions, and there exists a solution $V^j_{\lin}$ of the linear wave equation such that
$$ \lim_{t\to -\infty} \big\|V^j_{\lin}(t)-U^j(t)\big\|_{\hdot}+\big\|\partial_tV^j_{\lin}(t)-\partial_t U^j(t)\big\|_{L^2}=0.$$
Noting 
$V^j_{\lin,n}(t,x)=\frac{1}{\lambda_{j,n}^{1/2}}V^j_{\lin}\left(\frac{t-t_{j,n}}{\lambda_{j,n}},\frac{x}{\lambda_{j,n}}\right),$ we get, by the conservation of the energy for the linear wave equation and \eqref{channel_bis}
\begin{multline}
\label{ortho_channel1}
\int_{|x|\geq \tau_n+2\trho_n}\nabla_{t,x} U^{j}_n\lf(-\tau_n\rg)\cdot \nabla_{t,x} U^{k}_n\lf(-\tau_n\rg)dx
=\int_{\RR^3}\nabla_{t,x} U^{j}_n\lf(-\tau_n\rg)\cdot \nabla_{t,x} U^{k}_n\lf(-\tau_n\rg)dx
+o_n(1)\\
=\int_{\RR^3}\nabla_{t,x} V^{j}_{\lin,n}\lf(0\rg)\cdot \nabla_{t,x} U^{k}_{\lin,n}\lf(0\rg)dx+o_n(1)=o_n(1).
\end{multline}
where we used the orthogonality of the parameters $(\lambda_{j,n},t_{j,n})$ and $(\lambda_{k,n},t_{k,n})$. Similarly, if $J>k$,
\begin{multline}
\label{ortho_channel2}
\int_{|x|\geq \tau_n+2\trho_n}\nabla_{t,x} U^{k}_n\lf(-\tau_n,x\rg)\cdot \nabla_{t,x} \tilde{w}^{J}_n\lf(-\tau_n,x\rg)dx\\
=\int \nabla_{t,x} U^{k}_n\lf(-\tau_n,x\rg)\cdot \nabla_{t,x} \tilde{w}^{J}_n\lf(-\tau_n,x\rg)dx+o_n(1)\\
=\int \nabla_{t,x} V^{k}_{\lin,n}\lf(-\tau_n,x\rg)\cdot \nabla_{t,x} \tilde{w}^{J}_n\lf(-\tau_n,x\rg)dx+o_n(1)\\
=\int \nabla_{t,x} V^{k}_{\lin,n}\lf(-t_{k,n},x\rg)\cdot \nabla_{t,x} \tilde{w}^{J}_n\lf(0,x\rg)dx+o_n(1)=o_n(1).
\end{multline}
At the last line, we used the conservation of the energy, and the fact that by \eqref{small_tw}, the $\tilde{w}^J_n$ are the remainders of the profile decomposition \eqref{dev_chi0'}, \eqref{dev_chi1'} and thus by \eqref{weak_CV_wJ}, 
$$ \lambda_{k,n}^{N/2}\nabla_{t,x}\tilde{w}^J_n\lf(t_{k,n},\lambda_{k,n}x\rg)\xrightharpoonup[n\to \infty]{}0\text{ in }(L^2)^{N+1}.$$
%
Combining \eqref{approx_tu} with $t=-\tau_n$, \eqref{eq_tau_n}, \eqref{channel}, \eqref{ortho_channel1} and \eqref{ortho_channel2} we get, if $n$ is large enough,
$$\int_{|x|\geq \tau_n+2\trho_n}\lf[\lf|\nabla_{t,x} u(0,x)\rg|^2-\lf|\nabla_{t,x} v(0,x)\rg|^2\rg]\geq \frac{\eps_k}{2}. $$
Using that the function $x\mapsto \lf[\lf|\nabla_{t,x} u(0,x)\rg|^2-\lf|\nabla_{t,x} v(0,x)\rg|^2\rg]$ is supported in the set $\{|x|\leq 1\}$, we get
$$\int_{2\trho_n+|\tau_n|\leq |x|\leq 1}\Big|\lf|\nabla_{t,x} u(0,x)\rg|^2-\lf|\nabla_{t,x} v(0,x)\rg|^2\Big|\geq \frac{\eps_k}{2}. $$
Letting $n\to \infty$ we have $2\trho_n+|\tau_n|\to 1$ which yields a contradiction.

\EMPH{Step 2. No profile dispersing forward in time}
We next show by contradiction that if $k\in \JJJ_{\ext}$ we cannot have $\lim_{n\to +\infty}\frac{-t_{k,n}}{\lambda_{k,n}}\to +\infty$. Let $\sigma_n=(1-\tau_n)/2$. Using that $U^k$ scatters forward in time, we get by Lemma \ref{L:lin_odd} that if $M$ is large enough, there exists $\eps_k>0$ such that  for all large $n$,
\begin{equation}
\label{channel2}
\int_{|x|\geq |t_{k,n}|+\sigma_n-M\lambda_{k,n}} \lf|\nabla_{t,x} U^{k}_n\lf(\sigma_n,x\rg)\rg|^2dx\geq \eps_k.
\end{equation} 
By Lemma \ref{L:lin_odd}, we also have (using that $\lambda_{k,n}\ll |t_{k,n}|$), 
$$ \lim_{n\to \infty}\int_{|x|\leq \sigma_n} \lf|\nabla_{t,x} U^{k}_n\lf(\sigma_n,x\rg)\rg|^2dx=0,$$
from which we can deduce the analogues of the ortogonality conditions \eqref{ortho_channel1} and \eqref{ortho_channel2} with $\tau_n+2\trho_n$ replaced by $\sigma_n$.
As in the preceding case, using \eqref{approx_tu} with $t=\sigma_n$  we deduce from \eqref{channel2} that for large $n$,
$$\int_{|x|\geq \sigma_n}\lf[\lf|\nabla_{t,x} u\lf(\frac{1+\tau_n}{2},x\rg)\rg|^2-\lf|\nabla_{t,x} v\lf(\frac{1+\tau_n}{2},x\rg)\rg|^2\rg]\geq \frac{\eps_k}{2}. $$
As $1-\frac{1+\tau_n}{2}=\sigma_n$, this contradicts the fact that on the support of $u-v$, $|x|\leq 1-t$.

\EMPH{Step 3. No compact profile} In this step we conclude the proof, showing that $\JJJ_{\ext}$ is empty. 
According to Steps 1 and 2, for all $j\in\JJJ_{\ext}$, and all $n$, $t_{j,n}=0$, and  we can rewrite \eqref{approx_tu} as
\begin{equation}
\label{approx_tu'}
\tu_n(\tau_n+t)=v(\tau_n+t,x)+\sum_{\substack{j\in \JJJ_{\ext}\\j\leq J}} \frac{1}{\lambda_{j,n}^{1/2}}U^j\lf(\frac{t}{\lambda_{j,n}},\frac{x}{\lambda_{j,n}}\rg)+\tw_{n}^J(t,x)+r_n^J(t,x),
\end{equation} 
Furthermore, we know that for $j\in \JJJ_{\ext}$, $j>J_0$ and thus by the definition of $J_0$ we have 
\begin{equation}
\label{smallU0j}
\sqrt{\big\|\nabla U_0^j\big\|^2+\big\|U_1^j\big\|^2}\leq \delta_1(C_0).
\end{equation}

Assume that $\JJJ_{\ext}$ is not empty. Then by assumption \eqref{bound_nabla'} for large $J$,
$$ 0<\eta^2=\sum_{\substack{j\in \JJJ_{\ext}\\ j\leq J}} \lf\|\nabla U^j_0\rg\|^2_{L^2}+\lf\|U^j_1\rg\|^2_{L^2}\leq C_0^2.$$
Choose $J$ such that 
\begin{equation}
\label{small_remainder}
 \sup_{t\in \RR} \sqrt{\left\|\nabla r_n^J(t)\right\|^2_{L^2}+\lf\|\partial_t r_n^J(t)\rg\|_{L^2}}\leq \frac{\eta}{8}.
\end{equation} 
Let $k\in \JJJ_{\ext}$, such that $k\leq J$ and
$$ \lambda_k=\inf_{\substack{j\in \JJJ_{\ext}\\ j=1\ldots J}} \lambda_j.$$
By \eqref{smallU0j}, we can use Corollary \ref{C:small_data_behavior'}, which implies that there exists $r_0>0$ such that the following occurs for all $t\in [-\tau_n,0)$ or for all $t\in (0,1-\tau_n)$,
$$\int_{|x|\geq \lambda_{k,n}r_0+|t|}|\nabla U_n(t,x)|^2+|\partial_t U_n(t,x)|^2dx\geq \frac{\eta^2}{16},$$
where
$$ U_n=\sum_{j\in \JJJ_{\ext}} \frac{1}{\lambda_{j,n}^{1/2}}U^j\left(\frac{t}{\lambda_{j,n}},\frac{x}{\lambda_{j,n}}\right)+\tilde{w}^J_n(t,x).$$
And thus by \eqref{small_remainder}, for all $t\in [-\tau_n,0)$ or for all $t\in (0,1-\tau_n)$,
\begin{equation}
\label{channel3}
\int_{|x|\geq \lambda_{k,n}r_0+|t|} \lf|\nabla(\tu_n(\tau_n+t)-v(\tau_n+t))\rg|^2+\lf|\partial_t\tu_n(\tau_n+t)-\partial_tv(\tau_n+t)\rg|^2\geq \frac{\eta^2}{64}.
\end{equation} 
First assume that it holds for all $t\in (0,1-\tau_n)$. 
Letting $t_n=\frac{1-\tau_n-r_0\lambda_{k,n}}{2}$ in \eqref{channel3}, we obtain that for large $n$,
\begin{equation*}
\int_{|x|\geq 1-\tau_n-t_n} \lf|\nabla(\tu_n(\tau_n+t_n)-v(\tau_n+t_n))\rg|^2+\lf|\partial_t\tu_n(\tau_n+t_n)-\partial_tv(\tau_n+t_n)\rg|^2\geq \frac{\eta^2}{64}.
\end{equation*}
Furthermore, as $k\in \JJJ_{\ext}$, we have $2\trho_n+t_n\leq 1-\tau_n-t_n$, and thus by \eqref{eq_tau_n+t}, $|x|\leq 1-\tau_n-t_n$ on the support of $u(\tau_n+t_n,\cdot)-v(\tau_n+t_n,\cdot)$ a contradiction.

It remains to treat the case when \eqref{channel3} holds for all $t\in[-\tau_n,0)$. Then \eqref{channel3} with $t=-\tau_n$ yields
\begin{equation*}
\int_{|x|\geq \lambda_{k,n}r_0+\tau_n} \lf|\nabla(\tu_n(0)-v(0))\rg|^2+\lf|\partial_t\tu_n(0)-\partial_tv(0))\rg|^2\geq \frac{\eta^2}{64},
\end{equation*} 
which is again a contradiction, recalling that $\lf(\tilde{u}_n(0,x),\partial_t \tilde{u}_n(0,x)\rg)$ and $\lf(u(0,x),\partial_t u(0,x)\rg)$ coincide for $|x|\geq \tau_n+2\trho_n$, and thus for $|x|\geq \tau_n+\lambda_{k,n}r_0$ for large $n$. The proof of Claim \ref{C:exterior} is complete.
\end{proof}
To finish the proof of Proposition \ref{P:ext_profile}, we must show that if  $\overline{\rho}_n$ is as in \eqref{exterior}, and $J$ is large, then
\begin{equation}
\label{no_wn_outside}
\lim_{n\to\infty}
\int_{|x|\geq \overline{\rho}_n} \left(|\nabla w_{0,n}^J|^2+(w_{1,n}^J)^2\right)dx=0.
\end{equation}

 We will use that $w_n^J$ is a radial solution of the linear wave equation \eqref{lin_wave}. By \eqref{int_f_1}, we have
\begin{equation}
\label{IPPwn}
\int_{\overline{\rho}_n}^{+\infty} (\partial_r (rw_{0,n}^J))^2dr=\int_{|x|\geq \rhob_n}\lf|\nabla w_n^J\rg|^2dx-\rhob_nw_{0,n}^J(\rhob_n)^2
\end{equation} 
By the construction of the profile decomposition (see \eqref{weak_CV_wJ}), we can choose $J$ so large that 
\begin{equation}
\label{weakCVw0n}
\rhob_n^{1/2} w_{0,n}^J(\rhob_n \cdot)\xrightharpoonup[n\to\infty]{} 0 \text{ in }\hdot.
\end{equation} 
The map $u\mapsto u(1)$ is a continuous linear form on the vector space of radial functions in $\hdot$. Thus \eqref{weakCVw0n} implies
\begin{equation}
\label{trace0}
\lim_{n\to +\infty} \rhob_n^{1/2}w_{0,n}^J(\rhob_n)=0.
\end{equation} 
To show \eqref{no_wn_outside}, we argue by contradiction. Assume after extraction (in $n$) that for large $n$
\begin{equation*}
\int_{|x|\geq \overline{\rho}_n} \left(|\nabla w_{0,n}^J|^2+(w_{1,n}^J)^2\right)dx\geq \eps_0.
\end{equation*}
Then by \eqref{IPPwn} and \eqref{trace0}, for large $n$,
\begin{equation}
\label{minor_wnJ}
\int_{\overline{\rho}_n}^{+\infty} (\partial_r (rw_{0,n}^J))^2dr+(rw_{1,n}^J)^2dr\geq \frac{\eps_0}{2}.
\end{equation}
By Lemma \ref{L:linear_behavior}, and still extracting subsequences, the following holds for all $t>0$ or all $t<0$, and for all large $n$,
\begin{equation*}
\int_{\overline{\rho}_n+|t|}^{+\infty} (\partial_r (rw_{n}^J(t)))^2+(r\partial_tw_{n}^J(t))^2\,dr\geq \frac{\eps_0}{4}.
\end{equation*}
By  \eqref{int_f_1}, this implies that for all $t>0$ or for all $t<0$,
\begin{equation*}
\int_{|x|\geq \overline{\rho}_n+|t|}\left|\nabla w_{n}^J(t)\right|^2+(\partial_t w_{n}^J(t))^2\,dx\geq \frac{\eps_0}{4}.
\end{equation*}
By finite speed of propagation, we have
$$ \tw_n^J(t,x)=w_n^J(t,x),\quad |x|\geq 2\trho_n+|t|.$$
As $\trho_n\ll \rhob_n$, we obtain for large $n$,
\begin{equation}
\label{channelwn}
\int_{|x|\geq \overline{\rho}_n+|t|}\left|\nabla \tw_{n}^J(t)\right|^2+(\partial_t \tw_{n}^J(t))^2\,dx\geq \frac{\eps_0}{4},
\end{equation}
for all $t\geq 0$ or for all $t\leq 0$. In view of Claim \ref{C:exterior}, the equality \eqref{approx_tu'} can be rewritten
\begin{equation}
\label{approx_tu''}
\tu_n(\tau_n+t)=v(\tau_n+t,x)+\tw_{n}^J(t,x)+r_n^J(t,x),\quad 0\leq \tau_n+t<1.
\end{equation} 
Taking $t=-\tau_n$ if \eqref{channelwn} holds for all $t<0$, and $t=\frac{1-\tau_n}{2}$ if \eqref{channelwn} holds for all $t>0$, we get a contradiction as in the proof of Claim \ref{C:exterior}, concluding the proof of Proposition \ref{P:ext_profile}.
\end{proof}

\section{Rigidity argument for non-self-similar blow-up}
\label{S:non_self_similar}
In the section we consider, as in the preceding one, a radial solution in space dimension $N=3$ that blows up at time $T=1$ and satisfies \eqref{bound_nabla'}. We assume in addition that there exist sequences $\{\tau_n\}$, $\{\lambda_n\}$ such that $\tau_n\in(0,1)$, $\tau_n\to 1$ and
\begin{gather}
 \label{lambdan_small}
\lambda_n \ll 1-\tau_n\\
\label{concentration_a}
\lim_{n\to +\infty}\int_{|x|\geq \lambda_n} |\nabla a(\tau_n)|^2+(\partial_t a(\tau_n))^2+\frac{1}{|x|^2}(a(\tau_n))^2=0, 
\end{gather} 
where $a=u-v$ is as usual the singular part of $u$. The main result of this section is the following
\begin{prop}
\label{P:sum_of_W}
Assume that $u$ is radial and that \eqref{bound_nabla'}, \eqref{lambdan_small} and \eqref{concentration_a} hold. Then there exist a sequence $(t_n)$, $J_0>0$, $(\iota_j)_{j=1\ldots J_0}\in \{\pm 1\}^{J_0}$, and, for $j=1\ldots J_0$, sequences $\{\lambda_{j,n}\}_n$ of positive numbers, such that
\begin{align}
\label{sum_of_W1}
u(t_n,x)&=v(t_n,x)+\sum_{j=1}^{J_0} \frac{\iota_j}{\lambda_{j,n}^{1/2}}W\lf(\frac{x}{\lambda_{j,n}}\rg)+w_{0,n}&\text{ in }\hdot\\
\label{sum_of_W2}
\partial_t u(t_n,x)&=\partial_t v(t_n,x)+o_n(1)&\text{ in }L^2,
\end{align}
where, denoting by $w_n$ the solution of \eqref{lin_wave} with initial data $(w_{0,n},0)$,
$$\lim_{n\to \infty}\|w_n\|_{S(\RR)}=0,$$
\end{prop}
Let us mention that the assumption $N=3$ is not essential for the arguments of this section.
In \S \ref{SS:virial} we show that assumptions \eqref{lambdan_small}, \eqref{concentration_a} imply that $\partial_ta$ is small in $L^2$ for a sequence of times. Proposition \ref{P:sum_of_W} is proven in \S \ref{SS:proof_sum_of_W}.
\subsection{Smallness of the time-derivative of the solution}
\label{SS:virial}
\begin{lemma}
\label{L:dt=0}
Assume \eqref{concentration_a}. Then
$$ \lim_{n\to \infty}\frac{1}{1-\tau_n}\int_{\tau_n}^1 \int_{\RR^3}(\partial_t a)^2dx\, dt=0.$$
\end{lemma}
\begin{corol}
\label{C:dt=0}
Under the assumption of Lemma \ref{L:dt=0}, there exists an increasing sequence $t_n\to 1$, $t_n\in(0,1)$ such that
\begin{gather}
\label{CVdt=0}
\forall n,\quad \int_{\RR^3} |\partial_t a(t_n,x)|^2dx\leq \frac{1}{n},\\
\label{CVintdt=0}
\forall n,\;\forall \sigma\in(0,1-t_n),\quad \frac{1}{\sigma}\int_{t_n}^{t_n+\sigma} \int_{\RR^3} |\partial_t a(t,x)|^2dx\,dt\leq \frac{1}{n}.
\end{gather}
\end{corol}
Let us first assume Lemma \ref{L:dt=0} and prove Corollary \ref{C:dt=0}.
\begin{proof}
Using that the map $t\mapsto \partial_t a(t,\cdot)$ is continuous from $(0,1)$ to $L^2(\RR^3)$ we get \eqref{CVdt=0} from \eqref{CVintdt=0} letting $\sigma\to 0$.

To show \eqref{CVintdt=0}, we argue by contradiction. The existence of a sequence $\{t_n\}$ satisfying \eqref{CVintdt=0} is equivalent to
\begin{equation*}
\forall \eps>0,\; \forall t_*\in (0,1),\;\exists t_0\in(t_*,1),\; \forall \sigma \in(0,1-t_0),\quad \frac{1}{\sigma}\int_{t_0}^{t_0+\sigma}\int_{\RR^3}(\partial_t a)^2dx\,dt\leq \eps.
\end{equation*}
Assume
\begin{equation}
\label{absurd_T}
\exists \eps>0,\; \exists t_*\in (0,1),\;\forall t_0\in(t_*,1),\; \exists \sigma  \in(0,1-t_0),\quad \frac{1}{\sigma}\int_{t_0}^{t_0+\sigma}\int_{\RR^3}(\partial_t a)^2dx\,dt> \eps.
\end{equation}
By Lemma \ref{L:dt=0} we can fix a large  $n$ such that $\tau_n>t_*$ and
$$ \frac{1}{1-\tau_n}\int_{\tau_n}^1\int_{\RR^3}(\partial_t a)^2dx\,dt\leq \frac{\eps}{2}.$$
Let 
$$ A=\left\{\sigma\in(0,1-\tau_n)\;\big|\; \frac{1}{\sigma}\int_{\tau_n}^{\tau_n+\sigma}\int_{\RR^3}(\partial_t a)^2dx\,dt\geq \eps\right\}.$$
By \eqref{absurd_T}, $A$ is not empty. Furthermore, it is closed in $(0,1-\tau_n)$. Let $\theta_0=\sup A$. By the choice of $n$, $\theta_0\neq 1-\tau_n$. Furthermore,
$$ \int_{\tau_n}^{\tau_n+\theta_0} \int_{\RR^3}(\partial_t a)^2dx\,dt\geq \eps\theta_0.$$
By \eqref{absurd_T}, using that $t_*<\tau_n+\theta_0<1$, there exists $\sigma\in (0,1-\tau_n-\theta_0)$ such that
$$ \int_{\tau_n+\theta_0}^{\tau_n+\theta_0+\sigma} \int_{\RR^3}(\partial_t a)^2dx\,dt> \eps\sigma.$$
Summing up the two preceding inequalities, we get
$$ \int_{\tau_n}^{\tau_n+\theta_0+\sigma} \int_{\RR^3}(\partial_t a)^2dx\,dt>\eps (\theta_0+\sigma),$$
with $\theta_0+\sigma\in (\theta_0,1-\tau_n)$. Thus $\theta_0+\sigma\in A$, and $\theta_0+\sigma>\theta_0$ contradicting the fact that $\theta_0=\sup A$.
\end{proof}
It remains to prove Lemma \ref{L:dt=0}.
\begin{proof}[Proof of Lemma \ref{L:dt=0}]
Let 
$$ z_1(t)=\int_{\RR^3} \big[(u\partial_t u)-(v\partial_t v)\big]dx,\quad z_2(t)=\int_{\RR^3}\big[x\cdot \nabla u\,\partial_tu-x\cdot \nabla v\, \partial_tv\big]dx.$$ 
As $u-v$ and $\partial_t(u-v)$ are compactly supported in the space variable, both integrals are well-defined. We first show
\begin{equation}
\label{z_tau_n0}
\lim_{n\to +\infty} \frac{|z_1(\tau_n)|+|z_2(\tau_n)|}{1-\tau_n}=0.
\end{equation} 
Indeed, write
$$ z_1(\tau_n)=\int_{\RR^3} a(\tau_n)\partial_t u(\tau_n)+\int_{\RR^3} v(\tau_n)\partial_t a(\tau_n).$$
Then, using that on the supports of $a(\tau_n)$ and $\partial_t a(\tau_n)$, $|x|\leq 1-\tau_n$,
\begin{multline*}
\left|\int a(\tau_n)\partial_t u(\tau_n)\right|\leq \int_{|x|\leq \lambda_n}\left|a(\tau_n)\partial_t u(\tau_n)\right|+\int_{|x|\geq \lambda_n}\left|a(\tau_n)\partial_t u(\tau_n)\right|\\
\leq \lambda_n\int_{|x|\leq \lambda_n}\frac{1}{|x|}\left|a(\tau_n)\partial_t u(\tau_n)\right|+(1-\tau_n)\int_{|x|\geq \lambda_n}\frac{1}{|x|}\left|a(\tau_n)\partial_t u(\tau_n)\right|.
\end{multline*}
By \eqref{lambdan_small} and \eqref{concentration_a},
\begin{equation*}
\left|\int a(\tau_n)\partial_t u(\tau_n)\right|=o(1-\tau_n),\text{ as }n\to \infty.
\end{equation*}
Estimating the other terms in the same way we get \eqref{z_tau_n0}.

Differentiating the definitions of $z_1$ and $z_2$ and using that both $u$ and $v$ are solutions of \eqref{CP}, we get
\begin{align*}
z'_1(t)&=\int (\partial_t u)^2-\int|\nabla u|^2+\int u^{6}-\left[\int (\partial_t v)^2dx-\int|\nabla v|^2dx+\int v^{6}\right]\\
z_2'(t)&=-\frac{3}{2}\int (\partial_t u)^2+\frac{1}{2}\left(\int |\nabla u|^2-\int u^{6}\right)\\
&\qquad-\left[-\frac{3}{2}\int (\partial_t v)^2+\frac{1}{2}\left(\int |\nabla v|^2-\int v^{6}\right)\right].
\end{align*}
Noting that $|x|\leq 1-t$ on the support of $a$, that $v$ converges in $\hdot\times L^2$ as $t\to 1$ and that $u$ is bounded in $\hdot \times L^2$, we get, as $t\to 1^-$,
\begin{align*}
z_1'(t)&=\int (\partial_t a)^2dx-\int|\nabla a|^2dx+\int a^{6}dx+o(1)\\
z_2'(t)&=-\frac{3}{2}\int (\partial_t a)^2+\frac{1}{2}\left(\int |\nabla a|^2-\int a^{6}\right)+o(1).
\end{align*}
Let 
$$Z(t)=\frac{1}{2}z_1(t)+z_2(t).$$
Then
\begin{equation*}
Z'(t)=-\int (\partial_t a)^2+o(1)\text{ as } t\to 1.
\end{equation*}
Let $\eps>0$, and $m,n$ be two large integers with $n<m$. Integrating the preceding inequality we get
$$ \int_{\tau_n}^{\tau_m}\int_{\RR^3}(\partial_t a)^2\leq \left|Z(\tau_m)-Z(\tau_n)\right|+\eps|\tau_n-\tau_m|$$
Letting $m$ tends to infinity we obtain
$$ \int_{\tau_n}^1\int_{\RR^3}(\partial_t a)^2\leq |Z(\tau_n)|+o(1-\tau_n)\text{ as }n\to\infty.$$
From \eqref{z_tau_n0} we deduce
$$ \int_{\tau_n}^1\int_{\RR^3}(\partial_t a)^2=o(1-\tau_n)\text{ as }n\to\infty.$$
\end{proof}
\subsection{Decomposition into a sum of rescaled stationary solutions}
\label{SS:proof_sum_of_W}
The proof of Proposition \ref{P:sum_of_W}, is divided into four steps.

\EMPH{Step 1. Extraction of a sequence and profile decomposition}
Extracting a subsequence from $\{t_n\}$, we assume that $\{(a(t_n),\partial_ta(t_n))\}_n$ admits a profile decomposition with profiles $U^j$ and parameters $\lambda_{j,n}$, $t_{j,n}$. By the Pythagorean expansion
$$ \|\partial_t a(t_n)\|^2_{L^2}=\sum_{j=1}^J \left\|\partial_t U^j(-t_{j,n}/\lambda_{j,n})\right\|_{L^2}^2+\lf\|w_n^J\rg\|_{L^2}^2+o(1)\text{ as }n\to\infty$$
and using \eqref{CVdt=0}, we get that for all $j$ (here $U^j_n$ is the rescaled profiled, defined in Notation \ref{N:rescaled}),
\begin{equation}
\label{small_partialt}
\lim_{n\to +\infty} \|\partial_t U_n^j(0)\|_{L^2}=\lim_{n\to +\infty} \lf\|\partial_t U^j(-t_{j,n}/\lambda_{j,n})\rg\|_{L^2}=0. 
\end{equation} 
We deduce that for all $j$ such that $U^j\neq 0$, $\lf\{-t_{j,n}/\lambda_{j,n}\rg\}_n$ is bounded. Indeed assume that there exists a subsequence in $n$ such that $-t_{j,n}/\lambda_{j,n}\to\pm \infty$. Then by definition of $U^j$ and the equipartition of the energy for solutions of the linear equation \eqref{lin_wave} as $t\to \pm\infty$,
$$ \frac 12 \lf\|U^j_{\lin}(0)\rg\|_{\hdot}^2+\frac 12 \lf\|\partial_t U^j_{\lin}(0)\rg\|_{L^2}^2=\lim_{n\to \infty}\lf\|\partial_t U^j_{\lin}(-t_{j,n}/\lambda_{j,n})\rg\|_{L^2}^2=0,$$
showing that $U^j=0$, a contradiction.

Translating in time the profiles, we may assume
\begin{equation}
\label{tjn0}
\forall j,\; \forall n,\quad t_{j,n}=0.
\end{equation} 
As a consequence of \eqref{small_partialt}, $U_1^j:=\partial_tU^j(0)=0$ for all $j$.
Let $\delta_0>0$ be a small parameter (given by the small data theory for \eqref{CP}). There exists a finite number $J_0$ of profiles $U^j$ such that $ \big\| U_0^j\big\|_{\hdot}+\big\| U_1^j\big\|_{L^2}=\big\| U_0^j\big\|_{\hdot}\geq   \delta_0$. Reordering the profiles, we may assume
$$\big\| U_0^j\big\|_{\hdot}\geq \delta_0
\iff 1\leq j\leq J_0.$$
In view of \eqref{tjn0} and the orthogonality of the profiles, we obtain, after a new extraction in $n$, 
$$ \forall j,k \quad j\neq k\Longrightarrow \lambda_{j,n}\ll \lambda_{k,n}\text{ or }\lambda_{k,n}\ll \lambda_{j,n}.$$
Thus we may reorder the first profiles so that
$$ \lambda_{J_0,n}\ll \lambda_{J_0-1,n}\ll \ldots \ll \lambda_{1,n}.$$
We show by contradiction that $U^j\in \{W,-W\}$ if $1\leq j\leq J_0$ and $U^j=0$ if $j>J_0$. This is equivalent to the fact that the set of indexes $j$ such that $U^j\notin \{0,W,-W\}$ is empty. Assume that this set is not empty and let
$$ k_0=\min\big\{j\geq 1,\; U^j\notin\{0,W,-W\}\big\}.$$ 
Let $k_1=\min\big\{1\leq j\leq J_0\;\big|\;\lambda_{j,n}\ll \lambda_{k_0,n}\big\}$. If this set is empty, $k_1$ is not defined, and we will make the convention $\lambda_{k_1,n}=0$.
By Claim \ref{C:faraway}, there exists a sequence $\tlambda_n\to 0$ such that
\begin{gather}
\label{tlambda_lambdak0}
\lambda_{k_1,n}\ll \tlambda_n\ll \lambda_{k_0,n}\\
\label{tlambda_n_ortho}
\forall j,\quad \tlambda_n \ll \lambda_{j,n} \text{ or }\lambda_{j,n}\ll \tlambda_n.
\end{gather}
Let 
$$\JJJ_{\ext}= \Big\{j\geq 1,\; \tlambda_n\ll \lambda_{j,n}\Big\}.$$
Note that by the first inequality in \eqref{tlambda_lambdak0},
\begin{equation}
\label{tlambda_ext}
\forall j,\quad \Big(j\in \JJJ_{\ext} \text{ and }\lambda_{j,n}\ll \lambda_{k_0,n}\Big)\Longrightarrow j>J_0.
\end{equation} 

\EMPH{Step 2} 
Let $T>0$ be in the domain of existence of $U^{k_0}$. Using that $\lambda_{k_0,n}\lesssim 1-\tau_n$, we can choose $T$ small enough so that for large $n$, $\lambda_{k_0,n}T<1-t_n$. 
In this step we show
\begin{equation}
\label{profile_step}
\frac{1}{\lambda_{k_0,n}T}\int_{0}^{\lambda_{k_0,n}T} \int_{|x|\geq \lambda_{k_0,n}\eps+|t|}\bigg|  \partial_tU^{k_0}_n\lf(t,x\rg)+\partial_tw_n^J(t,x)\bigg|^2 dx\,dt=o_n^J,
\end{equation}
where by definition $\lim_{J\to \infty} \limsup_{n\to \infty} o_{n}^J=0$. More precisely, we will show the following two estimates which directly imply \eqref{profile_step}:
\begin{gather}
\label{small_ext}
\frac{1}{\lambda_{k_0,n}T}\int_{0}^{\lambda_{k_0,n}T} \int_{|x|\geq \lambda_{k_0,n}\eps+|t|}\bigg| \sum_{\substack{j\in \JJJ_{\ext}\\ j\leq J}} \partial_tU^j_n\lf(t,x\rg)+\partial_tw_n^J(t,x)\bigg|^2 dx\,dt=o_{n}^J\\
\label{contribution_0}
j\in \JJJ_{\ext} \text{ and }j\neq k_0\Longrightarrow
\lim_{n\to \infty}\frac{1}{\lambda_{k_0,n}T}\int_{0}^{\lambda_{k_0,n}T} \int_{|x|\geq \eps\lambda_{k_0,n}+|t|} \lf|\partial_tU^j_n\lf(t,x\rg)\rg|^2dx\,dt=0.
\end{gather} 

\EMPH{Proof of \eqref{small_ext}}

Consider a radial function $\chi\in C^{\infty}(\RR^3)$, such that $\chi(x)=1$ for $|x|\geq 2$ and $\chi(x)=0$ for $|x|\leq 1$. Let $\tilde{u}_n$ be the solution of \eqref{CP} with initial data
$$ \tilde{u}_{n\restriction t=t_n}=\chi\lf(\frac{x}{\tlambda_n}\rg)u(t_n,x),\quad \partial_t\tu_{n\restriction t=t_n}=\chi\lf(\frac{x}{\tlambda_{n}}\rg)\partial_t u(t_n,x).$$
Then, by finite speed of propagation, as long as $t_n+s$ is in the domain of existence of $u$ and $\tu_n$, 
$$ \tu_n(t_n+s,x)=u(t_n+s,x),\quad |x|\geq |s|+2\tlambda_n.$$
Furthermore, letting 
$$ \tilde{w}_{0,n}^J(x)=\chi\lf(\frac{x}{\tlambda_n}\rg) w_{0,n}^J(x),\quad \tilde{w}_{1,n}^J(x)=\chi\lf(\frac{x}{\tlambda_n}\rg) w_{1,n}^J(x),$$
we obtain (recall that $U_{1,n}^j=0$ for all $j$)
\begin{align}
\label{dev_chi0a}
\tu_n(t_n,x)-v(t_n,x)&=\sum_{\substack{j \in \JJJ_{\ext}\\j\leq J}} U_{0,n}^j(x)+\tilde{w}_{0,n}^J(x)+o_n(1)\text{ in }\hdot\\
\label{dev_chi1a}
\partial_t \tu_n(t_n,x)-\partial_t v(t_n,x)&=\tw_{1,n}^J(x)+o_n(1)\text{ in }L^2,
\end{align}
By Claim \ref{C:dispersive_tw} together with the argument than we used to show \eqref{small_tw},
\begin{equation}
\label{dispersive_tw'}
\lim_{n\to\infty} \limsup_{J\to+\infty}\left\|\tw_{n}^{J}\right\|_{S(\RR)}=0.
 \end{equation} 
By \eqref{dispersive_tw'}, the two equations \eqref{dev_chi0a}, \eqref{dev_chi1a} yield a profile decomposition of the sequence $$\lf\{\tu(t_n,x)-v(t_n,x),\partial_t \tu(t_n,x)-\partial_t v(t_n,x)\rg\}_n.$$

The development \eqref{dev_chi0a}, \eqref{dev_chi1a} satisfies the assumptions of Proposition \ref{P:lin_NL} with $\theta_n=\lambda_{k_0,n}T$. Indeed for $j>J_0$ the solution $U^j$ scatters both forward and backward in time. Furthermore by \eqref{tlambda_ext}, 
$$\left( j\in \{1,\ldots,J_0\} \cap \JJJ_{\ext} \text{ and } j\neq k_0\right) \Longrightarrow \lambda_{k_0,n}T\ll \lambda_{j,n}.$$
Thus by Proposition \ref{P:lin_NL}, for $s\in  [0,\lambda_{k_0,n}T]$
\begin{equation}
\label{dev_tun}
\tu_n(t_n+s,x)=v(t_n+s,x)+\sum_{\substack{j\in \JJJ_{\ext}\\j\leq J}} U^j_n\lf(s,x\rg)+\tilde{w}_n^J(s,x)+r_n^J(s,x),
\end{equation}
where $r_n^J$ satisfies \eqref{cond_rJn} with $\theta_n=\lambda_{k_0,n}T$.
Let $\eps>0$. We have, for large $n$ (so that $\eps\lambda_{k_0,n}\geq 2\tlambda_n$),
\begin{multline*}
o_n(1)=\frac{1}{\lambda_{k_0,n}T}\int_{t_n}^{t_n+\lambda_{k_0,n}T} \int_{\RR^3}|\partial_t a(t,x)|^2 dx\,dt\\ \geq\frac{1}{\lambda_{k_0,n}T}\int_{t_n}^{t_n+\lambda_{k_0,n}T} \int_{|x|\geq \lambda_{k_0,n}\eps+|t|}|\partial_t a(t,x)|^2 dx\,dt \\
=\frac{1}{\lambda_{k_0,n}T}\int_{t_n}^{t_n+\lambda_{k_0,n}T} \int_{|x|\geq \lambda_{k_0,n}\eps+|t|}|(\partial_t u_n-\partial_t v)(t,x)|^2 dx\,dt,
\end{multline*}
which yields \eqref{small_ext} in view of \eqref{dev_tun}.

\EMPH{Proof of \eqref{contribution_0}}

 Let $R_{j,n}=\lambda_{k_0,n}/\lambda_{j,n}$. We have
\begin{multline}
\label{coucou}
\frac{1}{\lambda_{k_0,n}T}\int_{0}^{\lambda_{k_0,n}T} \int_{|x|\geq \eps\lambda_{k_0,n}+|t|} \lf|\partial_tU^j_n\lf(t,x\rg)\rg|^2dx dt\\
=\frac{1}{TR_{j,n}}\int_{0}^{TR_{j,n}}\int_{|y|\geq \eps R_{j,n}+|s|}\lf|\partial_tU^j\lf(s,y\rg)\rg|^2dy ds.
\end{multline}

If $\lambda_{j,n}\ll \lambda_{k_0,n}$ (and thus $j>J_0$), we have that $R_{j,n}\to +\infty$. 
By finite speed of propagation, for all $\eta>0$, there exists $M>0$ such that 
$$\forall s\in \RR,\quad \int_{|y|\geq M+|s|}\lf|\partial_t U^j\lf(s,y\rg)\rg|^2dy\leq \eta,$$ 
which implies that the right-hand member of \eqref{coucou} tends to $0$ as $n\to \infty$. 

If $\lambda_{k_0,n}\ll \lambda_{j,n}$, $R_{j,n}\to 0$, and thus
\begin{multline*}
\frac{1}{TR_{j,n}}\int_{0}^{TR_{j,n}}\int_{|y|\geq R_{j,n}\eps+|s|}\lf|\partial_tU^j\lf(s,y\rg)\rg|^2dy ds
\leq \frac{1}{TR_{j,n}}\int_{0}^{TR_{j,n}}\int\lf|\partial_tU^j\lf(s,y\rg)\rg|^2dy ds\\
\xrightarrow[n\to\infty]{}\int\lf|\partial_tU^j\lf(0,y\rg)\rg|^2dy=0,
\end{multline*}
concluding the proof of \eqref{contribution_0}.

\EMPH{Step 3. Uniqueness argument and conclusion of the proof}

By \eqref{profile_step},
\begin{equation*}
\frac{1}{T}\int_{0}^{T} \int_{|x|\geq \eps+|t|}\bigg|  \partial_tU^{k_0}\lf(t,x\rg)+\lambda_{k_0,n}^{3/2}\partial_tw_n^J(\lambda_{k_0,n}t,\lambda_{k_0,n}x)\bigg|^2 dx\,dt=o_n^J.
\end{equation*}
Consider the mapping $\hdot\times L^2\longrightarrow  \RR$
\begin{equation*}
(f_0,f_1)\longmapsto \frac{1}{T}\int_{0}^{T} \int_{|x|\geq \eps+|t|}  \partial_tU^{k_0}\lf(t,x\rg)\partial_tf(t,x)dx\,dt,
\end{equation*}
where $f(t,x)$ is the solution of the linear wave equation with initial conditions $(f_0,f_1)$. This is a continuous linear form on $\hdot\times L^2$. By \eqref{weak_CV_wJ},
$$ \lf(\lambda_{k_0,n}^{1/2}w_{0,n}^J\lf(\lambda_{k_0,n}\cdot\rg),\lambda_{k_0,n}^{3/2}w_{1,n}^J\lf(\lambda_{k_0,n}\cdot\rg)\rg)\xrightharpoonup[n\to \infty]{} 0\text{ weakly in }\hdot \times L^2.$$
Hence
\begin{equation*}
\lim_{n\to \infty}
\frac{1}{T}\int_{0}^{T} \int_{|x|\geq \eps+|t|}  \partial_tU^{k_0}\lf(t,x\rg)\lambda_{k_0,n}^{3/2}\partial_tw_n^J(\lambda_{k_0,n}t,\lambda_{k_0,n}x)dx\,dt=0,
\end{equation*}
and we conclude that for all $\eps>0$,
\begin{equation*}
\frac{1}{T}\int_{0}^{T} \int_{|x|\geq \eps+|t|}\bigg|\partial_tU^{k_0}\lf(t,x\rg)\bigg|^2dx\,dt=0
\end{equation*} 
Letting $\eps\to 0$ we get
\begin{equation*}
\frac{1}{T}\int_{0}^{T} \int_{|x|\geq |t|}\bigg|\partial_tU^{k_0}\lf(t,x\rg)\bigg|^2dx\,dt=0
\end{equation*}
This shows that $\partial_t U^{k_0}(t,x)=0$ if $t\leq |x|$ and $0\leq t\leq T$. Let 
$$\Omega=\Big\{(t,x)\in [0,T]\times \RR^3,\; |x|\geq t\Big\}$$
Then 
$$ (t,x)\in\Omega\Longrightarrow U^{k_0}(t,x)=U^{k_0}_0(x).$$
In $\Omega$, the non-linear wave equation $\partial_t^2 U^{k_0}-\Delta U^{k_0}-\big(U^{k_0}\big)^{5}=0$ becomes $\Delta U^{k_0}=-\big(U^{k_0}\big)^{5}$. Thus $U^{k_0}$ satisfies in the sense of distributions the elliptic equation  
$$ \Delta U^{k_0}_{0}=-\big(U^{k_0}_0\big)^{5}\text{ in }\RR^3\setminus \{0\}.$$
This shows that $U^{k_0}$ is smooth in $\RR^3\setminus \{0\}$ and satisfies the preceding equation in the classical sense in $\RR^3\setminus\{0\}$. 
As a consequence $\Delta U^{k_0}_{0}+\big(U^{k_0}_0\big)^{5}$ is a distribution in $H^{-1}(\RR^3)$, supported at the origin. The only distribution with these properties in dimension $3$ is $0$ and we deduce
$$ \Delta U^{k_0}_{0}+\big(U^{k_0}_0\big)^{5}=0$$
in the sense of distributions on $\RR^3$ and thus by Claim \ref{C:uniqueness_W}, as $U^{k_0}$ is radial,
$$ U^{k_0}(x)= \frac{1}{\lambda_0^{1/2}}W\left(\frac{x}{\lambda_0}\right)\text{ or }U^{k_0}(x)=-\frac{1}{\lambda_0^{1/2}}W\left(\frac{x}{\lambda_0}\right)\text{ or }U^{k_0}=0,$$
for some $\lambda_0>0$, which yields the desired contradiction. The proof is complete.
\qed

\section{All radial compact solutions are stationary}
\label{S:radial_compact}
In this section we show Theorem \ref{T:compact=W}.

We will assume without loss of generality that $\lambda$ is continuous on $(T_-(u),T_+(u))$ (see \cite[Remark 5.4]{KeMe06}).

\EMPH{Step 1} 

We show in this step that the solution is globally defined. Assume that $T_+(u)<\infty$. For the sake of simplicity, we will assume that $T_+(u)=1$. By standard argument (see Section \ref{S:reg_sing}), $\lambda(t)\leq C(1-t)$.
By \cite[Section 6]{KeMe08}, self-similar, compact blow-up is excluded, which implies that there exists a sequence $\{\tau_n\}_n$ such that
$$\tau_n\in(0,1),\quad \lim_{n\to\infty} \tau_n=1,\quad \lim_{n\to +\infty}\frac{\lambda(\tau_n)}{1-\tau_n}=0.$$
Using that the regular part of $v$ at the blow-up point $t=1$ is $0$, we get, arguing as in Corollary \ref{C:dt=0} that there exists a sequence $\{t_n\}_n$ such that
\begin{equation}
\label{small_time_deriv}
\forall n,\; \forall \sigma\in (0,1-t_n),\quad \frac{1}{\sigma}\int_{t_n}^{t_n+\sigma}\int_{\RR^N}|\partial_tu(t,x)|^2\,dx\,dt\leq \frac{1}{n}.
\end{equation} 
Consider $(U_0,U_1)\in \hdot \times L^2$ such that for a subsequence,
$$ \lim_{n\to \infty} \left(\lambda^{\frac{N-2}{2}} (t_n)u\left(t_n,\lambda(t_n)x\right),\lambda^{\frac{N}{2}} (t_n)u\left(t_n,\lambda(t_n)x\right)\right)=(U_0,U_1).$$
Let $U$ be the solution of \eqref{CP} with initial condition $(U_0,U_1)$ and $\tau_0\in (0,T_+(U))$. Then by Theorem \ref{T:LTPT}
\begin{equation*}
 \lim_{n\to \infty}\int_0^{\tau_0}\int_{\RR^N} \lambda^N(t_n)\left(\partial_t u\left(t_n+\lambda(t_n)s,\lambda(t_n)x\right)\right)^2\,dx\,ds
=\int_0^{\tau_0}\int_{\RR^N} (\partial_t U(t))^2dt.
\end{equation*}
By \eqref{small_time_deriv}, we obtain
\begin{multline*}
\int_0^{\tau_0}\int_{\RR^N} \lambda^N(t_n)\left(\partial_t u\left(t_n+\lambda(t_n)s,\lambda(t_n)x\right)\right)^2\,dx\,ds\\
=\frac{1}{\lambda(t_n)} \int_0^{\tau_0\lambda(t_n)} \int_{\RR^N} (\partial_t u(t_n+t,x))^2\,dx\, dt\xrightarrow[n\to\infty]{}0.
\end{multline*}
As a consequence, $\partial_t U=0$ for $t\in [0,\tau_0]$. By Claim \ref{C:uniqueness_W}, $U=0$ or $U=W$ up to the invariances of the equation. If $U=0$, then $E(u_0,u_1)=0$, and as the $\|u(t_n)\|_{\hdot}$ tends to $0$, this implies by Claim \ref{C:variationnal} that $u=0$, contradicting our assumption. Thus $U=W$ up to the invariances, and by conservation of the energy we get that $E(u_0,u_1)=E(W,0)$.

The solution $u$ of \eqref{CP} has threshold energy $E(W,0)$, is not globally defined and satisfies $u_0\in L^2$. By the Theorem 2 of \cite{DuMe08}, $N=5$ and $u$ has to be the special solution $W^+$ constructed in this paper, which satisfies $\|u(t)-W\|_{\hdot}\leq e^{ct}$ as $t\to -\infty$. This contradicts the fact that $u$ has compact support in space, concluding step $1$.

\EMPH{Step 2} We assume in this step that $\lambda$ is bounded on $[0,+\infty)$ or on $(-\infty,0]$, and show that $E(u_0,u_1)=E(W,0)$. By time symmetry we can assume that $\lambda$ is bounded on $[0,+\infty)$. By the preceding step,
$$
T_+(u)=+\infty.
$$

Let us fix $\phi\in C^{\infty}_0(\RR^3)$ such that $\phi \equiv 1$ for $|x|\leq 1$, $\phi \equiv 0$ for $|x|\geq 2$. For $R>1$, consider $\phi_R=\phi(x/R)$, $\psi_R=x\phi(x/R)$ and
\begin{equation}
\label{def_rho}
\rho(R)=\sup_{t\in (T_-(u),T_+(u))} \int_{|x|\geq R}\frac{|u|^2}{|x|^2}+|\nabla u|^2+|\partial_t u|^2+|u|^{6}dx.
\end{equation} 
The compactness of $\overline{K}$ and the boundedness of $\lambda$ implies that $\rho(R)$ is finite, and tends to $0$ as $R$ goes to infinity. Let 
$$ y_R(t)=\int_{\RR^3} \psi_R\cdot\nabla u \partial_t udx+\frac{1}{2}\int_{\RR^3} \varphi_R u \partial_t udx.$$
Then (see \cite[Lemma 5.3]{KeMe08})
\begin{equation}
\label{viriel}
y_R'(t)=-\int_{\RR^3} (\partial_t u)^2dx+\OOO(\rho(R)).
\end{equation} 
Integrating with respect to time, we get that there exists a constant $C>0$, independent of $R$, such that for all $T>0$,
$$ \int_{0}^{T} \|\partial_t u(t)\|_{L^2}^2dt\leq |y_R(T)-y_R(0)|+C T\rho(R).$$
using that, for any fixed $R>0$, $y_R(t)$ is bounded independently of $t$, we get
\begin{equation}
\label{lim_ptu1}
\lim_{T\rightarrow+\infty}     \frac{1}{T}\int_{0}^{T} \|\partial_t u(t)\|_{L^2}^2dt=0.
\end{equation}
We next show that there exists a sequence $t_n$ that tends to infinity and such that 
\begin{equation}
\label{lim_partial_tu}
\lim_{n\rightarrow+\infty} \frac{1}{\lambda(t_n)}\int_{t_n}^{t_n+\lambda(t_n)} \|\partial_t u(t)\|^2_{L^2}dt=0.
\end{equation}
Indeed, define a sequence $\tau_n$ by
$$ \tau_0=0,\quad \tau_{n+1}=\tau_n+\lambda(\tau_n).$$
We first show that $\tau_n\rightarrow+\infty$. If not, $\tau_n$ has a finite limit $\tau_{\infty}=\sum_{n\geq 0}\lambda(\tau_n)$, which shows by continuity of $\lambda$ that $\lambda(\tau_{\infty})=0$ a contradiction with the assumption that $\lambda$ takes strictly positive values.

To show \eqref{lim_partial_tu}, we argue by contradiction.
Assuming that no subsequence $\{t_n\}$ of $\{\tau_n\}$ satisfies \eqref{lim_partial_tu}, we get that there exists $\eps>0$ such that
$$ \forall n, \quad \int_{\tau_n}^{\tau_{n+1}} \|\partial_t u(t)\|^2_{L^2}dt\geq \eps \lambda(\tau_n).$$
Summing up, and using that $\tau_{n+1}=\sum_{k=1}^n \lambda(\tau_k)$, we get
$$ \forall n,\quad \frac{1}{\tau_{n+1}}\int_{0}^{\tau_{n+1}} |\partial_t u(t)|^2_{L^2}dt \geq \eps, $$
contradicting \eqref{lim_ptu1}. Hence \eqref{lim_partial_tu}.

Extracting subsequences, we get $(U_0,U_1)\in \hdot\times L^2$ such that
$$ \left(\lambda_n^{N/2-1} u(t_n,\lambda(t_n)x),\lambda_n^{N/2} \partial_t u(t_n,\lambda(t_n)x)\right)\underset{n \to \infty}{\longrightarrow} (U_0,U_1).$$
Let $U$ be the solution of \eqref{CP} with initial conditions $(U_0,U_1)$. Let $\theta_0\in (0,T_+(U))$ such that $\theta_0\leq 1$. Then by Theorem \ref{T:LTPT},
\begin{multline*}
\frac{1}{\theta_0\lambda(t_n)}\int_{t_n}^{t_n+\lambda(t_n)} \big\|\partial_t u(t)\big\|_{L^2}^2dt\geq \frac{1}{\theta_0\lambda(t_n)}\int_{t_n}^{t_n+\theta_0\lambda(t_n)} \big\|\partial_t u(t)\big\|_{L^2}^2dt\\
=\frac{1}{\theta_0\lambda(t_n)}\int_0^{\theta_0\lambda(t_n)} \big\|\partial_t U\left(t/\lambda(t_n)\right)\big\|_{L^2}^2dt+o_n(1)\\
=\int_0^{\theta_0} \big\|\partial_t U\left(s\right)\big\|_{L^2}^2ds+o_n(1).
\end{multline*}
By \eqref{lim_partial_tu}, we get that $\partial_t U=0$ on $[0,\theta_0]$. By Claim \ref{C:uniqueness_W}, $U=W$, which shows that $E(U_0,U_1)=E(u_0,u_1)=E(W,0)$. This concludes Step 2.

\EMPH{Step 3}
We next show that $E(u_0,u_1)=E(W,0)$ also if $\lambda$ is unbounded on both intervals $[0,+\infty)$ and  $(-\infty,0]$. We will use an argument of \cite{KeMe06} to reduce to the previous case. We sketch the argument for the sake of completness. Consider the sequence $\{t_n\}_n$, 
$$ t_n=\inf\big\{t\in [0,+\infty)\;|\;\lambda(t_n)=n\big\}.$$
By continuity of $\lambda$ and the fact that $\lambda(t)$ tends to $+\infty$ as $t$ tends to $+\infty$, we get that $t_n$ is well-defined for large $n$ and
\begin{equation}
\label{prop_tn}
\lim_{n\to\infty} t_n=+\infty,\quad \forall t\in [0,t_n],\; \lambda(t)\leq \lambda(t_n).
\end{equation}
Extracting subsequences, consider $(U_0,U_1)$ such that
$$ \lim_{n\to \infty} \left(\lambda^{N/2-1}(t_n)u(t_n,\lambda(t_n)x),\lambda^{N/2}(t_n)\partial_tu(t_n,\lambda(t_n)x)\right)=(U_0,U_1).$$
Note that we cannot have $(U_0,U_1)=(0,0)$ (this would imply, by Claim \ref{C:variationnal} that $u=0$).
Let $U$ be the solution of \eqref{CP} with initial conditions $(U_0,U_1)$.
By the arguments of \cite[Proof of Theorem 7.1]{KeMe06}, we can show, as a consequence of the compactness of $\overline{K}$ and \eqref{prop_tn}, that there exists a continuous function $\tlambda$ on $(T_-(U),T_+(U))$, bounded on $(T_-(U),0]$ and such that 
$$ \widetilde{K}=\left\{\left(\tlambda^{N/2-1}(t)U\left(t,\tlambda(t)x\right),\tlambda^{N/2}(t)\partial_t U\left(t,\tlambda(t)x\right)\right),\; t\in \left(T_-(U),T_+(U)\right)\right\}$$
has compact closure in $\hdot\times L^2$. By Step 1, $U$ is globally defined. By Step 2, as $\tlambda$ is bounded on $(-\infty,0]$, we get that $E(U_0,U_1)=E(W,0)$. Thus by conservation of the energy of $u$, $E(u_0,u_1)=E(W,0)$ which concludes this step.

\EMPH{Step 4} Convergence in mean to $W$.
By \cite[Theorem 2]{DuMe08}, $\|\nabla u(t)\|_{L^2}^2\geq \|\nabla W\|_{L^2}^2$ for all $t$: if not, $u$ would scatter at least in one time direction, contradicting the compactness of $K$. 

To show the that $u=W$, we will use the arguments of \cite[Section 3]{DuMe08}\footnote{In the cited paper, the notation $\lambda(t)$ stands for the function $1/\lambda(t)$ of the present paper}. In this section, it is shown in particular that a globally defined solution $u$ of \eqref{CP} of energy $E(W,0)$, satisfying $\|\nabla u_0\|_{L^2}^2\leq \|\nabla W\|^2_{L^2}$ and such that there exists $\lambda(t)$ with $K$ compact must be equal to $W$ up to the symmetries of \eqref{CP}. We will quickly check here that the same proof works with a slight modification in the case $\|\nabla u(t)\|_{L^2}^2\geq \|\nabla W\|^2_{L^2}$. As usual, we may assume that $\lambda(t)$ is a continuous function of $t$. Let 
$$\dd(t)=8\int (\partial_tu )^2+4\left(\int |\nabla u|^2-\int |\nabla W|^2\rg)\geq 0.$$
By the characterization of $W$ (\cite{Au76}, \cite{Ta76}), for any $t_0$, $\dd(t_0)=0$ if and only if $u(t_0)\equiv W$ up to the symmetries of the equation. In this case, by uniqueness of the Cauchy problem, $u(t)$ is a stationary solution identically equal to $W$ up to the symmetries. 

In this step we show that 
\begin{equation}
\label{CVmean}
\lim_{T\to +\infty}\frac{1}{T}\int_{-T}^{+T} \dd(t)dt=0
\end{equation} 
Consider a function $\varphi\in C_0^{\infty}$ such that $\varphi=1$ if $|x|\leq 1$, and denote by $\varphi_R(x)=\varphi(x/R)$. 
Let $g_R(t)=2 \int u\partial_t u  \varphi_R$ and note that $|g_R(t)|\leq C_0R$, for a constant $C_0>0$ depending only on $\sup_t \|\partial_tu(t)\|_{L^2}+\|\nabla u(t)\|_{L^2}$. Using that $u$ is solution of \eqref{CP}, we get
\begin{equation}
\label{hR'}
g_R'(t)=\dd(t)+A_R(t),
\end{equation} 
where 
\begin{equation}
\label{bound_AR}
|A_R(t)|\leq \int_{|x|\geq R}\frac{1}{|x|^2}u^2+u^6+|\nabla u|^2+(\partial_t u)^2.
\end{equation} 
As in the case $\|\nabla u(t)\|_{L^2}<\|\nabla W\|_{L^2}$ we will use that $g_R$ and $g_R'$ vanish for $u=W$, and that $|g_R'|$ is larger than $\dd(t)$ up to the remainder term $A_R$. In our case, the definition of $g_R$ is slightly different but it will not affect the proof.

Fix a small $\eps>0$. Using as in the proof of Lemma 3.3 of \cite{DuMe08} that $\lambda(t)/t\to 0$ as $t\to \pm \infty$, we get that there exists a constant $C_1$, independent of $\eps$, and a time $t_1=t_1(\eps)$ such that
$$ \forall T>2t_1(\eps),\; \forall t\in [t_1(\eps),T],\quad g_{\eps T}'(t)\geq \dd(t)-C_1\eps,$$
integrating between $t_1$ and $T$ we get that $\frac{1}{T}\int_0^T \dd(t)dt$ tends to $0$. The same proof works for negative time, yielding \eqref{CVmean}.

\EMPH{Step 5}
In view of \eqref{hR'}, and refining the bound on $g_R(t)$ and the estimate \eqref{bound_AR} on $A_R(t)$ by modulating the solution around $W$ for small $\dd(t)$, we get that there is a constant $C>0$ (depending only on the set $K$) such that
\begin{equation}
\label{virial}
\forall \sigma,\tau\in \RR,\quad \sigma<\tau\Longrightarrow  \int_{\sigma}^{\tau} \dd(t)dt\leq C\left(\sup_{\sigma\leq t\leq \tau} \lambda(t)\right)(\dd(\sigma)+\dd(\tau))
\end{equation} 
(see the proof of Lemma 3.8 in \cite{DuMe08}). Using compactness and modulation arguments, we get the following control on $\lambda(t)$ (see Lemma 3.10 in \cite{DuMe08} and its proof)
\begin{equation}
\label{param_control}
\sigma+\lambda(\sigma)\leq\tau\Longrightarrow \left|\lambda(\sigma)-\lambda(\tau)\right|\leq \int_{\sigma}^{\tau}\dd(t)dt.
\end{equation}  
Consider two sequences $\sigma_n\to -\infty$ and $\tau_n\to +\infty$ such that $\dd(\sigma_n)\to 0$ and $\dd(\tau_n)\to 0$ as $n\to \infty$. The existence of $\{\sigma_n\}_ n$ and $\{\tau_n\}_n$ is given by \eqref{CVmean} in Step 4. Let $n_0$ such that $\dd(\tau_{n_0})\leq \frac{1}{2}$. 
Let us prove by contradiction that $\lambda$ is bounded. 
For large $n$, let $t_n\in[\tau_{n_0},\tau_n]$ such that 
$$\lambda(t_n)=\max_{\tau_{n_0}\leq t\leq \tau_n} \lambda(t).$$ 
If $\lambda(t_n)\to \infty$, then by continuity of $\lambda$, $t_n\to \infty$.  In particular for large $n$, $\tau_{n_0}+\lambda(\tau_{n_0})\leq t_n$, and we can deduce from \eqref{virial} and \eqref{param_control} that
$$   \lambda(t_n)\leq \lambda(\tau_{n_0})+\lambda(t_n)\lf(\frac{1}{2}+\dd(\tau_n)\rg),$$
a contradiction if $\lambda(t_n)\to +\infty$. Thus $\lambda$ is bounded on $[0,+\infty)$ and a similar proof yields the boundedness of $\lambda$ on $(-\infty,0]$. As a consequence of \eqref{virial}, we get
\begin{equation*}
\int_{\sigma_n}^{\tau_n} \dd(t)dt\leq C(\dd(\sigma_n)+\dd(\tau_n)),
\end{equation*} 
which implies that $\dd(t)=0$ for all $t$, concluding the sketch of the proof.

\qed

\section{Bounded globally defined solutions are not self-similar}
\label{S:global}
This section is dedicated to the proof of the following proposition, which will be needed in Section \ref{S:radial} and uses some of the material of Section \ref{S:reg_sing}:
\begin{prop}
\label{P:norm_infinity}
Assume that $N=3$. There exists a constant $\eta_1>0$ with the following property. Let $u$ be a spherically symmetric solution of \eqref{CP} such that $T^+(u)=+\infty$, which does not scatter for positive time and such that
\begin{equation}
\label{bound_nabla_global}
\sup_{t\geq 0} \|\nabla u(t)\|_{L^2}^2+\|\partial_t u(t)\|_{L^2}^2\leq \|\nabla W\|_{L^2}^2+\eta_1.
\end{equation}
Define
\begin{equation}
\label{def_mu2}
\nu(t)=\inf\left\{\mu\;:\; \int_{|x|\geq \mu}|\partial_t u(t)|^2+|\nabla u(t)|^2\leq \frac{1}{2}\int|\nabla W|^2\right\}.
\end{equation} 
Then there exists a sequence $t_n\to \infty$ such that
\begin{equation}
\label{non_self_similar}
\lim_{n\to \infty} \frac{\nu(t_n)}{t_n}=0.
\end{equation} 
\end{prop}
\begin{proof}
We argue by contradiction. Assume that \eqref{non_self_similar} does not hold. Taking into account the finite speed of propagation, we deduce that there exist $c_0,C_0$ such that 
\begin{equation}
\label{mu=t}
\forall t\geq 1,\quad c_0 t\leq \nu(t)\leq C_0t.
\end{equation} 

\EMPH{Step 1}
Let $\AAA$ be the set of $\coltwo{U_0}{U_1}$ such that there exists $t_n\rightarrow +\infty$ with
$$\coltwo{t_n^{1/2}u(t_n,t_nx)}{t_n^{3/2}\partial_tu(t_n,t_nx)} \xrightharpoonup[n\rightarrow +\infty]{} \coltwo{U_0}{U_1} \text{ weakly in }\hdot\times L^2.$$
In this step we show that there is a $(A_0,A_1)\in \AAA$ with minimal energy, that is such that
\begin{equation}
\label{min_E}
\forall (U_0,U_1)\in \AAA, \quad E(A_0,A_1)\leq E(U_0,U_1).
\end{equation} 

We first show that $\AAA$ is sequentially closed in $\hdot\times L^2$ for the weak topology. Indeed, let $(U_{0,n},U_{1,n})\rightharpoonup (U_0,U_1)$, with $(U_{0,n},U_{1,n})\in \AAA$. Consider a countable family of smooth compactly supported functions $(\varphi_j,\psi_j)_{j\in\NN}$ which is dense in $\dot{H}^{-1}\times L^2$. Then for all $k$, there exists $n_k$ such that
$$ \left|\int (U_{0,n_k}-U_0) \varphi_j\right|+\left|\int (U_{1,n_k}-U_1) \psi_j\right|\leq \frac 1k, \quad j=0\ldots k.$$
Thus there exists $t_{k}\geq k$ such that
$$ \left|\int (t_{k}^{1/2}u(t_{k},t_{k}x)-U_0) \varphi_j\right|+\left|\int (t_{k}^{3/2}\partial_t u(t_{k},t_{k}x)-U_1) \psi_j\right|\leq \frac 2k, \quad j=0\ldots k.$$
This shows that $(t_{k}^{1/2}u(t_{k},t_{k}x),t_{k}^{3/2}\partial_t u(t_{k},t_{k}x))$ converges weakly to $(U_0,U_1)$ and thus that $(U_0,U_1)\in \AAA$.

We next construct the minimizing element $(A_0,A_1)$ of $\AAA$. Let $\left\{(U_{0,n},U_{1,n})\right\}_n$ be a sequence in $\AAA$ minimizing the energy. As $\left\{(U_{0,n},U_{1,n})\right\}_n$ is bounded in $\hdot\times L^2$, we can extract a subsequence from $\left\{(U_{0,n},U_{1,n})\right\}_n$ such that
$$ (U_{0,n},U_{1,n})\xrightharpoonup[n \to\infty]{}(A_{0},A_1)\in \AAA.$$
Denote by $\tw_{0,n}=U_{0,n}-A_0$, $\tw_{1,n}=U_{1,n}-A_1$. 
Writing after extraction of a subsequence the profile decomposition of the sequence $(U_{0,n},U_{1,n})$ and using the Pythagorean expansions \eqref{pythagore1a}, \eqref{pythagore1b} and \eqref{pythagore2}, we get 
\begin{gather}
\label{pythagore_nabla}
\|\nabla U_{0,n}\|^2_{L^2}+\|U_{1,n}\|^2_{L^2}=\|\nabla A_{0}\|^2_{L^2}+\|A_1\|^2_{L^2}+\|\nabla \tw_{0,n}\|^2_{L^2}+\|\tw_{1,n}\|_{L^2}^2+o_n(1),\\
\label{pythagore_E}
E(U_{0,n},U_{1,n})=E(A_{0},A_1)+E(\tw_{0,n},\tw_{1,n})+o_n(1).
\end{gather}
By \eqref{pythagore_nabla} and assumption \eqref{bound_nabla_global}, we obtain, for large $n$, 
$$\|\nabla A_{0}\|^2_{L^2}+\|A_1\|^2_{L^2}+\|\nabla \tw_{0,n}\|^2_{L^2}+\|\nabla \tw_{1,n}\|_{L^2}^2\leq \|\nabla W\|_{L^2}^2+2\eta_1,$$
which shows by Claim \ref{C:variationnal} that in \eqref{pythagore_E}, all the energies are positive. 
Thus
$$\inf_{(V_0,V_1)\in \AAA} E(V_0,V_1)=\lim_{n\rightarrow+\infty} E(U_{0,n},U_{1,n})\geq E(A_{0},A_1),$$
implying that $(A_0,A_1)$ satisfies \eqref{min_E}.

\EMPH{Step 2. Profile decomposition}

Consider an arbitrary positive sequence $\{\tau_n\}_n$ that tends to $+\infty$ and such that
\begin{equation}
\label{CV_to_A}
\lf( \tau_n^{1/2}u(\tau_n,\tau_nx),\tau_n^{3/2}\partial_tu(\tau_n,\tau_nx)\rg)\xrightharpoonup[n\to +\infty]{} (A_0,A_1),\text{ weakly in }\hdot\times L^2,
\end{equation} 
where $(A_0,A_1)$ is the minimal element of $\AAA$ defined in Step 1.

Extracting a subsequence from $\{\tau_n\}_n$, we can assume that their exists a profile decomposition $\lf\{U^j_{\lin}\rg\}$, $\lf\{\lambda_{j,n},t_{j,n}\rg\}$ associated to the sequence $(u(\tau_n),\partial_tu(\tau_n))_n$. 

Reordering the profiles, we may assume
\begin{equation}
\label{U1_biggest}
\big\|\nabla U^1_0\big\|_{L^2}^2+\big\|U^1_1\big\|_{L^2}^2=\sup_{j\geq 1} \big\|\nabla U^j_0\big\|_{L^2}^2+\big\|U^j_1\big\|_{L^2}^2.
\end{equation} 
We remark that 
\begin{equation}
\label{bound_U1}
\left\|\nabla U^1_0\right\|_{L^2}^2+\left\|U^1_1\right\|_{L^2}^2\geq \frac{2}{3}\|\nabla W\|_{L^2}^2.
\end{equation}
If not, the result of \cite{KeMe08} would imply that all nonlinear profiles $U^j$ scatter showing by Proposition \ref{P:lin_NL} that $u$ scatters for both positive and negative times, which contradicts our assumption.

As a consequence, we get from \eqref{bound_nabla_global} and again the result of \cite{KeMe08} that for all $j\geq 2$, the nonlinear profile $U^j$ scatters both for positive and negative time.


Extracting and time translating $U^1_{\lin}$ if necessary, we may distinguish three cases
\begin{enumerate}
\item \label{U1scatter+} $\ds \lim_{n\to \infty} \frac{-t_{1,n}}{\lambda_{1,n}}=+\infty$.
\item \label{U1scatter-} $\ds \lim_{n\to \infty} \frac{-t_{1,n}}{\lambda_{1,n}}=-\infty$.
\item \label{U1compact} $\ds \forall n,\quad t_{1,n}=0$.
\end{enumerate}

Case \eqref{U1scatter+} is clearly excluded, as it would imply by Proposition \ref{P:lin_NL} that $u$ scatters for positive time, contradicting our assumptions.

Assume that \eqref{U1scatter-} holds. Then the nonlinear solution $U^1$ scatters for negative time. Precisely, by definition of $U^1$,
$$ \lim_{t\to -\infty}\left\|U^1_{\lin}(t)-U^1(t)\right\|_{\hdot}+\left\|\partial_t U^1_{\lin}(t)-\partial_t U^1(t)\right\|_{L^2}=0.$$
Furthermore, by Proposition \ref{P:lin_NL}, denoting as usual by $U_n^j$ the rescaled profiles (see Notation \ref{N:nonlinearprofiles}),
\begin{align}
\label{bad_dev_1}
u(0)&=\frac{1}{\lambda_{1,n}^{1/2}}U^1_{\lin}\left(\frac{-t_{1,n}-\tau_n}{\lambda_{1,n}},\frac{x}{\lambda_{1,n}}\right)+\sum_{j=2}^J U^j_n\left(-\tau_n,x\right)+w_n^J(-\tau_n)+r_n^J(-\tau_n),\\
\label{bad_dev_2}
\partial_t u(0)&=\frac{1}{\lambda_{1,n}^{3/2}}\partial_t U^1_{\lin}\left(\frac{-t_{1,n}-\tau_n}{\lambda_{1,n}},\frac{x}{\lambda_{1,n}}\right)+\sum_{j=2}^J \partial_t U^j_n\left(-\tau_n,x\right)+\partial_t w_n^J(-\tau_n)+\partial_t r_n^J(-\tau_n),
\end{align} 
Let $v_n(t)=S_{\lin}(t)\left(\lambda_{1,n}^{1/2}u\left(0,\lambda_{1,n}x\right),\lambda_{1,n}^{3/2}u\left(0,\lambda_{1,n}x\right)\right)$.
By orthogonality of the parameters $\{\lambda_{j,n}\}$, $\{t_{j,n}\}$, the developments \eqref{bad_dev_1}, \eqref{bad_dev_2} imply
$$ \left(v_n\left(\frac{t_{1,n}+\tau_n}{\lambda_{1,n}}\right),\partial_t v_n \left(\frac{t_{1,n}+\tau_n}{\lambda_{1,n}}\right)\right) \xrightharpoonup[n\to\infty]{} (U^1_{0},U^1_1) \text{ in }\hdot\times L^2,$$
since $\frac{t_{1,n}+\tau_n}{\lambda_{1,n}}\to +\infty$ this would imply $(U^1_0,U^1_1)=(0,0)$, a contradiction.

\EMPH{Step 3. Compact main profile}
It remains to consider case \eqref{U1compact}.
 By \eqref{mu=t},
\begin{multline}
\label{c0C0}
\int_{C_0 \tau_n\leq |x|}|\nabla u(\tau_n,x)|^2+|\partial_t u(\tau_n,x)|^2dx\\ \leq \frac{1}{2}\int|\nabla W|^2
\leq\int_{c_0 \tau_n\leq |x|}|\nabla u(\tau_n,x)|^2+|\partial_t u(\tau_n,x)|^2dx.
\end{multline}
This shows by assumption \eqref{bound_nabla_global},
$$
\int_{|x|\leq c_0\tau_n}|\nabla u(\tau_n,x)|^2+|\partial_t u(\tau_n,x)|^2dx\leq \frac{1}{2}\int|\nabla W|^2+\eta_1
$$
and thus by \eqref{bound_U1} (using that $t_{1,n}=0$),
$\lambda_{1,n}\approx\tau_n$. Extracting subsequences and rescaling $U^1$ we may assume that $\lambda_{1,n}=\tau_n$. Then by \eqref{CV_to_A},
\begin{equation}
\label{U1min}
U^1_0=A_0,\quad U^1_1=A_1.
\end{equation} 
We will show that $T^-(U^1)=-1$ and that 
\begin{equation}
\label{def_K_U1}
K=\lf\{\coltwo{(1+t)^{1/2}U^1\lf(t,(1+t)x\rg)}{(1+t)^{3/2}\partial_t U^1\lf(t,(1+t)x\rg)},\; t\in (-1,0]\rg\}
\end{equation}
has compact closure $\hdot \times L^2$. This type of self-similar solution is excluded by \cite[Section 6]{KeMe08}.
Let $\sigma\in(T^-(U_1),0)$. 
Then by Proposition \ref{P:lin_NL},
\begin{align}
\label{devu1}
u(\tau_n+\sigma\tau_n)&=\frac{1}{\tau_n^{1/2}}U^1\left(\sigma,\frac{x}{\tau_n}\right)+\sum_{j=2}^J U^j_n\left(\sigma\tau_n\right)+w_n^J(\sigma\tau_n)+r_n^J(\sigma\tau_n),\\
\label{devu2}
\partial_t u(\tau_n+\sigma\tau_n)&=\frac{1}{\tau_n^{3/2}}\partial_t U^1\left(\sigma,\frac{x}{\tau_n}\right)+\sum_{j=2}^J \partial_t U^j_n\left(\sigma\tau_n\right)+\partial_t w_n^J(\sigma\tau_n)+\partial_t r_n^J(\sigma\tau_n).
\end{align} 
Let
$$ Z^j_n(t)=\sum_{j=2}^J U^j_n(t)+w_n^J(t).$$
By Remark \ref{R:pythag}, we have for large $J$ and $n$,
$$ \lf\|\nabla_{t,x}U^1\lf(\sigma\rg)\rg\|^2_{L^2}+\lf\|\nabla_{t,x}Z_n^j(\sigma\tau_n)\rg\|_{L^2}^2\leq\|\nabla_{t,x}u(\tau_n+\sigma\tau_n)\|^2_{L^2}+\eta_1.
$$
By assumption \eqref{bound_nabla_global} and using \eqref{bound_U1}, we get
\begin{equation}
\label{bound_other_profiles}
\lf\|\nabla_{t,x} Z^j_n(\sigma\tau_n)\rg\|_{L^2}^2\leq \frac 13 \|\nabla W\|_{L^2}^2+2\eta_1.
\end{equation} 
By \eqref{mu=t} and the triangle inequality, we deduce  for large $J$ and $n$,
\begin{multline*}
\sqrt{\frac{1}{2}\int|\nabla W|^2}\leq\sqrt{\int_{c_0(1+\sigma) \tau_n\leq |x|}|\nabla_{t,x} u((1+\sigma)\tau_n,x)|^2dx}
\\
\leq \sqrt{\int_{c_0(1+\sigma)\tau_n\leq |x|}\frac{1}{\tau_n^{3}}\lf|\nabla_{t,x}U^1\lf(\sigma,\frac{x}{\tau_n}\rg)\rg|^2dx}+\sqrt{\int_{c_0 (1+\sigma)\tau_n\leq |x|}\lf|\nabla_{t,x}Z_n^j(\sigma\tau_n,x)\rg|^2dx}+\eta_1.
\end{multline*}
Thus by \eqref{bound_other_profiles}, and if $\eta_1$ is chosen small enough so that the left hand side inequality holds,
\begin{equation}
\label{self-similar}
(2\eta_1)^{1/2}\leq \sqrt{\frac{1}{2}\int|\nabla W|^2}-\sqrt{\frac 13 \int|\nabla W|^2+2\eta_1}-\eta_1\leq \sqrt{\int_{c_0(1+\sigma)\leq |x|} \lf|\nabla_{t,x} U^1(\sigma)\rg|^2}.
\end{equation} 
Using again assumption \eqref{bound_nabla_global}, we obtain
$$\forall \sigma \in (-1,0),\quad \int_{|x|\leq c_0(1+\sigma)} \lf|\nabla_{t,x} U^1(\sigma)\rg|^2\leq \int_{\RR^3} |\nabla W|^2-\eta_1.$$
In view of \eqref{limsup} in Theorem \ref{T:general_blowup}, we must have $T^-\lf(U^1\rg)\leq -1$. We cannot have $T^-(U^1)<-1$ because \eqref{devu1}, \eqref{devu2} with $\sigma=-1$ would give a nontrivial profile decomposition for $(u(0),\partial_t u(0))$, a contradiction. Thus $T^-(U^1)=-1$.

Next, note that by the development \eqref{devu1},\eqref{devu2}, we have
$$ \lf(\tau_n^{1/2}u((1+\sigma)\tau_n,\tau_n\cdot),\tau_n^{3/2}u((1+\sigma)\tau_n,\tau_n\cdot)\rg)\xrightharpoonup[n\to \infty]{}\lf(U^1(\sigma),\partial_t U^1(\sigma)\rg).$$
This shows that 
\begin{equation}
\label{U1inA}
\forall \sigma \in (-1,T_+(U^1)),\quad
\lf((1+\sigma)^{1/2}U^1(\sigma,(1+\sigma)\cdot),(1+\sigma)^{3/2}\partial_t U^1(\sigma,(1+\sigma)\cdot)\rg)\in \AAA.
\end{equation} 
We next show that $K$ defined by \eqref{def_K_U1}  
has compact closure in $\hdot\times L^2$. Indeed, let $t_n$ be a sequence that goes to $-1$ and assume after extraction that (weakly in $\hdot\times L^2$)
\begin{equation}
\label{weak_CV_U1}
\lf((1+t_n)^{1/2}U^1(t_n,(1+t_n)\cdot),(1+t_n)^{3/2}\partial_t U^1(t_n,(1+t_n)\cdot)\rg) \xrightharpoonup[]{n\to\infty} \lf(\tU_0,\tU_1\rg).
\end{equation}
Then by \eqref{U1inA} and the fact that $\AAA$ is closed for the weak topology, $\lf(\tU_0,\tU_1\rg)\in\AAA$ . In particular, using that $(U^1_0,U^1_1)=(A_0,A_1)$ has minimal energy in $\AAA$,
\begin{equation}
\label{ineg_energy}
0<
 E(U^1_0,U^1_1)\leq E\lf(\tU_0,\tU_1\rg).
\end{equation} 
We must show that \eqref{weak_CV_U1} is (at least for a subsequence) a strong convergence. For this, consider, after extraction, a profile decomposition for the sequence
$$ \left(U^1(t_n,x)-\frac{1}{(1+t_n)^{1/2}}\tU_0\left(\frac{x}{1+t_n}\right),\partial_t U^1(t_n,x)-\frac{1}{(1+t_n)^{3/2}}\tU_1\left(\frac{x}{1+t_n}\right) \right).$$ 
Denote the profiles by $V^j_{\lin}$, the parameters by $s_{j,n}$ and $\nu_{j,n}$ and the remainders by $\tw_n^J$. By the Pythagorean expansion of the energy
$$ E(U^1_0,U^1_1)=E(\tU_0,\tU_1)+\sum_{j=1}^J E\lf(V^j_{\lin}\lf(\frac{-s_{jn}}{\nu_{j,n}}\rg),\partial_t V^j_{\lin}\lf(\frac{-s_{jn}}{\nu_{j,n}}\rg)\rg)+E\lf(\tw_{0,n}^J,\tw_{1,n}^J\rg)+o_n(1).$$
By Claim \ref{C:variationnal}, all the energies are positive in this expansion.
By \eqref{ineg_energy}, $E(\tU_0,\tU_1)=E(U_0^1,U_1^1)$, and thus using Claim \ref{C:variationnal} again, $V^j_{\lin}=0$ for all $j\geq 1$ and $\lf\|\nabla \tw_{0,n}^J \rg\|_{L^2}+\lf\|\tw_{1,n}^J\rg\|_{L^2}$ tends to $0$ as $n$ tends to infinity, concluding the proof of the compactness of $\overline{K}$ in $\hdot\times L^2$ and yielding the desired contradiction. Note that in this last argument, we only needed the profile decomposition, for a fixed $J$, to show that the weak convergence \eqref{weak_CV_U1} and the inequality \eqref{ineg_energy} imply the strong convergence. The proof of Proposition \ref{P:norm_infinity} is complete.
\end{proof}

\section{Proof of the main result}
\label{S:radial}
In this section we show Theorem \ref{T:classification}. 

Assume that $N=3$ and that $u$ is a spherically symmetric type II blow-up solution such that:
\begin{equation}
\label{bound_nabla}
\sup_{\tau_0\leq t<1} \|\nabla u(t)\|_{L^2}^2+\|\partial_t u(t)\|_{L^2}^2\leq \|\nabla W\|_{L^2}^2+\eta_0,
\end{equation} 
The proof of Theorem \ref{T:classification} takes several steps. Consider the singular part $a$ of $u$ given by Definition \ref{D:regular_part}. In \S \ref{ss:compact_U1}, we show that a profile decomposition of a sequence $(a(\tau_n),\partial_t a(\tau_n))$, $\tau_n\to 1^-$ admits a large profile which is compact up to scaling. In \S \ref{SS:no_ss}, we show that, at least for a time sequence, the concentration is not self-similar, i.e that $u$ satisfy the assumptions of Section \ref{S:non_self_similar}. In \S \ref{SS:proof_compact}, we show that $a(t)$ is compact in the energy space up to a scaling parameter. In \S \ref{SS:proof_W} it is proven that the only limit as $t$ tends to $1$, of $a(t)$ up to scaling is $W$. We then conclude the proof of the theorem. 

\subsection{Compactness of the main profile}
\label{ss:compact_U1}
\begin{lemma}
\label{L:compactness_U1}
Assume that $N=3$ and that \eqref{bound_nabla} holds. Consider a sequence $\tau_n\to 1^-$, a profile decomposition $\left\{U^j\right\}$, $\left\{\lambda_{j,n}\right\}$, $\left\{t_{j,n}\right\}$ associated to $\left(a(\tau_n),\partial_ta(\tau_n)\right)$ and reorder the profiles (after extraction) so that \eqref{reorder} holds. Then all the profiles $U^j$, $j\geq 2$ scatter. Furthermore $U^1$ does not scatter for positive nor negative time,
\begin{equation} 
\label{U1}
\|U^1_{0}\|^2_{\hdot}+\|U^1_{1}\|^2_{\hdot}\geq \frac{2}{3}\|\nabla W\|_{L^2}^2,
\end{equation}
and the sequence $\left\{\frac{-t_{1,n}}{\lambda_{1,n}}\right\}_n$ is bounded.
\end{lemma}
In other words, the largest profile is compact up to modulation and we may assume that $t_{1,n}=0$ for all $n$.
\begin{proof}
The inequality \eqref{U1} follows from Lemma \ref{L:profiles}. The assumption \eqref{bound_nabla} implies that for all $j\geq 2$, $\|U_1^{j}\|_{L^2}^2+\|\nabla U_0^{j}\|_{L^2}^2\leq \frac 13 \|\nabla W\|_{L^2}^2+\eta_0$. Thus all nonlinear profiles $U^j$, $j\geq 2$, scatter both forward and backward in time. To conclude the proof, it is sufficient to show that $U^1$ does not scatter forward nor backward in time, which would imply that $\left\{\frac{-t_{1,n}}{\lambda_{1,n}}\right\}_n$ is bounded.
Assume that $U^1$ is globally defined and scatters forward in time. Then, by, Proposition \ref{P:lin_NL}, $u$ is globally defined and scatters forward in time, a contradiction.
It remains to exclude the case when $U^1$ is globally defined and scatters backward in time. 
By Proposition \ref{P:lin_NL} again, we obtain that for $t<0$, 
$$ u(\tau_n+t,x)=v(\tau_n+t,x)+\sum_{j=1}^J U^j_n\left(t,x\right)+w^J_{n}(t,x)+r^J_n(t,x),$$
where 
\begin{equation*}
\lim_{J\rightarrow +\infty} \limsup_{n\rightarrow +\infty} \|r^J_n\|_{S(-\infty,0)}+\sup_{t\in (-\infty,0)} \left(\|\nabla r^J_n(t)\|_{L^2}+\|\partial_t r^J_n(t)\|_{L^2}\right)=0.
\end{equation*}
The solution $U^1$ scatters backward, but not forward in time. By \cite{KeMe08}, this implies that $E(U^1_0,U^1_1)\geq E(W,0)$. As a consequence, for all $t$ in the domain existence of $U^1$,
\begin{equation}
\label{nabla_grand}
\|\nabla U^1(t)\|_{L^2}^2+\|\partial_t U^1(t)\|_{L^2}^2\geq 2E(U^1_0,U^1_1)\geq 2E(W,0)=\frac{2}{3}\|\nabla W\|_{L^2}^2.
\end{equation}
Let $t_0\in (\tau_0,1)$, where $\tau_0$ is defined in \eqref{bound_nabla}. Taking $t=t_0-\tau_n<0$ in the preceding decomposition, we obtain that for large $n$,
\begin{align*}
u(t_0,x)=\frac{1}{\lambda_{1,n}^{1/2}} U^1\left(\frac{t_0-\tau_n-t_{1,n}}{\lambda_{1,n}},\frac{x-x_{1,n}}{\lambda_{1,n}}\right)+R_{0,n}(x),\\
\partial_t u(t_0,x)=\frac{1}{\lambda_{1,n}^{3/2}} \partial_t U^1\left(\frac{t_0-\tau_n-t_{1,n}}{\lambda_{1,n}},\frac{x-x_{1,n}}{\lambda_{1,n}}\right)+R_{1,n}(x),
\end{align*}
where by Pythagorean expansion, $\|\nabla R_{0,n}\|_{L^2}^2+\|R_{1,n}\|_{L^2}^2\leq \frac{1}{3}\|\nabla W\|_{L^2}^2+\eta_0$. 
By \eqref{nabla_grand},
we get that $\left\{\big(u(t_0),\partial_t u(t_0)\big)\right\}_n$, considered as a sequence  in $n$, admits a nontrivial profile decomposition (recall that $\lambda_{1,n}\to 0$), a contradiction.
The proof is complete.
\end{proof}

\subsection{Existence of a sequence avoiding self-similar blow-up}
\label{SS:no_ss}
\begin{prop}
\label{P:noss}
Assume $N=3$ and let $u$ be a radial solution satisfying \eqref{bound_nabla}. 
Then there exists $\{\tau_n\}_{n}$, $\{\mu_n\}_n$ with 
$$ \tau_n\to 1^-,\quad 0<\mu_n\ll 1-\tau_n\text{ as }n\to\infty$$
such that
$$ \lim_{n\to \infty} \int_{|x|\geq \mu_n}(\partial_t a(\tau_n,x))^2+|\nabla a(\tau_n,x)|^2+\frac{1}{|x|^2}(a(\tau_n,x))^2dx=0$$
\end{prop}
\begin{corol}
\label{C:energy}
$$\lim_{t\to 1^-} E(a(t),\partial_t a(t))=E(W,0).$$
\end{corol}
\begin{proof}[Proof of Corollary \ref{C:energy}]
The result of Proposition \ref{P:noss} implies by Proposition \ref{P:sum_of_W} that (replacing $u$ by $-u$ if necessary), there exists a sequence $\tau_n\to 1^-$, a sequence $\lambda_n\to 0$ such that
\begin{align}
a(\tau_n,x)&=\frac{1}{\lambda_n^{1/2}}W\left(\frac{x}{\lambda_n}\right)+w_{0,n}\\
\partial_t a(\tau_n,x)&=o(1)\text{ in }L^2\text{ as }n\to\infty,
\end{align}
where, denoting by $w_n$ the solution of \eqref{lin_wave} with initial condition $(w_{0,n},0)$,
$$ \lim_{n\to \infty} \|w_n\|_{S(-\infty,+\infty)}=0.$$
\EMPH{Step 1} We first show 
\begin{equation}
\label{no_dispersive}
\lim_{n\to\infty} \|w_{0,n}\|_{\hdot}=0.
\end{equation} 
Let us mention that this step still works, with a small refinement, replacing the assumption \eqref{bound_nabla} by the more general \eqref{bound_nabla'}.

Assume that \eqref{no_dispersive} does not hold. Extracting a subsequence in $n$, we can assume that there exists $\eps_0>0$ and, for all $n$, $r_n>0$ such that
$$ \int_{|x|\geq r_n} |\nabla w_{0,n}(x)|^2\,dx\geq \eps_0.$$
Then, arguing as in the proof of Proposition \ref{P:ext_profile} (see \eqref{minor_wnJ}), we get that for large $n$,
$$ \int_{r_n}^{+\infty} |\partial_r (rw_{0,n})(r)|^2\,dr\geq \frac{\eps_0}{2}.$$
Next, the fact that $w_n(t)=w_n(-t)$ and Lemma \ref{L:linear_behavior} imply that for large $n$, for all $T>0$
\begin{equation}
\label{wn_large_r}
\int_{|x|\geq r_n+T\lambda_n} |\nabla_{t,x}w_n(-T\lambda_n,x)|^2dx\geq \frac{\eps_0}{4}.
\end{equation} 
By Proposition \ref{P:lin_NL}, we have
\begin{equation}
\label{other_dev}
\left\{
\begin{aligned} a(\tau_n-T\lambda_n)&=\frac{1}{\lambda_n^{1/2}}W\left(\frac{x}{\lambda_n}\right)+w_n(-T\lambda_n)+o_n(1)\text{ in }\hdot\\
\partial_ta(\tau_n-T\lambda_n)&=\partial_t w_n(-T\lambda_n)+o_n(1)\text{ in }L^2.
\end{aligned}\right.
\end{equation}
Combining with \eqref{wn_large_r} we get that there exists an increasing sequence $\{n_k\}$ such that $\tau_{n_k}-k\lambda_{n_k}\geq 0$ and
$$ \int_{|x|\geq r_{n_k}+k\lambda_{n_k}} |\nabla_{t,x}a(\tau_{n_k}-k\lambda_{n_k},x)|^2\,dx\geq \frac{\eps_0}{8}.$$
In view of \eqref{other_dev}, this contradicts Proposition \ref{P:ext_profile} (here $\rho_{1,n}=\lambda_{n}$). Step 1 is complete.

\EMPH{Step 2}
By Step 1,
$$ \lim_{n\to \infty} E(a(\tau_n),\partial_ta(\tau_n))=E(W,0).$$
Note that 
$$E(u(t),\partial_tu(t))=E(v(t),\partial_t v(t))+E(a(t),\partial_t a(t))+o(1)\text{ as }t\to 1^-,$$
which shows by conservation of the energy for $u$ and $v$ that $E(a(t),\partial_t a(t))$ has a limit as $t\to 1^-$, concluding the proof of Corollary \ref{C:energy}.
\end{proof}
\begin{proof}[Proof of Proposition \ref{P:noss}]
We argue by contradiction. By Hardy's inequality
$$ \int_{|x|\geq R} \frac{1}{|x|^2}(a(t,x))^2dx\leq C\int_{|x|\geq R}|\nabla a(t,x)|^2dx,$$
so that we only need to show that there exist sequences $\mu_n$ and $\tau_n$ as in the proposition such that
$$ \lim_{n\to \infty} \int_{|x|\geq \mu_n}(\partial_t a(\tau_n,x))^2+|\nabla a(\tau_n,x)|^2dx=0$$
If this does not hold, there exists $\alpha>0$ and $\eps_0>0$ such that
\begin{equation}
\label{absurd_ss}
\forall t\in (0,1),\quad \int_{|x|\geq \alpha (1-t)}|\partial_t a(t,x)|^2+|\nabla a(t,x)|^2dx\geq \eps_0.
\end{equation} 

\EMPH{Step 1} We first show that there exists $\beta>0$ such that\footnote{we could replace $\frac{2}{3}\|\nabla W\|_{L^2}^2$ by $\|\nabla W\|_{L^2}^2-C\eta_0$ for some large positive constant $C$, where $\eta_0$ is given by \eqref{bound_nabla}.}
\begin{equation}
\label{big_ss}
\liminf_{t\to 1^-}\int_{|x|\geq \beta (1-t)}|\partial_t a(t,x)|^2+|\nabla a(t,x)|^2dx\geq \frac{2}{3}\|\nabla W\|_{L^2}^2.
\end{equation} 
Indeed, assume that \eqref{big_ss} does not hold, i.e. that there exists sequences $\tau_n\to 1^-$, $\beta_n\to 0^+$ such that
\begin{equation}
\label{absurd_2}
\int_{|x|\geq \beta_n(1-\tau_n)}|\partial_t a(\tau_n,x)|^2+|\nabla a(\tau_n,x)|^2dx\leq \frac{2}{3}\|\nabla W\|_{L^2}^2-\eps_1.
\end{equation} 
After extraction, consider a profile decomposition $\{U^j\}$, $\{t_{j,n}, \lambda_{j,n}\}$ for $\left\{(a(\tau_n),\partial_t a(\tau_n))\right\}_n$. Reordering the profiles, we assume \eqref{reorder}, i.e that $U^1$ is the largest profile in the energy space. By Lemma \ref{L:compactness_U1}, we may assume that $t_{1,n}=0$, and the norm of $(U^1_0,U^1_1)$ in the energy space is bounded from below (see \eqref{U1}).

Let $\eps_2>0$ to be specified later. By Proposition \ref{P:scatteringKeMe}, there exists $T\in  (0,T_+(U^1))$ such that 
$$ \|\nabla U^1(T)\|_{L^2}^2+\|\partial_t U^1(T)\|_{L^2}^2\geq \|\nabla W\|_{L^2}^2-\eps_2.$$
Then by Proposition \ref{P:lin_NL}
\begin{gather}
\label{dev_aa}
a(\tau_n+\lambda_{1,n}T)=\sum_{j=1}^JU^j_n\left(\lambda_{1,n} T\right)+w_{n}^J(\lambda_{1,n}T)+r_{n}^J(\lambda_{1,n}T)\\
\label{dev_ab}
\partial_t a(\tau_n+\lambda_{1,n}T)=\sum_{j=1}^J\partial_t U^j_n\left(\lambda_{1,n} T\right)+\partial_t w_{n}^J(\lambda_{1,n}T)+\partial_t r_{n}^J(\lambda_{1,n} T).
\end{gather} 
The rescaled profiles $U^j_n$ are defined as usual (see Notation \ref{N:nonlinearprofiles}).
Note that 
$$ \|\nabla U^1(\lambda_{1,n}T)\|_{L^2}^2+\|\partial_t U^1(\lambda_{1,n}T)\|_{L^2}^2=\|\nabla U^1_n(T)\|_{L^2}^2+\|\partial_t U^1_n(T)\|_{L^2}^2\geq \|\nabla W\|_{L^2}^2-\eps_2.$$
Combining with \eqref{bound_nabla}, \eqref{dev_aa}, \eqref{dev_ab} and the orthogonality of the parameters, we get 
\begin{equation*}
\sum_{j=2}^J \lf\|\nabla_{t,x} U^j\lf(\frac{\lambda_{1,n} T-t_{j,n}}{\lambda_{j,n}}\rg)\rg\|_{L^2}^2+\lf\|\nabla_{t,x}w_n^J(\lambda_{1,n}T)\rg\|_{L^2}^2\leq \eta_0+\eps_2.
\end{equation*} 
And thus using the conservation of the energy
\begin{equation*}
\sum_{j=2}^J E(U^j_0,U_j^1)\leq \frac{1}{2}(\eta_0+\eps_2).
\end{equation*} 
Take $\eta_0$ and $\eps_2$ so small that $\frac{1}{2}(\eta_0+\eps_2)\leq \frac{1}{3}\delta_1^2$, where $\delta_1=\delta_1(2\|\nabla W\|_{L^2}^2)$ is given by Corollary \ref{C:small_data_behavior'}. Then
$U^1$ is the only one large profile, i.e., with the notations of \S \ref{SS:No_exterior}, $J_0=1$. Assume that $\lambda_{1,n}=o_n(1-\tau_n)$. Then by Proposition \ref{P:ext_profile} we would obtain
\begin{equation*}
\lim_{R\rightarrow +\infty} \limsup_{n\rightarrow+\infty} \int_{|x|\geq R\lambda_{1,n}} \left(|\nabla a(\tau_n)|^2+(\partial_t a(\tau_n))^2\right)dx=0.
\end{equation*} 
a contradiction with \eqref{absurd_ss}. Thus 
$$ (1-\tau_n)\approx \lambda_{1,n}.$$
Consider a sequence $\{\tbeta_n\}$ such that 
$$ \beta_n\ll\tbeta_n\ll 1$$
Let $\chi\in C^{\infty}(\RR^3)$ such that $\chi(x)=1$ if $|x|\geq 2$ and $\chi(x)=0$ is $|x|\leq 1$. By Remark \ref{R:pythag}
\begin{multline*}
\int_{|x|\geq \beta_n(1-\tau_n)} |\nabla a(\tau_n,x)|^2+|\partial_t a(\tau_n,x)|^2\geq \int\chi\left(\frac{x}{\tbeta_n(1-\tau_n)}\right)\left(|\nabla a(\tau_n,x)|^2+|\partial_t a(\tau_n,x)|^2\right)\\
\geq \int \chi\left(\frac{\lambda_{1,n}y}{\tbeta_n(1-\tau_n)}\right)\left(|\nabla U^1_0(y)|^2+|U^1_1(y)|^2\right)dy\\
\underset{n\to\infty}{\longrightarrow} \int \left(|\nabla U^1_0(y)|^2+|U^1_1(y)|^2\right)dy
\geq  \frac{2}{3}\|\nabla W\|^2_{L^2}.
\end{multline*}
This contradicts \eqref{absurd_2} and concludes Step 1.

\EMPH{Step 2. End of the argument}
Let, for $t\in [0,1)$,
\begin{equation}
\label{def_mu1}
\mu(t)=\inf\left\{\mu\;:\; \int_{|x|\leq \mu}|\partial_t a(t)|^2+|\nabla a(t)|^2\geq \frac 25 \int |\nabla W|^2\right\}.
\end{equation} 
By Step 1 and assumption \eqref{bound_nabla},
\begin{equation}
\label{estim_mu}
\beta(1-t)\leq \mu(t)\leq 1-t.
\end{equation} 
Take any sequence $\tau_n\to 1^-$ such that $(a(\tau_n),\partial_t a(\tau_n))$ admits a profile decomposition. By Lemma \ref{L:compactness_U1} and Step 1, we may assume, after extraction,
$$ \lf\|\nabla U^1_0\rg\|_{L^2}^2+\lf\|U^1_0\rg\|_{L^2}^2\geq \frac{2}{3}\|\nabla W\|_{L^2}^2,
\quad t_{1,n}=0,\quad \lambda_{1,n}=\mu(\tau_n)\approx 1-\tau_n.$$
Furthermore, the solution $U^1$ does not scatter forward nor backward in time. Let $\eps_3=\frac{\beta^2}{2}$, where $\beta$ is given by \eqref{estim_mu}. By Proposition \ref{P:norm_infinity} (if $T_-(U^1)=-\infty$) or Section \ref{S:reg_sing} (if $T_-(U^1)\in(-\infty,0)$), there exists $-\theta\in (T_-(U^1),0)$ such that
\begin{equation}
\label{compression_time+}
\int_{|x|\leq \eps_3\theta} |\partial_t U^1(-\theta)|^2+|\nabla U^1(-\theta)|^2\geq \frac{1}{2}\int |\nabla W|^2.
\end{equation} 
Let us show that for large $n$
\begin{equation}
\label{concentration_E}
\mu(\tau_n-\theta\mu(\tau_n))\leq \eps_3\theta\mu(\tau_n).
\end{equation} 
If this holds, we would get by \eqref{estim_mu} 
$$ \beta^2\theta(1-\tau_n)\leq \beta\theta\mu(\tau_n)\leq \beta(1-\tau_n+\theta\mu(\tau_n))\leq \mu\left(\tau_n-\theta\mu(\tau_n)\right)\leq \eps_3\theta\mu(\tau_n)\leq\frac{\beta^2}{2}\theta(1-\tau_n),$$
a contradiction. The inequality \eqref{concentration_E} is equivalent to the following 
\begin{equation}
\label{concentration_E1}
\int_{|x|\leq \eps_3\theta\mu(\tau_n)} |\nabla_{t,x} a(\tau_n-\theta\mu(\tau_n))|^2\geq \frac{2}{5}\int |\nabla W|^2.
\end{equation} 
We have, denoting by  $\theta_{j,n}=-\theta\mu(\tau_n)-t_{j,n}$,
\begin{align*}
 u(\tau_n-\theta\mu(\tau_n))&=v(\tau_n-\theta\mu(\tau_n))+\frac{1}{\mu^{\frac 12}(\tau_n)}U^1\left(-\theta,\frac{x}{\mu(\tau_n)}\right)\\
&\qquad +\sum_{j=2}^J \frac{1}{\lambda_{jn}^{1/2}}U^j\left(\frac{\theta_{j,n}}{\lambda_{j,n}},\frac{x}{\lambda_{j,n}}\right)+w_{n}^J(-\theta\mu(\tau_n))+r_{n}^J(-\theta\mu(\tau_n))\\
 \partial_t u(\tau_n-\theta\mu(\tau_n))&=\partial_t v(\tau_n-\theta\mu(\tau_n))+\frac{1}{\mu^{\frac 32}(\tau_n)}\partial_t U^1\left(-\theta,\frac{x}{\mu(\tau_n)}\right)\\
&\qquad +\sum_{j=2}^J \frac{1}{\lambda_{j,n}^{3/2}}\partial_t U^j\left(\frac{\theta_{j,n}}{\lambda_{j,n}} ,\frac{x}{\lambda_{j,n}}\right)+\partial_t w_{n}^J(-\theta\mu(\tau_n))+\partial_t r_{n}^J(-\theta\mu(\tau_n)),
\end{align*}
where $r_{n}^J$ satisfies 
$$ \lim_{J\rightarrow +\infty} \limsup_{n\rightarrow+\infty} \|\nabla r_{n}^J(\theta\mu(\tau_n))\|_{L^2}+\|\partial_t r_{n}^J(\theta\mu(\tau_n)\|_{L^2}=0.$$
Let $J_0$ such that for all $J\geq J_0$, for large $n$, $\lf\|\partial_{t,x}r_n^J\rg\|^2\leq \frac{1}{40}\int |\nabla W|^2$.
Let $\psi\in C_0^{\infty}(\RR^3)$ with $\psi(x)=1$ for  $|x|\leq 1$ and $\psi(x)=0$ for $|x|\geq 2$. Then by Remark \ref{R:pythag},
$$ \int\psi\lf(\frac{x}{\eps_3\theta\mu(\tau_n)}\rg) |\nabla_{t,x} a(\tau_n-\theta\mu(\tau_n))|^2\geq \int\psi\lf(\frac{x}{\eps_3\theta}\rg)|\nabla_{t,x} U^1(-\theta)|^2\geq\frac{2}{5}\int |\nabla W|^2,$$
hence \eqref{concentration_E}. The proof is complete.
\end{proof}

\subsection{Compactness of the singular part}
\label{SS:proof_compact}
\begin{prop}
\label{P:compactness}
Under the assumptions of Theorem \ref{T:classification} (in particular $N=3$ and $u$ is spherically symmetric),
$a$ is compact in the energy space up to a scaling parameter: there exists a continuous function $\lambda(t)$, $t\in (0,1)$ such that the closure of
$$K=\left\{\left(\lambda^{1/2}(t)a\left(t,\lambda(t)x\right),\lambda^{3/2}(t)\partial_t a\left(t,\lambda(t)x\right)\right),\;t\in (0,1)\right\}$$
is compact in $\hdot\times L^2$.
\end{prop}
\begin{proof}
It is sufficient to show that for any time sequence $\tau_n \xrightarrow{\ssstyle <}1$, there exists a subsequence of $\tau_n$ and a sequence $\lambda_n$ such that $\left(\lambda^{1/2}_na\left(\tau_n,\lambda_nx\right),\lambda^{3/2}_n\partial_t a\left(\tau_n,\lambda_nx\right)\right)$ converges in $\hdot\times L^2$.

Let $\tau_n \xrightarrow{\ssstyle <}1$. After extraction of a subsequence (in $n$), assume that $(a(\tau_n),\partial_t a(\tau_n))$ as a profile decomposition with profiles $U^j_{\lin}$ and parameters $\lambda_{j,n}$, $t_{j,n}$. Let $U^1$ be the largest profile. By Lemma \ref{L:compactness_U1}, $\|\nabla U_0^1\|_{L^2}+\|U_1^1\|_{L^2}\geq \frac{2}{3}\|\nabla W\|_{L^2}^2$. By \eqref{bound_nabla} and the Pythagorean expansions \eqref{pythagore1a} and \eqref{pythagore1b}, we get
$$\|\nabla w_n^J(\tau_n)\|_{L^2}^2\leq \frac{1}{3}\|\nabla W\|_{L^2}^2+\eta_0$$
and
$$\forall j\geq 2, \quad \|\nabla U^j_0\|_{L^2}^2+\|U^j_1\|_{L^2}^2\leq \frac{1}{3}\|\nabla W\|_{L^2}^2+\eta_0.$$
This implies that the energies of $U^j$, $j\geq 2$ and of $w_n^J$ are all positive (see Claim \ref{C:variationnal}).
We distinguish three cases:
\begin{itemize}
\item If $E(U^1,\partial_t U^1)\geq E(W,0)$, then by Corollary \ref{C:energy} and the Pythagorean expansion of the energy (using that all energy are positive), we obtain immediately that $E(U^1,\partial_tU^1)=E(W,0)$, that there are no nonzero other profile and that $(w_{0,n}^J,w_{1,n}^J)$ tends to $0$ as $n\to \infty$, hence the compactness property.
\item If $E(U^1,\partial_t U^1)<E(W,0)$, and $\|\nabla U^1_0\|_{L^2}^2+\|U^1_1\|_{L^2}^2<\|\nabla W\|^2_{L^2}$, the profile $U^1$ scatters yielding immediately a contradiction. 
\item If $E(U^1,\partial_t U^1)<E(W,0)$, and $\|\nabla U^1_0\|_{L^2}^2+\|U^1_1\|_{L^2}^2>\|\nabla W\|^2_{L^2}$. The nonlinear solution $U^1$ blows up in both time directions. By Proposition \ref{P:lin_NL}, $U^1$ is a type II blow-up solution of \eqref{CP} such that $E(U^1,\partial_t U^1)<E(W,0)$. 
Furthermore, as $(a,\partial_t a)$ converges weakly to $0$ and $(v,\partial_t v)$ converges strongly in $\hdot\times L^2$ as $t\to 1$, we have
$$ \int |\nabla_{t,x} u(t,x)|^2=\int |\nabla_{t,x} a(t,x)|^2+\int |\nabla_{t,x} v(t,x)|^2+\underset{t\to 1^-}{o(1)}.$$
Thus $U^1$ also satisfies \eqref{bound_nabla}, which shows that $U^1$ contradicts Corollary \ref{C:energy}.
\end{itemize}
The proof if complete.
\end{proof}

\subsection{Convergence to the stationary solution up to the scaling}
\label{SS:proof_W}
In this section we conclude the proof of Theorem \ref{T:classification}.
Consider a solution $u$ of \eqref{CP} satisfying the assumptions of Theorem \ref{T:classification}. By Corollary \ref{C:energy},
\begin{equation}
\label{lim_a}
 \lim_{\substack{t\to 1^-}}E(a(t),\partial_t a(t))=E(W,0).
\end{equation} 
By Proposition \ref{P:compactness}, there exists $\lambda(t)$ such that the closure of 
$$ K=\left\{\left(\lambda^{1/2}(t)a\left(t,\lambda(t)x\right),\lambda^{3/2}(t)\partial_t a\left(t,\lambda(t)x\right)\right),\;t\in (0,1)\right\}$$
is compact in $\hdot\times L^2$. The following result is classical in this setting.
\begin{lemma}
\label{L:compactness_profile}
Let $\tau_n$ be a sequence that tends to $1$, and such that
$$ \left(\lambda^{1/2}(\tau_n)a\left(\tau_n,\lambda(\tau_n)x\right),\lambda^{3/2}(\tau_n)\partial_t a\left(\tau_n,\lambda(\tau_n)x\right)\right)\xrightarrow[n\to\infty]{} (U_0,U_1)$$ 
in $\hdot\times L^2$. Consider the solution $U$ of \eqref{CP} such that
$$ U_{\restriction t=0}=U_0,\quad \partial_t U_{\restriction t=0}=U_1.$$
Then there exists a continuous function $\tlambda$ defined on $(T_-(U),T_+(U))$ such that
$$ \tK=\left\{\coltwo{\tlambda(t)^{1/2}U\left(t,\tlambda(t)x\right)}{\tlambda(t)^{3/2}\partial_t U\left(t,\tlambda(t)x\right)},\quad t\in (T_-(U),T_+(U))\right\}$$
has compact closure in $\hdot\times L^2$.
\end{lemma}
\begin{proof}[Sketch of proof]
 We have
\begin{align*}
 u(\tau_n,x)&=v(\tau_n,x)+\frac{1}{\lambda(\tau_n)^{1/2}}U_0\left(\frac{x}{\lambda(\tau_n)}\right)+o_n(1)&&\text{ in }\hdot\\
\partial_tu(\tau_n,x)&=\partial_tv(\tau_n,x)+\frac{1}{\lambda(\tau_n)^{3/2}}U_1\left(\frac{x}{\lambda(\tau_n)}\right)+o_n(1)&&\text{ in }L^2.
\end{align*}
Let $T\in \left(T_-(U),T_+(U)\right)$. By Proposition \ref{P:lin_NL},
\begin{align*}
 u(\tau_n+\lambda(\tau_n)T,x)&=v(\tau_n+\lambda(\tau_n)T,x)+\frac{1}{\lambda(\tau_n)^{1/2}}U\left(T,\frac{x}{\lambda(\tau_n)}\right)+o_n(1)&&\text{ in }\hdot\\
\partial_tu(\tau_n+\lambda(\tau_n)T,x)&=\partial_tv(\tau_n+\lambda(\tau_n)T,x)+\frac{1}{\lambda(\tau_n)^{3/2}}\partial_t U\left(T,\frac{x}{\lambda(\tau_n)}\right)+o_n(1)&&\text{ in }L^2.
\end{align*}
Letting $\sigma_n=\tau_n+\lambda(\tau_n)T$, we get
\begin{align*}
 \lambda(\sigma_n)^{1/2}a(\sigma_n,\lambda(\sigma_n)x)&=\left(\frac{\lambda(\sigma_n)}{\lambda(\tau_n)}\right)^{1/2}U\left(T,\frac{\lambda(\sigma_n)}{\lambda(\tau_n)}x\right)+o_n(1)&&\text{ in }\hdot\\
\lambda(\sigma_n)^{3/2}\partial_ta(\sigma_n,x)&=\left(\frac{\lambda(\sigma_n)}{\lambda(\tau_n)}\right)^{3/2}\partial_t U\left(T,\frac{\lambda(\sigma_n)}{\lambda(\tau_n)}x\right)+o_n(1)&&\text{ in }L^2.
\end{align*}
Extracting subsequences, we obtain by compactness of $\overline{K}$ that there exists $(V_0,V_1)\in \overline{K}$ such that
$$ \lim_{n\to \infty}\left(\left(\frac{\lambda(\sigma_n)}{\lambda(\tau_n)}\right)^{1/2}U\left(T,\frac{\lambda(\sigma_n)}{\lambda(\tau_n)}x\right),\left(\frac{\lambda(\sigma_n)}{\lambda(\tau_n)}\right)^{3/2}\partial_t U\left(T,\frac{\lambda(\sigma_n)}{\lambda(\tau_n)}x\right)\right)=(V_0,V_1)\text{ in }\hdot\times L^2.$$
This shows that $\frac{\lambda(\sigma_n)}{\lambda(\tau_n)}$ has a limit $\tilde{\lambda}(T)\in (0,+\infty)$ (by conservation of the energy $(0,0)\notin \overline{K}$) and that
$$ \left(\tilde{\lambda}(T)^{1/2}U\left(T,\tilde{\lambda}(T)x\right),\tilde{\lambda}(T)^{3/2}\partial_t U\left(T,\tilde{\lambda}(T)x\right)\right)\in \overline{K}.$$
The proof is complete, up to the proof of the known fact that the function $T\mapsto \tilde{\lambda}(T)$ may be taken continuous for which we refer to \cite[Remark 5.4]{KeMe06}.
\end{proof}


We next prove Theorem \ref{T:classification}.

\EMPH{Step 1. Convergence to $W$ for sequences}
Let $\{t_n\}_n$ be a sequence in $(0,1)$ such that $t_n\to 1$ and 
\begin{equation*}
\lim_{n\to \infty} (\lambda^{1/2}(t_n) a(t_n,\lambda(t_n)x),\lambda(t_n)^{3/2}\partial_t a(t_n,\lambda(t_n)x))=(U_0,U_1)\text{ in }\hdot\times L^2.
\end{equation*}
In this step we show that for some $\lambda_0>0$ and some sign $+$ or $-$, $(U_0,U_1)=\pm \lf(\lambda_0^{1/2}W(\lambda_0\cdot),0\rg)$. 

 Let $U$ be the solution of \eqref{CP} with initial condition $(U_0,U_1)$. By Lemma \ref{L:compactness_profile}, $U$ is compact up to scaling. By Theorem \ref{T:compact=W},
$U=W$ up to the symmetries, concluding step 1.
%
%

\EMPH{Step 2. Estimate on the scaling parameter} 
Let 
$$\lambda_1(t)=\inf\left\{\mu>0\;:\; \int_{|x|\leq \mu}|\nabla u(t,x)-\nabla v(t,x)|^2dx\geq \int_{|x|\geq 1} |\nabla W|^2dx\right\}.$$
By Step 1, $\int|\nabla a(t,x)|^2dx\to \int |\nabla W|^2$ as $t\to 1$, which shows that $\lambda_1(t)$ is well-defined for $t<1$, close to $1$. Consider a sequence $t_n \xrightarrow{\ssstyle <} 1$. By Step $1$, for $\iota_0=-1$ or $+1$ and some sequence of positive numbers $\{\lambda_n\}_n$,
$$ a(t_n,x)=\iota_0\frac{1}{\lambda_n^{1/2}}W\left(\frac{x}{\lambda_n}\right)+o_n(1)\text{ in }\hdot.$$
Thus if $\mu>0$, 
$$ \int_{|x|\leq \mu} |\nabla a(t_n,x)|^2dx=\int_{|x|\leq \mu} \frac{1}{\lambda_n^3}\left|\nabla W\left(\frac{x}{\lambda_n}\right)\right|^2+o_n(1)=\int_{|y|\leq \mu/\lambda_n}|\nabla W(y)|^2dy+o_n(1),$$
which shows that
\begin{equation*}
\lim_{n\to \infty}\frac{\lambda_n}{\lambda_1(t_n)}= 1.
 \end{equation*}
Thus
\begin{equation}
\label{good_param}
 a(t_n,x)=\iota_0\frac{1}{\lambda_1(t_n)^{1/2}}W\left(\frac{x}{\lambda_1(t_n)}\right)+o_n(1)\text{ in }\hdot.
\end{equation} 

\EMPH{Step 3. Choice of the sign}
Let 
$$f(t)=\int \nabla a(t,x)\cdot\frac{1}{\lambda_1(t)^{1/2}}\nabla W\left(\frac{x}{\lambda_1(t)}\right)dx.$$
Then by Step 2, for each sequence $t_n \xrightarrow{\ssstyle <} 1$, there exists a subsequence such that $f(t_n)\to \pm \int |\nabla W|^2$. As $f$ is a continuous function, the intermediate value theorem implies that the value must be the same for all the sequences $\{t_n\}$. Changing $u$ into $-u$ if necessary, we can assume
$$ \lim_{t\to 1^-} f(t)=\int |\nabla W|^2.$$
By Step 2, we get that for all sequences $\{t_n\}$, 
$$ u(t_n,x)=v(t_n,x)+\frac{1}{\lambda_1(t_n)^{1/2}}W\left(\frac{x}{\lambda_1(t_n)}\right)+o_n(1),$$
which concludes the proof of the development \eqref{dev_u}.

\EMPH{Step 4. Estimate on $\lambda_1$}
Recalling that $u-v$ is supported in the cone $\{|x|\leq 1-t\}$, we get, for $t$ close to $1$
\begin{multline*}
0=\int_{|x|\geq 1-t}|\nabla u(t)-\nabla v(t)|^2dx=\int_{|x|\geq 1-t}\frac{1}{\lambda_1^{3}(t)}\lf|\nabla W\lf(\frac{x}{\lambda_1(t)}\rg)\rg|^2dx+\underset{t\to 1^-}{o(1)}\\
=\int_{|y|\geq \frac{1-t}{\lambda_1(t)}}|\nabla W(y)|^2dy+\underset{t\to 1^-}{o(1)},
\end{multline*}
which shows that $\frac{1-t}{\lambda_1(t)}\to +\infty$, concluding the proof of Theorem \ref{T:classification}.
\qed

\appendix

\section{Properties of profiles}
\label{A:ortho}
In this appendix we prove a pseudo-orthogonality property (Claim \ref{C:ortho}) and Claim \ref{C:dispersive_tw}.
\begin{claim}
\label{C:ortho}
Assume that $N\geq 3$ is odd.
Let $\{w_n\}$ be a sequence of finite energy solutions of the linear wave equation \eqref{lin_wave}, bounded in the energy space and $U$ be a finite energy solution of \eqref{lin_wave}. Consider real sequences $\{\lambda_n\}$, $\{\mu_n\}$, $\{t_n\}$, $\{\theta_n\}$ with $\lambda_n>0$, $\mu_n>0$. Assume that 
\begin{equation}
\label{weakCVwn}
\lambda_n^{N/2}\nabla_{t,x}w_n(t_n,\lambda_n \cdot) \xrightharpoonup[n\to \infty]{}0 \text{ in }L^2(\RR^{N+1}).
\end{equation}
Then, if $\varphi=1$, or if $\varphi$ is a radial, continuous, compactly supported function on $\RR^N$ and such that $\varphi(r)=1$ if $r$ is small, there exist subsequences such that
\begin{equation}
\label{orthogonality}
\lim_{n\to +\infty} \int \varphi\lf(\frac{|x|}{\mu_n}\rg)\nabla_{t,x}w_n(\theta_n,x)\cdot\frac{1}{\lambda_n^{N/2}}\nabla_{t,x}U\lf(\frac{\theta_n-t_n}{\lambda_n},\frac{x}{\lambda_n}\rg)dx=0,
\end{equation} 
and
\begin{equation}
\label{orthogonality'}
\lim_{n\to +\infty} \int \lf(1-\varphi\lf(\frac{|x|}{\mu_n}\rg)\rg)\nabla_{t,x}w_n(\theta_n,x)\cdot\frac{1}{\lambda_n^{N/2}}\nabla_{t,x}U\lf(\frac{\theta_n-t_n}{\lambda_n},\frac{x}{\lambda_n}\rg)dx=0.
\end{equation} 
\end{claim}
\begin{proof}
We start to show \eqref{orthogonality} when $\varphi=1$. By conservation of the energy for solutions of \eqref{lin_wave},
$$  \int \nabla_{t,x}w_n(\theta_n,x)\frac{1}{\lambda_n^{N/2}}\nabla_{t,x}U\lf(\frac{\theta_n-t_n}{\lambda_n},\frac{x}{\lambda_n}\rg)dx=\int \nabla_{t,x}w_n(t_n,x)\frac{1}{\lambda_n^{N/2}}\nabla_{t,x}U\lf(0,\frac{x}{\lambda_n}\rg)dx.$$
By the change of variable $\lambda_n y=x$, we see that \eqref{weakCVwn} implies \eqref{orthogonality} for $\varphi=1$. 

We next consider the case when $\varphi\in C^0(\RR^N)$ is compactly supported and satisfies $\varphi=1$ around $0$. Because of the case $\varphi=1$, one of the estimates \eqref{orthogonality} or \eqref{orthogonality'} implies the other. By the change of variable $\mu_n y=x$,
\begin{multline*}
\int \varphi\lf(\frac{|x|}{\mu_n}\rg)\nabla_{t,x}w_n(\theta_n,x)\cdot\frac{1}{\lambda_n^{N/2}}\nabla_{t,x}U\lf(\frac{\theta_n-t_n}{\lambda_n},\frac{x}{\lambda_n}\rg)dx\\
=\int \varphi\lf(|y|\rg)\mu_n^{N/2}\nabla_{t,x}w_n(\mu_n\tilde{\theta}_n,\mu_n y)\cdot\frac{1}{\tilde{\lambda}_n^{N/2}}\nabla_{t,x}U\lf(\frac{\tilde{\theta}_n-\tilde{t}_n}{\tilde{\lambda}_n},\frac{y}{\tilde{\lambda}_n}\rg)dy,
\end{multline*}
where $\tilde{\theta}_n=\frac{\theta_n}{\mu_n}$, $\tilde{\lambda}_n=\frac{\lambda_n}{\mu_n}$, $\tilde{t}_n=\frac{t_n}{\mu_n}$.
Replacing $w_n$ by the solution $(t,y)\mapsto \mu_n^{N/2}w_n(\mu_nt,\mu_n y)$ of \eqref{lin_wave}, $\theta_n$ by $\tilde{\theta}_n$, $t_n$ by $\tilde{t}_n$ and $\lambda_n$ by $\tilde{\lambda}_n$, we will assume in the sequel, in addition to \eqref{weakCVwn},
that $\mu_n=1$ for all $n$.

Extracting subsequences, we distinguish two cases.

\EMPH{Case $1$} Assume
\begin{equation}
\label{Case_1}
\lim_{n\to +\infty} \frac{\theta_n-t_n}{\lambda_n}=\pm \infty.
\end{equation}

Then, by Lemma \ref{L:lin_odd}, the energy of $\frac{1}{\lambda_n^{N/2}}U\lf(\frac{\theta_n-t_n}{\lambda_n},\frac{x}{\lambda_n}\rg)$ concentrates in sets of the form
$$ \Big\{|\theta_n-t_n|-C\lambda_n\leq |x|\leq |\theta_n-t_n|+C\lambda_n\Big\}$$
Recalling that $\mu_n=1$, we deduce that if $|\theta_n-t_n|\to +\infty$, \eqref{orthogonality} holds, and if $|\theta_n-t_n|\to 0$, \eqref{orthogonality'} holds. In both cases, the proof is complete. 

We next assume, after extraction, that 
\begin{equation}
\label{finite_limit}
\lim_{n\to +\infty} \theta_n-t_n=T\in \RR^*.
\end{equation} 

Let $\eps>0$, and $R$ (given by Lemma \ref{L:lin_odd}) such that 
$$ \limsup_{n\to \infty}\int_{\complement \CCC_n(R)}\frac{1}{\lambda_n^{N}}\left|\nabla_{t,x}U\lf(\frac{\theta_n-t_n}{\lambda_n},\frac{x}{\lambda_n}\rg)\rg|^2dx\leq \eps^2,$$
where $\CCC_n(R)=\left\{x\in \RR^N,\text{ s.t } |\theta_n-t_n|-R\lambda_n\leq |x|\leq |\theta_n-t_n|+R\lambda_n\right\}$ and $\complement\CCC_n(R)$ is its complement in $\RR^N$. Using the boundedness of $\nabla_{t,x}w_n$ in $(L^2)^{N+1}$, we get for large $n$
\begin{multline*}
\lf|\int \lf(\varphi\lf(x\rg)-\varphi(|\theta_n-t_n|)\rg)\nabla_{t,x}w_n(\theta_n,x)\frac{1}{\lambda_n^{N/2}}\nabla_{t,x}U\lf(\frac{\theta_n-t_n}{\lambda_n},\frac{x}{\lambda_n}\rg)dx\rg|\leq \\
C \max_{x\in\CCC_n(R)}\Big\{\big|\varphi\lf(x\rg)-\varphi(|\theta_n-t_n|)\big|\Big\}  +C\eps,
\end{multline*} 
where the constant $C$ depends only on the energy of $U$ and the bound of $\nabla_{t,x}w_n$ in $(L^2)^{N+1}$. As $\varphi$ is uniformly continuous, and $\lambda_n\to 0$, we get by \eqref{Case_1} and \eqref{finite_limit}
$$ \limsup_{n\to +\infty}\lf|\int \lf(\varphi\lf(x\rg)-\varphi(|\theta_n-t_n|)\rg)\nabla_{t,x}w_n(\theta_n,x)\frac{1}{\lambda_n^{N/2}}\nabla_{t,x}U\lf(\frac{\theta_n-t_n}{\lambda_n},\frac{x}{\lambda_n}\rg)dx\rg|\leq C\eps,$$
and hence (using the case $n=1$),
$$ \limsup_{n\to +\infty}\lf|\int \varphi\lf(x\rg)\nabla_{t,x}w_n(\theta_n,x)\frac{1}{\lambda_n^{N/2}}\nabla_{t,x}U\lf(\frac{\theta_n-t_n}{\lambda_n},\frac{x}{\lambda_n}\rg)dx\rg|\leq C\eps.$$
The proof is complete if \eqref{Case_1} holds. 

\EMPH{Case $2$} Assume
\begin{equation}
\label{Case_2}
\lim_{n\to +\infty} \frac{\theta_n-t_n}{\lambda_n}=t_0\in \RR.
\end{equation}
Then by Lemma \ref{L:lin_odd} the $L^2$ norm of  $\frac{1}{\lambda_n^{N/2}}\nabla_{t,x}U\lf(\frac{\theta_n-t_n}{\lambda_n},\frac{x}{\lambda_n}\rg)$ is localized in sets of the form
$$ \Big\{ C^{-1}\lambda_n \leq |x|\leq C\lambda_n\Big\}.$$
If $\lambda_n\to +\infty$ or $\lambda_n\to 0$, the argument of Case $1$ yields \eqref{orthogonality} and \eqref{orthogonality'}. Let us assume
$$ \lim_{n\to +\infty}\lambda_n=\lambda_{\infty}\in (0,+\infty).$$
Then
$$\frac{1}{\lambda_n^{N/2}}\nabla_{t,x}U\lf(\frac{\theta_n-t_n}{\lambda_n},\frac{x}{\lambda_n}\rg)=\frac{1}{\lambda_{\infty}^{N/2}}\nabla_{t,x}U\lf(t_0,\frac{x}{\lambda_{\infty}}\rg)+o_n(1)\text{ in }\lf(L^2\lf(\RR^N\rg)\rg)^{N+1}.$$
Thus we must show
\begin{equation}
\label{orthogonality_bis}
\lim_{n\to +\infty} \int \varphi\lf(x\rg)\nabla_{t,x}w_n(\theta_n,x)\frac{1}{\lambda_{\infty}^{N/2}}\nabla_{t,x}U\lf(t_0,\frac{x}{\lambda_{\infty}}\rg)dx=0.
\end{equation} 
First notice that if $\Phi\in (L^2)^{N+1}$,
\begin{align*}
\int \nabla_{t,x}w_n(t_n,x) \cdot \Phi(x)dx&=\int \lambda_n^{N/2}\nabla_{t,x}w_n(t_n,\lambda_n y)\cdot\lambda_n^{N/2}\Phi\lf(\lambda_ny\rg)dy\\
&=\int \lambda_n^{N/2}\nabla_{t,x}w_n(t_n,\lambda_n y)\cdot\lambda_{\infty}^{N/2}\Phi\lf(\lambda_{\infty}y\rg)dy+o_n(1).
\end{align*}
At the last line we used that $\lambda_n^{N/2}\Phi\lf(\lambda_ny\rg)$ converges strongly to $\lambda_{\infty}^{N/2}\Phi\lf(\lambda_{\infty}y\rg)$ in $(L^2)^{N+1}$.
Thus by \eqref{weakCVwn}
\begin{equation}
\label{weakCVwnter}
\nabla_{t,x}w_n(t_n,x)\xrightharpoonup[n\to \infty]{}0 \text{ in }(L^2)^{N+1}.
\end{equation} 
Next, consider the solution $v$ of \eqref{lin_wave} with initial data $(v_0,v_1)\in \hdot\times L^2$ such that
$$\Delta v_0(x)=\frac{1}{\lambda_{\infty}^{N/2}}\Div\left(\varphi(x)\nabla_x U\lf(t_0,\frac{x}{\lambda_{\infty}}\rg)\rg),\quad v_1(x)=\frac{1}{\lambda_{\infty}^{N/2}}\varphi(x)\partial_t U\lf(t_0,\frac{x}{\lambda_{\infty}}\rg)$$
Write $\theta_n=\lambda_{\infty}t_0+t_n+\eps_n$, with $\eps_n\to 0^+$. Then by conservation of the energy
\begin{multline*}
\int \varphi\lf(x\rg)\nabla_{t,x}w_n(\theta_n,x)\frac{1}{\lambda_{\infty}^{N/2}}\nabla_{t,x}U\lf(t_0,\frac{x}{\lambda_{\infty}}\rg)dx=\int \nabla_{t,x}w_n(\theta_n,x)\nabla_{t,x}v(0,x)dx\\
=\int \nabla_{t,x}w_n(t_n,x)\nabla_{t,x}v(-\lambda_{\infty}t_0-\eps_n,x)dx=\int \nabla_{t,x}w_n(t_n,x)\nabla_{t,x}v(-\lambda_{\infty}t_0,x)dx+o_n(1),
\end{multline*} 
which shows \eqref{orthogonality_bis} in view of \eqref{weakCVwnter}.
\end{proof}

We next prove Claim \ref{C:dispersive_tw}. 
\begin{proof}
We prove the result when $N$ is odd, although it should also hold when $N$ is even.
Rescaling if necessary, we will assume
\begin{equation}
\label{rescaling}
 \forall n,\quad \tlambda_n=1.
\end{equation} 
Note that the assumption \eqref{dispersive_w} implies that for any sequence $\{\lambda_n\}$, $\{t_n\}$, 
\begin{equation}
 \label{weak_CV_rescaled}
\left( \frac{1}{\lambda_n^{\frac{N-2}{2}}}w_n\left(\frac{-t_n}{\lambda_n},\frac{\cdot}{\lambda_n}\right), \frac{1}{\lambda_n^{\frac{N}{2}}}\partial_t w_n\left(\frac{-t_n}{\lambda_n},\frac{\cdot}{\lambda_n}\right)\right)\xrightharpoonup[n\to \infty]{}
(0,0)\text{ weakly in }\hdot\times L^2.
 \end{equation} 
Indeed if \eqref{weak_CV_rescaled} does not hold, the sequence $\{w_n\}$ would have a nontrivial profile decomposition, contradicting \eqref{dispersive_w}.

Conversely, we claim that \eqref{dispersive_tw} holds as soon as for all sequences $\{\lambda_n\}$, $\{t_n\}$,
\begin{equation}
\label{weak_CV_twn}
\left( \frac{1}{\lambda_n^{\frac{N-2}{2}}}\tw_n\left(\frac{-t_n}{\lambda_n},\frac{\cdot}{\lambda_n}\right), \frac{1}{\lambda_n^{\frac{N}{2}}}\partial_t \tw_n\left(\frac{-t_n}{\lambda_n},\frac{\cdot}{\lambda_n}\right)\right)\xrightharpoonup[n\to \infty]{}
(0,0)\text{ weakly in }\hdot\times L^2.
\end{equation} 
Again, if \eqref{dispersive_tw} does not hold, then the sequence $(\tw_{n}(0),\partial_t\tw_n(0))$ has a profile decomposition with at least one nonzero profile, which contradicts \eqref{weak_CV_twn}.

Let us show \eqref{weak_CV_twn}. Let $(Z_0,V_1)\in \dot{H}^{-1}\times L^2$ and $V_0\in \hdot$ such that 
$\Delta V_0=Z_0$. Let $V$ be the solution of \eqref{lin_wave} with initial conditions $(V_0,V_1)$. We have
\begin{multline}
\label{test_twn}
\int \frac{1}{\lambda_n^{\frac{N-2}{2}}}\tw_n\left(\frac{-t_n}{\lambda_n},\frac{x}{\lambda_n}\right)Z_0(x)dx+\int \frac{1}{\lambda_n^{\frac{N}{2}}}\partial_t\tw_n\left(\frac{-t_n}{\lambda_n},\frac{x}{\lambda_n}\right)V_1(x)dx\\
=\int \nabla_{x}\tw_n\left(0,x\right)\cdot \lambda_n^{\frac{N}{2}}\nabla_x V(t_n,\lambda_nx)dx+\int \partial_t\tw_n\left(0,x\right)\lambda_n^{\frac{N}{2}}\partial_tV(t_n,\lambda_nx)\\
=\int \nabla_{x}\left(\varphi(|x|)w_{0,n}\left(x\right)\right)\cdot \lambda_n^{\frac{N}{2}}\nabla_x V(t_n,\lambda_nx)dx+\int \varphi(|x|)w_{1,n}\left(x\right)\lambda_n^{\frac{N}{2}}\partial_tV(t_n,\lambda_nx).
\end{multline}
Thus it suffices to show
\begin{gather}
 \label{CV0_1}
\lim_{n\to \infty}\int \varphi(|x|)\nabla_{x}w_{0,n}\left(x\right)\cdot \lambda_n^{\frac{N}{2}}\nabla_x V(t_n,\lambda_nx)dx+\int \varphi(|x|)w_{1,n}\left(x\right)\lambda_n^{\frac{N}{2}}\partial_tV(t_n,\lambda_nx)=0\\
\label{CV0_2}
\lim_{n\to\infty}\int \left(\nabla_{x}\varphi(|x|)\right)w_{0,n}\left(x\right)\cdot \lambda_n^{\frac{N}{2}}\nabla_x V(t_n,\lambda_nx)dx=0.
\end{gather}
The first limit, \eqref{CV0_1}, follows immediately from Claim \ref{C:ortho}. To show \eqref{CV0_2}, we use that there exists $C>0$ such that $\nabla \varphi$ is supported in $\{1/C\leq |x|\leq C\}$, and distinguish several cases.

If $t_n$ is bounded, then one can assume after extraction that $t_n$ has a limit $T\in[0,+\infty)$. If $\lambda_n\to 0$ or $\lambda_n\to +\infty$ then by Lemma \ref{L:lin_odd}, 
\begin{equation}
\label{part_a_linfini}
\lim_{n\to +\infty} \int_{1/C\leq |x|\leq C} \lambda_n^N\left|\nabla V(t_n,\lambda_nx)\right|^2dx=0,
\end{equation} 
and \eqref{CV0_2} follows. If $\lambda_n$ has a limit $\lambda_{\infty}\in (0,+\infty)$, then $\lambda_n^{\frac{N}{2}}\nabla V(t_n,\lambda_nx)$ converges strongly to $\lambda_{\infty}^{\frac{N}{2}}\nabla V(T,\lambda_{\infty}x)$, and we are reduced to show
\begin{equation*}
\lim_{n\to\infty}\int_{1/C\leq |x|\leq C} \left(\nabla_{x}\varphi(|x|)\right)w_{0,n}\left(x\right)\cdot \lambda_{\infty}^{\frac{N}{2}}\nabla V(T,\lambda_{\infty}x)dx=0,
\end{equation*}
which follows from the fact that by \eqref{weak_CV_rescaled}, $w_{0,n}$ tends to $0$ weakly in $\hdot$ (and thus, by Hardy's inequality, that $\frac{1}{|x|}w_{0,n}$ tends to $0$ weakly in $L^2$).

We next treat the case when $t_n$ is not bounded. Extracting, we will assume that $t_n\to +\infty$ (the case $t_n\to -\infty$ is analoguous). If $t_n/\lambda_n\to 0$ or $t_n/\lambda_n\to +\infty$, Lemma \ref{L:lin_odd} implies again \eqref{part_a_linfini}, and \eqref{CV0_2} follows. It remains to consider the case when (after extraction), $t_n/\lambda_n\to \ell\in (0,+\infty)$. By Lemma \ref{L:lin_odd}, for all $\eps>0$ there exists $R_{\eps}$ such that for all $R\geq R_{\eps}$,
\begin{equation*}
 \limsup_{n\to \infty}\int_{\left||x|-\ell\right|\geq \frac{R}{\lambda_n}} \lambda_n^N\left|\nabla V(t_n,\lambda_nx)\right|^2dx\leq \eps
\end{equation*}
As a consequence,
\begin{equation}
\label{localized_profile}
 \limsup_{n\to \infty}\int_{\left||x|-\ell\right|\geq \frac{1}{\sqrt{\lambda_n}}} \lambda_n^N\left|\nabla V(t_n,\lambda_nx)\right|^2dx=0.
\end{equation}
It remains to show that
\begin{equation}
\label{lim0}
 \lim_{n\to \infty}\left|\int_{\left||x|-\ell\right|\leq \frac{1}{\sqrt{\lambda_n}}} \partial_r\varphi(|x|)w_{0,n}\left(x\right)\lambda_n^{\frac{N}{2}}\partial_r V(t_n,\lambda_nx)dx\right|=0.
\end{equation}
We have
\begin{multline}
\label{last_equality}
 \int_{\left||x|-\ell\right|\leq \frac{1}{\sqrt{\lambda_n}}} \partial_r\varphi(|x|)w_{0,n}\left(x\right)\lambda_n^{\frac{N}{2}}\partial_r V(t_n,\lambda_nx)dx
\\
=\int_{\left||x|-\ell\right|\leq \frac{1}{\sqrt{\lambda_n}}} |x|\partial_r\varphi(|x|)\frac{1}{|x|}w_{0,n}\left(x\right)\lambda_n^{\frac{N}{2}}\partial_r V(t_n,\lambda_nx)dx\\
=\int_{\left||x|-\ell\right|\leq \frac{1}{\sqrt{\lambda_n}}} \ell\partial_r\varphi(\ell)\frac{1}{|x|}w_{0,n}\left(x\right)\lambda_n^{\frac{N}{2}}\partial_r V(t_n,\lambda_nx)dx+o_n(1)\\
=\ell\partial_r\varphi(\ell)\int_{\RR^N} \frac{1}{|x|}w_{0,n}\left(x\right)\lambda_n^{\frac{N}{2}}\partial_r V(t_n,\lambda_nx)dx+o_n(1).
\end{multline}
At the third line, we have used that $r\partial_r\varphi$ is continuous and thus
$$ \lim_{n\to \infty}\sup_{|r-\ell|\leq\frac{1}{\sqrt{\lambda_n}}}\left|r  \partial_r \varphi(r)-\ell\partial_{r}\varphi(\ell)\right|=0.$$
At the last line we have used \eqref{localized_profile}.
By Hardy's inequality and assumption \eqref{weak_CV_rescaled}, $\frac{1}{|x|}w_{0,n}$ converges weakly to $0$ in $L^2$, and thus \eqref{last_equality} implies \eqref{lim0}, which concludes the proof of Claim \ref{C:dispersive_tw}.
\end{proof}

\section{Family of sequences of positive numbers}
\label{A:faraway}
\begin{claim}
\label{C:faraway}
Let $\{\lambda_n\}_n$, $\{\nu_n\}_n$ and for $j\in\NN$, $\{\rho_{j,n}\}_n$, be sequences of positive numbers and assume 
\begin{equation}
\label{lambda_nu}
\lambda_n\ll \nu_n.
\end{equation}
Then, after extraction of subsequences in $n$, there exists a sequence $\{\mu_n\}_n$ such that 
\begin{gather}
\label{lambda_mu_nu}
\lambda_n\ll \mu_n\ll \nu_n\\
\label{mu_rhoj}
\forall k,\quad \mu_n\ll \rho_{k,n}\text{ or }\rho_{k,n}\ll\mu_n.
\end{gather}
\end{claim}
\begin{proof}
Let for $s\in(0,1)$,
$$ \mu_n(s)=\lambda_n^{1-s}\nu_n^{s}.$$
Note that for any $s\in (0,1)$, $\lambda_n\ll\mu_n(s)\ll \nu_n$. Let $j\in \NN$. Then, extracting subsequences in $n$ if necessary, we are in one of the three following cases:
\begin{equation}
\label{cases_s}
\begin{cases}
\forall s\in (0,1)\quad \mu_n(s)\ll \rho_{j,n}\text{ or }\\
\forall s\in (0,1)\quad \rho_{j,n}\ll \mu_n(s)\text{ or }\\
\exists s_j\in (0,1),\quad \forall s\in (0,s_j),\; \mu_n(s)\ll \rho_{j,n}\text{ and } \forall s\in(s_j,1),\;  \rho_{j,n}\ll \mu_n(s).
\end{cases}
\end{equation}
Indeed let
$$ s_j=\inf\Big\{s\in [0,1]\;\big|\; \lf\{\rho_{j,n}/\mu_n(s)\rg\}_n\text{ is bounded}.\Big\}.$$
Note that $\mu_n(s)\ll \mu_n(s')$ is $s<s'$. As a consequence, if $s_j=0$, then $\rho_{j,n}/\mu_n(s)\to 0$ for all $s\in(0,1)$. 
Similarly if $s_j=1$, then $\lf\{\rho_{j,n}/\mu_n(s)\rg\}_n$ is never bounded for $s\in(0,1)$ and by diagonal extraction we can find a subsequence such that $\rho_{j,n}/\mu_n(s)\to +\infty$ for any $s\in (0,1)$. Finally if $s_j\in(0,1)$,
then $\rho_{j,n}/\mu_n(s)\to 0$ for all $s\in(s_j,1)$, and $\lf\{\rho_{j,n}/\mu_n(s)\rg\}_n$ is not bounded for $s\in (0,s_j)$. Using diagonal extraction again we can assume that $\rho_{j,n}/\mu_n(s)\to +\infty$ for all $s\in (0,s_j)$. Hence \eqref{cases_s}.

After another diagonal extraction, we can assume that \eqref{cases_s} holds for all $j\in \NN$. Chosing $s\in(0,1)$ distinct from all $s_j$, and letting $\mu_n=\mu_n(s)$ we get the desired properties \eqref{lambda_mu_nu} and \eqref{mu_rhoj}.
\end{proof}
\bibliographystyle{alpha} 
\bibliography{blowup}

\end{document}